\documentclass[11pt]{article}

\usepackage[margin=1in]{geometry}
\usepackage{amsmath,amssymb,amsfonts,amsthm,mathtools,mathrsfs}
\usepackage{booktabs,enumitem,float,array}
\usepackage[round,authoryear]{natbib}
\usepackage[colorlinks=true,linkcolor=blue,citecolor=blue,urlcolor=blue]{hyperref}

\theoremstyle{plain}
\newtheorem{theorem}{Theorem}[section]
\newtheorem{lemma}[theorem]{Lemma}
\newtheorem{proposition}[theorem]{Proposition}
\newtheorem{corollary}[theorem]{Corollary}
\newtheorem{assumption}[theorem]{Assumption}
\newtheorem{definition}[theorem]{Definition}
\theoremstyle{definition}
\newtheorem{remark}[theorem]{Remark}
\newtheorem{example}[theorem]{Example}

\newcommand{\E}{\mathbb{E}}
\newcommand{\Var}{\operatorname{Var}}
\newcommand{\Cov}{\operatorname{Cov}}

\newcommand{\diag}{\operatorname{diag}}
\newcommand{\F}{\mathcal F}
\newcommand{\hVar}{\widehat{\operatorname{Var}}}
\providecommand{\cE}{\mathcal E}
\providecommand{\cR}{\mathcal R}
\providecommand{\cU}{\mathcal U}
\providecommand{\cD}{\mathcal D}

\providecommand{\Hc}{\mathcal H}
\providecommand{\one}{\mathbf 1}
\providecommand{\avg}{\operatorname{avg}}
\newcommand{\dmin}{\underline{\mu}}
\newcommand{\Eg}[1]{\E_{g}[#1]}
\newcommand{\Labs}{{\Lambda^{|A|}}}
\newcommand{\gn}{\mathsf e}
\newcommand{\nb}[1]{\langle #1\rangle}
\newcommand{\echo}{{\mathsf d^{A}}}
\newcommand{\src}{\mathsf a}
\newcommand{\Ebar}{\bar{\cE}}
\newcommand{\Odeg}{\Omega_{\deg}}
\newcommand{\Osp}{\Omega_{\mathrm{sp}}}
\newcommand{\Osup}{\Omega_{\mathrm{sup}}}
\newcommand{\Ods}{\Omega_{\mathrm{ds}}}
\newcommand{\Ohg}{\Omega_{\hat g}}
\newcommand{\Ohgone}{\Omega_{\hat g,1}}
\newcommand{\Odegone}{\Omega_{\deg,1}}
\newcommand{\OB}{\Omega_{\mathrm{B}}}
\newcommand{\OU}{\Omega_{\mathrm{U}}}
\newcommand{\oneU}{\one_{\mathrm{U}}}
\newcommand{\oneds}{\one_{\mathrm{ds}}}
\newcommand{\onedeg}{\one_{\deg}}

\title{Central limit theorems and network-free inference for heterogeneous spatial autoregressions on sparse graphons}
\author{Wei Li\thanks{Department of Statistics, University of Wisconsin--Madison. Email: \texttt{wli734@wisc.edu}.}
\and
Robert Lunde\thanks{Department of Statistics and Data Science, Washington University in St.\ Louis. Email: \texttt{lunde@wustl.edu}.}}
\date{}

\begin{document}
\maketitle

\begin{abstract}
Network dependence is widely regarded as a fundamental obstacle to inference: outcomes linked through a graph are not independent, the graph is rarely observed in full, and models for how dependence propagates are easily misspecified. We show that for a broad class of responses generated through a random network, the obstacle is smaller than it appears. We study a heterogeneous spatial autoregression on a sparse graphon, in which each node's response is a weighted average of its neighbors' responses with a node-specific, nonparametric coefficient plus a node-specific input, and we prove that sample means and smooth functions of vector means are asymptotically Gaussian down to average degree nearly logarithmic in the population size, allowing unbounded graphons with hubs and heavy-tailed inputs. The limiting variance interpolates exactly between two regimes according to the sampling fraction. When the observed nodes are a vanishing fraction of the population, the variance is the marginal variance of the response, the ordinary sample variance is consistent, and the usual $t$-interval is valid with no knowledge of the network. The proof averages over the graph rather than conditioning on it, isolating the behavior of the response through an exact identity for resolvents of Markov kernels and controlling the network noise by an Efron--Stein inequality; neither step requires the walk-by-walk covariance bounds of fixed-network analyses. Consequences for regression with node-level and network-derived covariates are derived.
\end{abstract}

\section{Introduction}

Network dependence is widely regarded as a serious obstacle to statistical inference. When outcomes are linked through a graph, the observations are not independent  and confidence intervals built as if they were can have coverage far from nominal \citep{Lee03072021}.
In applications, the situation looks worse still: the network is rarely observed in full, the mechanism by which dependence propagates is unknown, and any simple generative model is likely misspecified. A reasonable reader of this literature might conclude that valid inference for even a population mean is out of reach without strong, unverifiable assumptions about the graph.

There is a tension between this pessimism and representation theorems for network data. In \citet{li2025assumption} we showed that for a broad class of network-linked data satisfying a joint exchangeability condition, a representation in the spirit of \citet{aldous-representation-array} and \citet{hoover-exchangeability} reduces the problem to manageable building blocks: independent latent random variables and a graphon. The catch is that the representation requires infinite exchangeability. Outcomes that are generated through the network, by peer effects, contagion, or spillovers, form a finitely exchangeable array for each population size, but not the finite-dimensional margins of a single infinitely exchangeable sequence: adding a node changes everyone's response. The question this paper addresses is whether the optimism survives once infinite joint exchangeability is dropped.

Our answer is that it largely does. We study responses generated by a heterogeneous
spatial autoregression on a sparse graphon, chosen deliberately as a demanding test case:
every response depends on every other through feedback along the network, with no decay in
network distance and random degrees in every normalization, average degree may be only just
above $\log N$, and both the graphon and the inputs may be heavy-tailed. Despite all of
this, the dependence is weak once the network is averaged over: each response is, up to a
vanishing error, a function of its own latent position, and every pair of responses is
correlated only at order $1/N$, with a leading constant identified through an explicit
influence function. This correlation is invisible within a sample that is a vanishing
fraction of the population, so the ordinary sample variance is then consistent for the
variance of the mean, and the usual $t$-interval is valid with no knowledge of the network,
the graphon, or the propagation mechanism. This is in the spirit of the finite de Finetti
theorems of \citet{diaconis1980finite}: for general state spaces, the first $n$ coordinates
of a finitely exchangeable, extendable vector of length $N$ are close in total variation to
a mixture of i.i.d.\ laws provided $n^2/N\to0$, and this rate is sharp. Our requirement is
only $n/N\to0$, since we ask that a single statistic, the sample mean, behave as under
i.i.d.\ sampling, rather than that the joint law of the subsample do so. That network-free
inference holds even in this model suggests that the phenomenon is not an artifact of local
or limited propagation.

One key to our result is that the network is treated as random rather than fixed.
Nearly all prior central limit theory for network-dependent data conditions on the
realized graph and imposes weak-dependence conditions expressed in terms of it (see
Section~\ref{sec:related} for a literature review), so that the network must be known
for the conditions to be checked and, typically, for the variance to be estimated. We
regard averaging over the graph as the natural choice for the questions studied here.
Assumption-lean inference treats covariates as random \citep{li2025assumption}, but in network data the same latent attributes typically drive both the graph and the covariates, so holding the graph fixed while letting the covariates vary is internally inconsistent.
When both the attributes and the graph are random, both channels of network
dependence contribute to fluctuations: homophily, through the dependence of the graph
on the attributes, and contagion, through the dependence of the responses on the
graph. The two are generally difficult to distinguish in observational data
\citep{shalizi2011homophily}, and our results do not require distinguishing them.
A fixed-graph analysis remains appropriate when the realized network is itself the
object of interest; averaging targets a quantity that is stable across realizations
of the network-forming process. The phenomenon we describe emerges only under this
averaging.

The technical content of the paper is to show that this averaging can be carried out
in a sparse regime. We split the response into a signal component, driven by the
latent positions, and a graph component, driven by the edges given the positions. We
handle the signal component by viewing the response as a functional of the empirical
distribution of latent positions and isolating its linear part through an exact
identity. We handle the graph component with a leave-one-node-out identity and an
Efron--Stein inequality over vertex blocks, which replaces the combinatorial control
of walk covariances used in the fixed-network literature and brings the sparsity
requirement down to an average degree just above $\log N$. Since the argument
controls the effect of a single node's edges on a linear functional of a
degree-normalized resolvent, a form shared by other statistics such as random-walk
centralities, we expect these tools to apply beyond the present problem, with sparsity
requirements that depend on the statistic.

The paper is organized as follows. Section~\ref{sec:model} states the model and the sense
in which it is general. Section~\ref{sec:main} states the main theorem and the variance
formula, explains its two-stage structure and its relation to finite-population sampling,
and derives the consequences for network-free inference, together with the delta-method
extension to smooth functionals. Section~\ref{sec:proof-outline} outlines the proof:
Section~\ref{sec:operators} treats the component of the response driven by the latent
positions, through an exact identity for degree-normalized operators, and
Section~\ref{sec:fluctuation} treats the network fluctuation. Section~\ref{sec:simulations}
gives simulations, and Section~\ref{sec:conclusion} discusses generalizations. All proofs
are in the appendix.

\subsection{Model}\label{sec:model}

Let $\xi_1,\dots,\xi_N$ be i.i.d.\ $\mathrm{Unif}(0,1)$ latent positions, which we also call \emph{types}, and let $A=A^{(N)}\in\{0,1\}^{N\times N}$ be a symmetric adjacency matrix with $A_{ii}=0$ and degrees $d_i:=\sum_{j\ne i}A_{ij}$.  For the network, we consider a triangular array of sparse graphon models; variants of the sparse graphon model have been proposed in the mathematics literature by \citet{bollobas2011sparse} and  \citet{borgs-lp-part-one}, and in the statistics literature by \citet{Bickel-Chen-on-modularity}.

\begin{definition}[Sparse graphon]\label{def:graphon}
Let $\{U_{ij}\}_{i<j}$ be i.i.d.\ $\mathrm{Unif}(0,1)$, independent of $\{\xi_i\}$, and let $w_N:[0,1]^2\to[0,\infty)$ be symmetric with $\int_0^1\int_0^1 w_N(u,v) \ du \ dv=1$, $w_N(x,x) = 0$ for any $x \in [0,1]$, and $\rho_N\in(0,1)$ with $\rho_N\to0$. For $1 \leq i<j \leq N$\footnote{By assumption (M3) introduced later, $\rho_N w_N \leq 1/2$, so $w_N$ need not be truncated.},
\begin{equation}\label{eq:graphon}
A_{ij}= \one\{U_{ij}\le\rho_Nw_N(\xi_i,\xi_j)\}.
\end{equation}
\end{definition}
The convention $w_N(x,x)=0$ changes $w_N$ only on a null set, so it does not affect any integral. The results of \citet{borgs-lp-part-one} suggest that sparse graphons are the natural generalization of graphons \citep{lovasz2012large} to sparse graph sequences, in the sense that they are the limit objects of appropriately normalized sparse sequences converging in cut metric. In the statistics literature, it is common to fix a single graphon $w$ and target fixed functionals of it, such as subgraph frequencies \citep{Bickel-Chen-Levina-method-of-moments} or eigenvalues of integral operators \citep{lunde2023subsampling}. This convention is partly justified by this convergence theory, but it should be noted that these results are qualitative: they identify the limit of a sparse sequence converging in cut metric, but give no rate, and inference requires more than the existence of a limit. Moreover, this presumes that the sequence of networks is converging: that a network observed at different sizes has the same connection profile up to normalization. In many applications, this is not the right description. A school of two hundred students and one of two thousand are generally not two samples from one limiting object up to sparsity.  Allowing $w_N$ to depend on $N$ covers non-convergent sequences and growing hubs, and reduces to the fixed-graphon setting when $w_N\equiv w$.

In the present problem the triangular-array formulation is in any case forced on us rather than chosen. Even with $w_N\equiv w$, the responses in \eqref{eq:SAR} form a triangular array (Remark~\ref{rem:exch}), since the law of a node's response depends on the size of the population it is embedded in. The limit theory must therefore be stated row-wise, with centering $\mu_N$ and variance $\sigma^2_{n,N}$ depending on $N$, and once this is accepted, letting the graphon vary with $N$ costs nothing further. All constants in our conditions are uniform over the array, and that is the only sense in which uniformity across $N$ enters.

We now introduce the regression model, which is a generalization of spatial autoregressive models commonly studied in the econometrics literature (see, for example, \citet{manski1993identification} and \citet{bramoulle2009identification}), and is related to various models of opinion formation, which we will discuss shortly. In what follows, the response $Y_i^{(N)}$, the input
function $z_N$, and the coefficient function $\gamma_N$, may all vary with $N$; all
conditions below are uniform in $N$, for example $\sup_x|\gamma_N(x)|\le\bar\gamma<1$
for every $N$. Since all results concern the distribution at each fixed $N$, the latent positions may equally be taken as a triangular array, i.i.d.\ within each row. To lighten notation, we suppress the subscript $N$ when no confusion
arises, writing $w$, $z$ and $\gamma$.

\begin{definition}[Heterogeneous spatial autoregression]\label{def:SAR}
Let $z,\gamma:[0,1]\to\mathbb R$ be measurable with $\bar\gamma:=\sup_{x\in[0,1]}|\gamma(x)|<1$, and set $Z_i:=z(\xi_i)$, $\gamma_i:=\gamma(\xi_i)$. The response vector $Y=(Y_1,\dots,Y_N)^\top$ solves
\begin{equation}\label{eq:SAR}
Y_i=\gamma_i\,\frac{1}{d_i}\sum_{j\ne i}A_{ij}Y_j+Z_i,\qquad i=1,\dots,N,
\end{equation}
where $A$ is generated by the sparse graphon model \eqref{eq:graphon}, with the convention that the neighborhood average is zero when $d_i=0$. Writing $\Gamma=\mathrm{diag}(\gamma_1,\dots,\gamma_N)$, $Z=(Z_1,\dots,Z_N)^\top$, and $P=D^{-1}A$ with $D=\mathrm{diag}(d_1,\dots,d_N)$ and the $i$th row of $P$ set to zero when $d_i=0$, \eqref{eq:SAR} reads $Y=\Gamma PY+Z$, so that
\begin{equation}\label{eq:solution}
Y=(I-\Gamma P)^{-1}Z .
\end{equation}
The inverse exists for every realization of the network: $P$ is nonnegative with row sums at most one, so $\|\Gamma P\|_{\infty\to\infty}\le\bar\gamma<1$.
\end{definition}

The same bound gives the Neumann expansion
\begin{equation}\label{eq:neumann}
Y=\sum_{k\ge0}(\Gamma P)^kZ,\qquad
Y_i=\sum_{k\ge0}\ \sum_{i=i_0,i_1,\dots,i_k}\ \prod_{s=0}^{k-1}\frac{\gamma_{i_s}A_{i_si_{s+1}}}{d_{i_s}}\;Z_{i_k},
\end{equation}
where the inner sum runs over walks of length $k$ starting at $i$. A node's response is the sum of the inputs of all nodes it can reach, discounted along the walk by the autoregressive coefficient and the reciprocal degree at each step; equivalently, when $\gamma\ge0$, consider a random walk from $i$ that at each node continues to a uniformly chosen neighbor with probability $\gamma$ at that node and stops otherwise; then $Y_i$ is the expected value of $Z_\tau/(1-\gamma_\tau)$ at the node $\tau$ where the walk stops. Dependence is therefore global: every node depends on every other node in its component, and there is no decay with network distance to exploit.

The model is more general than its additive form suggests, and a reparametrization makes this more transparent. Take $\gamma_i \geq 0$ and write $Z_i=(1-\gamma_i)X_i$ for a node-level variable $X_i=X(\xi_i)$. Then \eqref{eq:SAR} becomes
\begin{equation}\label{eq:convex}
Y_i=\gamma_i\,\frac1{d_i}\sum_{j\ne i}A_{ij}Y_j+(1-\gamma_i)X_i ,
\end{equation}
so each response is a convex combination of the neighborhood average and the node's own input. Read as an opinion or behavior, $X_i$ is node $i$'s intrinsic position and $\gamma_i$ is its susceptibility to peers. In matrix form, $Y=(I-\Gamma P)^{-1}(I-\Gamma)X$, which is the equilibrium of the Friedkin-Johnsen model of opinion dynamics with heterogeneous stubbornness \citep{friedkin1990social}, a widely studied generalization of DeGroot averaging.

Three further points about the specification. First, because $Z$ and $\gamma$ are arbitrary measurable functions of the latent position, the joint law of the intrinsic input, susceptibility, and connection propensity is unrestricted. Susceptibility to peers may be correlated with the attributes that govern who one connects to, so homophily in the network and heterogeneity in influence are modeled jointly rather than assumed independent. Second, the equation is linear in the neighbors' responses but has random intercept and slope with an arbitrary joint distribution, which places it in the random-coefficient tradition of \citet{hastie1993varying}. Third, under the monotonicity condition that the structural response is increasing in the latent index, \eqref{eq:SAR} admits a conditional-quantile interpretation in the sense of \citet{zhu2019network}, so that quantile-type network autoregressions are also covered. None of this structure needs to be known or estimated: the target is the mean response, and the heterogeneity is averaged out in the results that follow.

\begin{remark}[Generality of the input]\label{rem:generality}
The requirement that $Z_i$ and $\gamma_i$ be functions of a single uniform random
variable is weaker than it appears. A uniform random variable can be split into
countably many independent uniforms, $\xi\mapsto(\xi',\xi'',\xi''',\dots)$, by a
measurable bijection, and any random element of a standard Borel space is a measurable
function of a uniform \citep[see, e.g.,][]{kallenberg2002foundations}. Hence the specification
\[
X_i=h_X(\xi_i',\xi_i''),\qquad \gamma_i=h_\gamma(\xi_i',\xi_i'''),\qquad
\varepsilon_i=h_\varepsilon(\xi_i',\xi_i''''),\qquad Z_i=f(X_i,\varepsilon_i,\gamma_i),
\]
in which the covariates $X_i$ and the noise $\varepsilon_i$ may take values in any
standard Borel spaces, while $\gamma_i$ and $Z_i$ are real-valued with
$|\gamma_i|\le\bar\gamma$, is a special case of Definition~\ref{def:SAR}. This is the same mechanism by which
the Aldous-Hoover representation, in which node and edge variables are functions of
uniform latent variables, is more general than its form suggests.
\end{remark}

\begin{remark}[Finite but not infinite exchangeability]\label{rem:exch}
For each $N$, the latent positions, the graph mechanism and \eqref{eq:SAR} are invariant under relabeling, so $Y^{(N)}$ is exchangeable. The array $\{Y^{(N)}_i\}$ is not, however, projectively consistent across $N$: adding a node adds a neighbor to other nodes and changes their responses, so the law of $(Y^{(N)}_1,\dots,Y^{(N)}_n)$ depends on $N$. The responses form a triangular array of finitely exchangeable vectors rather than the finite-dimensional margins of an infinitely exchangeable sequence, and the Aldous--Hoover representation of \citet{li2025assumption} is not available. This is the sense in which the present setting is outside the i.i.d.-triplet framework, and it is why the results below must be stated row-wise, with centering and variance depending on $N$.
\end{remark}

\subsection{Related work}\label{sec:related}

Conditional on the graph, the phenomenon we describe does not exist, and this is the organizing fact of the comparison below: every framework we are aware of either conditions on the realized network or requires a form of exchangeability that the autoregressive mechanism breaks.

\paragraph*{Limit theory for network-dependent data}
Central limit theorems for network-dependent variables descend from dependency-graph and Stein's-method arguments \citep{baldi1989normal,chen2004normal} and from the weak-dependence tradition for random fields \citep{jenish2009central}. The modern network versions replace a Euclidean index set by a graph: \citet{kojevnikov2021limit} measure dependence by covariances of nonlinearly transformed variables that decay with network distance ($\psi$-dependence) and trade off the decay rate against the growth of neighborhoods; \citet{xu2024functional,jiang2025limit} generalize the functional (physical) dependence measure of \citet{wu2005nonlinear} to network-indexed variables and verify it for spatial autoregressions on a given weight matrix; \citet{leung2026normal} formulate weak dependence through exponential stabilization and obtain primitive conditions from branching-process arguments. In the causal setting, \citet{ogburn2024causal} prove a CLT on a single growing network that allows each unit to depend on a growing number of others, with dependence from both transmission along ties and latent similarity across ties. For least squares, \citet{lei2024leastsquares} establish a central limit theorem for the ordinary
least squares estimator under network dependence on a fixed network, with
neighborhood growth conditions, and propose a consistent variance estimator that uses
the observed network, with tuning chosen adaptively to the hypothesis being tested. All of these results are stated conditionally on an observed network. For the variables generated by \eqref{eq:SAR}, these conditions are not the natural
tool: once the graph is averaged over, the responses are exchangeable, so all pairwise
covariances are equal and no pair of nodes is more distant than another. More fundamentally, the conclusion we draw, that the variance of the sample mean in the $n/N\to0$ regime is the marginal variance of the response, is not available conditionally on the graph. That conclusion requires averaging over the network, which none of these frameworks does. The closest exception is \citet{li2022random}, who study treatment effect estimation under interference on a graphon-generated random graph and average over the graph as we do; their outcome model is an exposure-based potential-outcomes specification rather than a propagation mechanism, treatment is randomized, and the sampling-fraction dichotomy does not arise.

\paragraph*{Exchangeable arrays}
A second line treats network data as an exchangeable array and uses the Aldous--Hoover representation to reduce dependence to independent latent variables. Laws of large numbers and CLTs for dissociated exchangeable arrays go back to \citet{silverman1976limit,eagleson1978limit}, with uniform and bootstrap versions in \citet{davezies2021empirical} and \citet{menzel2021bootstrap}, and subsampling and resampling for graphon-generated network statistics in \citet{bhattacharya2015subsampling,lunde2023subsampling}. Our paper \citet{li2025assumption} is in this tradition and is the point of departure for the present paper. The common requirement is infinite exchangeability: the finite array must be the restriction of an infinitely exchangeable one, so that the Aldous--Hoover representation applies and the $U$-statistic machinery follows. The SAR array is finitely exchangeable for each $N$ but not projectively consistent (Remark~\ref{rem:exch}), so this machinery does not apply directly. 
\paragraph*{Sampling from dependent populations.}
That a sparse sample from a dependent population behaves almost as an independent one has
precedents outside network analysis. In cluster sampling, the design effect depends on the
intraclass correlation and the number of sampled units per cluster \citep{kish1965survey}.
\citet{abadie2020sampling} show how the ratio of sample to population size determines
whether design-based or sampling-based uncertainty is relevant for regression, in a setting
without interference, where each unit's outcomes do not depend on the other units.
Incomplete $U$-statistics, computed from randomly sampled index tuples, have a variance
that interpolates between that of the complete statistic and that of the sampling stage
\citep{blom1976some,janson1984asymptotic,chen2019randomized}. Our setting shares the
two-stage form of these results but differs in mechanism: we observe outputs that each
depend on the entire population through the network, and the dependence that this
creates can make the variance of the population mean smaller than the marginal variance
as well as larger (Remark~\ref{rem:two-stage}).

\paragraph*{Inference for spatial autoregressive models}
The spatial autoregressive model of \citet{cliff1973spatial} has a mature inferential theory in econometrics, including instrumental-variable and GMM estimators \citep{kelejian1998generalized,kelejian1999generalized,kelejian2010specification} and the quasi-maximum-likelihood analysis of \citet{lee2004asymptotic}, who showed that convergence rates depend on the structure of the weight matrix and that irregular information arises when each unit has many neighbors. The linear-in-means model of peer effects \citep{manski1993identification,bramoulle2009identification,lin2010identifying} is the same model with a constant coefficient. In all of these the weight matrix is fixed and known, the autoregressive coefficient is a common scalar to be estimated, and asymptotics are conditional on the design. Heterogeneous spatial lag coefficients are studied by \citet{aquaro2021estimation}, who estimate unit-specific coefficients in a panel with both dimensions growing; identification of node-level coefficients requires repeated observation, and the coefficients are treated as fixed parameters rather than random effects. Endogeneity of the network, meaning dependence between link formation and the outcome through latent attributes, is addressed by control-function and latent-position approaches \citep{johnsson2021estimation,auerbach2022identification}, and partially observed networks by \citet{boucher2022estimating,lewbel2023social}; these model the network as random in order to correct the estimator, but the target remains a structural coefficient and the asymptotic theory remains conditional on the realized or estimated graph. Our setting differs in three ways: the coefficient and the input are nonparametric functions of the latent type with an unrestricted joint law, nothing about them is estimated, and the target is the mean response, for which we show that the network can be averaged out rather than conditioned on.

\paragraph*{Network time series}
A parallel literature studies autoregression on networks over time. The network vector autoregression of \citet{zhu2017network} regresses each node's response on its own lag and on the lagged neighborhood average, with a common influence coefficient; extensions include quantile \citep{zhu2019network}, grouped \citep{zhu2020grouped}, community-level \citep{chen2023community}, nonlinear and count-valued \citep{armillotta2023nonlinear,armillotta2024count}, and generalized \citep{knight2020generalized} versions. \citet{wu2023random} is closest to the present paper in that it places a graphon model on the network and uses latent positions to separate influence from homophily; its asymptotic theory is nonetheless conditional on the realized graph and requires $T\to\infty$, with the graph entering through spectral-gap and influential-node conditions on the realized transition matrix. These models are dynamic and their inference uses time as the source of replication: stationarity in $t$ and $T\to\infty$ are either required or heavily used, and the network is fixed throughout. The static model \eqref{eq:SAR} is the equilibrium of the noise-free dynamic $Y^{(t)}=\Gamma PY^{(t-1)}+Z$, which converges to $(I-\Gamma P)^{-1}Z$ under $\bar\gamma<1$; the mean of the stationary distribution of a network autoregression with temporal innovations is the same object. We observe a single snapshot, use $N\to\infty$ as the only source of replication, and let the network be random.

\paragraph*{Summary}
Table~\ref{tab:related} records, for each line of work, whether the network is averaged over, whether the array is required to be infinitely exchangeable, whether dependence is required to decay in network distance, and whether inference needs the network. The present paper is, to our knowledge, the first to combine a random network, finite exchangeability, global dependence, and network-free inference in a single result.
\begin{table}[t]
\centering
\footnotesize
\setlength{\tabcolsep}{5pt}
\begin{tabular}{lcccc}
\toprule
 & Network & Infinite & Decay in & Inference \\
 & averaged over & exchangeability & network distance & needs network \\
\midrule
Fixed-network CLTs        & no  & no  & yes & yes \\
Causal network inference  & no  & no  & yes & yes \\
Exchangeable arrays       & yes & yes & no  & no \\
SAR econometrics          & no  & no  & n/a & yes \\
Network time series       & no  & no  & n/a & yes \\
This paper                & yes & no  & no  & no ($n/N\to0$) \\
\bottomrule
\end{tabular}
\caption{Positioning relative to prior work. Fixed-network CLTs: \citet{kojevnikov2021limit,leung2026normal,xu2024functional}; causal network inference: \citet{ogburn2024causal}; exchangeable arrays: \citet{davezies2021empirical,li2025assumption}; SAR econometrics: \citet{lee2004asymptotic,aquaro2021estimation}; network time series: \citet{zhu2017network,wu2023random}. ``Decay in network distance'' refers to a weak-dependence condition bounding dependence between distant nodes; ``n/a'' indicates that the question does not arise because asymptotics are conditional on a fixed design.}
\label{tab:related}
\end{table}

\section{Main results}\label{sec:main}
In this section, we state our main theorems; novel technical tools used to derive this result are discussed in Section \ref{sec:proof-outline}.

\subsection{Conditions}\label{sec:conditions}

We write $g(x)=\int_0^1w_N(x,v)\,dv$ for the expected-degree profile, so that node $i$ has expected degree $(N-1)\rho_Ng(\xi_i)$, and $p(v\mid x)=w_N(x,v)/g(x)$ for the transition density of a step of the walk from type $x$. The mixed norm
\[
\mathsf M_N:=\int_0^1\Big(\int_0^1w_N^2(x,v)\,dv\Big)^2dx
\]
measures the second moment of a node's edge weights, averaged over nodes; it is finite for bounded graphons and for many unbounded ones, and it is the quantity through which hubs enter the conditions. All constants below are uniform over the array. Fix $\delta\in(0,1]$.

\begin{assumption}[Common conditions]\label{ass:common}
\leavevmode
\begin{enumerate}[label=(M\arabic*),leftmargin=*]
\item$\sup_x|\gamma(x)|\le\bar\gamma<1$.
\item $w_N\ge c_0>0$ almost everywhere.
\item $\rho_Nw_N\le1/2$.
\item $c_0(N-2)\rho_N\ge C_*\log N$ for a constant $C_*$ depending only on $(\bar\gamma,c_0,\delta)$.
\end{enumerate}
\end{assumption}

(M1) is the contraction condition that makes \eqref{eq:SAR} well posed; it is also the condition under which a walk has geometrically decaying memory, and it is the only condition on the (random) autoregressive coefficient. (M2) excludes pathologies with small degree nodes and is what stabilizes the degree normalization; see Remark \ref{rem:M2} below for further discussions and possible relaxations. (M3) is a mild condition that prevents the value of $w_N$ from being censored; here $1/2$ is arbitrary, and this condition can be relaxed to any constant less than 1 at the cost of a change in constants. (M4) is the sparsity condition: average degree of order $\log N$ is the regime in which degrees concentrate around their conditional means, which is the minimal regime in which one can expect the degree normalization to be stable.  In setting (B) below, the network fluctuation bound is of order $\mathsf M_N\log N/(N\rho_N)$, so (B3) implicitly requires $N\rho_N/\log N\to\infty$,
which is slightly stronger than (M4); in setting (U), (M4) alone suffices.

\begin{remark}[On the degree lower bound]\label{rem:M2}
In the regime $N\rho_N\asymp\log N$ of (M4), (M2) is not an additional restriction on degrees: it requires expected degrees of order $\log N$, below which $1/d_i$ does not concentrate and the degree normalization is unstable for any argument. In denser regimes (M2) could be weakened, in setting (B), to a minimum expected degree of order $\log N$, provided the degree disparity $\|w_N\|_\infty/\inf g$ is $o(\sqrt n)$, or $o(N/\sqrt n)$ when $n/N\to0$. The disparity bound is a Feller-type negligibility condition on the influence of a single node on the mean, since a node of type $v$ influences a node of type $x$ by at most $w_N(x,v)/g(x)$; (U1) is this condition under (M2). We do not pursue the weakening, as the single-regime statement is simpler and the gain is confined to dense graphs.
\end{remark}

The two settings below are complementary. The first has bounded inputs and asks almost nothing of the graphon beyond integrability; the second allows heavy-tailed inputs and pays for them with a mild growth condition on the largest edge weight and with fourth-moment-type conditions on the degree profile.

\begin{assumption}[Setting (B): bounded inputs]\label{ass:B}
\leavevmode
\begin{enumerate}[label=(B\arabic*),leftmargin=*]
\item $|Z|\le B$.
\item $\int_0^1 \int_0^1 w_N^{2+\delta} (u,v) \ du \ dv\le K$.
\item $\mathsf M_N\log N=o(N\rho_N)$.
\end{enumerate}
\end{assumption}

\begin{assumption}[Setting (U): heavy-tailed inputs]\label{ass:U}
\leavevmode
\begin{enumerate}[label=(U\arabic*),leftmargin=*]
\item $\|w_N\|_\infty^2=o(N)$.
\item $\E[g(\xi)Z(\xi)^4]\le K$.
\item $\int_0^1\big(\int_0^1 w_N^{2+\delta}(x,v)\,dx\big)^{4/(2-\delta)}dv\le K$.
\item $\E g(\xi)^4\le K$, $\mathsf M_N\le K$, and the two- and three-point functionals $\iint g(x)^2W_2(x,z)^2\,dx\,dz$, $\int_0^1 \mathfrak h^2(u)\, du $ and $\int_0^1 g^3(u)(Tg)^2(u) \ du$ are at most $K$, where $W_2(x,z)=\int_0^1 w_N(x,v)w_N(v,z)\,dv$ and $\mathfrak h(z)=\int_0^1\big(\int_0^1 w_N(v,u)^2du\big)w_N(v,z)\,dv$.
\end{enumerate}
\end{assumption}

Setting (B) places no condition on $\|w_N\|_\infty$ beyond (M3), so expected degrees of order $N$ are allowed for a vanishing fraction of nodes; (B3) holds whenever $\mathsf M_N\le K$ and $N\rho_N/\log N\to\infty$, and permits $\mathsf M_N$ to grow slowly. In setting (U), (U1) is a deterministic, Feller-type negligibility condition on the influence of a single node. Since $\int_0^1 \int_0^1 w_N \ du \ dv=1$, the ratio of the largest expected degree to the average expected degree is at most $\|w_N\|_\infty$, so (U1) allows this ratio to grow as fast as $o(\sqrt N)$ (Remark~\ref{rem:M2}); it is used in the signal remainder and in the feedback and coincidence terms of the network fluctuation. (U2) is a weighted fourth moment of the input in which the weight $g$ is the stationary density of the walk. The remaining conditions in (U3)--(U4) are integrated moments of fixed order.

Some control of the tails is needed, and the two settings control different ones. A large input at a hub is transmitted to many neighbors at once, so the tail of $Z$ and the tail of $w_N$ interact in the influence a single node can exert on the mean. Setting (B) bounds the input, and can then allow hubs with expected degree of order $N$; setting (U) bounds the hub size, $\|w_N\|_\infty=o(\sqrt N)$, and can then allow heavy-tailed inputs, with the interaction appearing directly in (U2), a fourth moment of the input weighted by the degree profile. In neither setting does the walk length add a further factor, even though \eqref{eq:neumann} involves products of arbitrary length; this is a consequence of the Markov chain structure exploited in our proofs.

\begin{remark}[(U1) and realized hubs]\label{rem:U1-whp}
Condition (U1) is a deterministic bound on $w_N$, but the proofs use $w_N$ only
through its values at pairs of realized types. If $w$ is unbounded and
$\max_{i\ne j}w(\xi_i,\xi_j)=o_P(\sqrt N)$, then the truncated graphons
$w_N\propto\min(w,W_N)$, with a suitable $W_N=o(\sqrt N)$, satisfy (U1), converge to
$w$ in $L^1$ and hence in cut norm, and generate the same graph as $w$ with
probability tending to one. The relevant restriction is therefore on the hubs that
are realized, not on the graphon itself. For example, in a rank-one graphon with
sociability $a(x)\asymp(1-x)^{-\alpha}$, the largest realized weight is of order
$N^{2\alpha}$, and this requirement amounts to $\alpha<1/4$, close to the fourth-moment
condition of (U4). Transferring the conclusions from $w_N$ to $w$ also requires that
the centering and the limiting variance change negligibly under truncation, which
holds under the moment conditions; we do not pursue the details.
\end{remark}

\subsection{The limiting variance}\label{sec:variance}

The limit theorem is stated in terms of the following population objects, which are defined for each $N$ from $(w_N,\gamma,z)$ alone.  See Section \ref{sec:operators} for further interpretation of these population quantities.

Let $T$ be the Markov operator of the walk in type space, $(Tf)(x)=\int_0^1 p(v\mid x)f(v)\,dv$, and let $T_\gamma=\gamma T$, that is, $(T_\gamma f)(x)=\gamma(x)(Tf)(x)$: the coefficient multiplies after averaging over the next step. The mean-field response is
\begin{equation}\label{eq:mean-field}
y^\circ:=(I-T_\gamma)^{-1}z=\sum_{k\ge0}T_\gamma^kz ,
\end{equation}
which is \eqref{eq:neumann} with the empirical walk replaced by the walk in the latent space; $Y_i^\circ :=y^\circ(\xi_i)$ is the response of a node of type $\xi_i$ in an infinite population. Its mean $\mu^\circ_N:=\int_0^1y^\circ(u) \ du$ is the mean-field target, which differs from $\mu_N=\E Y_1$ by a bias of order $\log N/(N\rho_N)$ up to lower-order terms, characterized in Section~\ref{sec:bias}.

A node of type $v$ influences the response of a node of type $x$ in two ways: directly, because $v$ may be a neighbor of $x$ and contribute $y^\circ(v)$ to the neighborhood average, and by dilution, because the presence of $v$ among $x$'s neighbors reduces the weight of all other neighbors. The influence kernel
\begin{equation}\label{eq:influence}
f_v(x):=\gamma(x)\,p(v\mid x)\big(y^\circ(v)-(Ty^\circ)(x)\big)
\end{equation}
records both: the first term is the direct contribution and the second is the dilution, which is the neighborhood average $x$ would have had without $v$. The influence propagates through the network, and the total effect of a node of type $v$ on a node of type $x$ is
\[
\phi(x,v):=\big((I-T_\gamma)^{-1}f_v\big)(x).
\]
Averaging over $x$ gives the influence of a type-$v$ node on the population mean. Define
\begin{equation}\label{eq:H}
H_{\mathrm{root}}(v):=y^\circ(v),\qquad H_{\mathrm{rest}}(v):=\int_0^1\phi(x,v)\,dx,\qquad H_{\mathrm{full}}:=H_{\mathrm{root}}+H_{\mathrm{rest}} ,
\end{equation}
and, for $n\le N$,
\begin{equation}\label{eq:sigma}
\sigma^2_{n,N}:=\Big(1-\frac nN\Big)\Var H_{\mathrm{root}}(\xi_1)+\frac nN\Var H_{\mathrm{full}}(\xi_1).
\end{equation}
$H_{\mathrm{root}}(\xi_i)$ is what node $i$ contributes to $\bar Y_S$ through its own response; $H_{\mathrm{rest}}(\xi_i)$ is what it contributes through everyone else's. A node in the sample contributes both; a node outside the sample contributes only the second, and \eqref{eq:sigma} is the resulting mixture. The formula is exact at every $(n,N)$ and is not an asymptotic approximation: it is the variance of the Hoeffding projection of the signal component of $\bar Y_S$ at population size $N$.

\subsection{The central limit theorem}\label{sec:clt}

Throughout, $S\subseteq[N]$ is a deterministic set of $n$ observed nodes, $\bar Y_S=n^{-1}\sum_{i\in S}Y_i$, and $\mu_N=\E Y_1$. By exchangeability, $\mu_N$ does not depend on $S$. The non-degeneracy condition
\begin{itemize}
\item[(M5)] $\liminf_N\sigma^2_{n,N}>0$
\end{itemize}
rules out the degenerate case in which the response carries no variation at the leading order.

\begin{theorem}[CLT for the SAR mean]\label{thm:clt}
Assume (M1)--(M4), (M5), and either setting (B) or setting (U). Then for every sequence of deterministic $S\subseteq[N]$ with $|S|=n\to\infty$, $n\le N$,
\[
\frac{\sqrt n\,(\bar Y_S-\mu_N)}{\sigma_{n,N}}\;\Rightarrow\;N(0,1).
\]
\end{theorem}

The proof is outlined in Sections~\ref{sec:operators}--\ref{sec:fluctuation}; the complete proof is in the appendix.

\begin{corollary}[Sampling regimes]\label{cor:regimes}
Under the conditions of Theorem~\ref{thm:clt}, if $n/N\to c\in[0,1]$ and $\Var H_{\mathrm{root}}(\xi_1)\to v_0$, $\Var H_{\mathrm{full}}(\xi_1)\to v_1$, then
\[
\sqrt n\,(\bar Y_S-\mu_N)\Rightarrow N\big(0,\,(1-c)v_0+cv_1\big).
\]
In particular, when $c=0$ the asymptotic variance is the marginal variance of the mean-field response, $v_0=\lim\Var y^\circ(\xi_1)$, and the network does not appear.
\end{corollary}

Three regimes are worth separating. When $n/N\to0$, the observed nodes are a vanishing fraction of the population, their influence on each other through the network is a vanishing part of the variance, and $\bar Y_S$ behaves as the mean of i.i.d.\ draws of $y^\circ(\xi)$. When $n/N\to c\in(0,1)$, the sample is a non-vanishing fraction, each sampled node also appears in others' neighborhoods, and the variance picks up the influence term $c\,(v_1-v_0)$; this term can be of either sign, since $H_{\mathrm{rest}}$ may be negatively correlated with $H_{\mathrm{root}}$, so that network dependence can reduce as well as inflate the variance of the mean. When $n=N$ the full population is observed and the variance is $v_1=\Var(H_{\mathrm{root}}+H_{\mathrm{rest}})$, the variance of a node's total influence on the population mean, which is the form a Hájek projection would take for a $U$-statistic of growing rank.

\begin{remark}[Weak dependence and the relation to two-stage sampling]\label{rem:two-stage}
Since the responses are exchangeable, the variance of their mean over any fixed set of
$n$ nodes satisfies $n\Var(\bar Y_S)=\Var(Y_1)+(n-1)\Cov(Y_1,Y_2)$, and
Theorem~\ref{thm:clt} amounts to identifying the two quantities on the right:
$\Var(Y_1)\to\Var H_{\rm root}$ and $N\Cov(Y_1,Y_2)\to\Var H_{\rm full}-\Var H_{\rm root}$.
The identity has the familiar two-stage form of finite-population sampling: since a fixed
$S$ has the same law as a uniformly random one, $n\Var(\bar Y_S)$ equals the expected
variance of sampling $n$ of the $N$ realized responses plus $n/N$ times $N\Var(\bar Y_N)$, the
fluctuation of the population mean. The connection is one of form rather than content. No
sampling takes place in Theorem~\ref{thm:clt}: $S$ is fixed, and the variance is taken over
the network-forming process. Even the first stage is not classical, since the spread of the
realized responses includes graph noise, and its convergence to a network-free limit is part of Proposition~\ref{prop:svar}. The content of the
result is that, in spite of feedback along every path, sparsity at the level of $\log N$,
and heavy-tailed graphons and inputs, the dependence between responses
is weak, of order $1/N$ for every pair, and is identified exactly. Under strong dependence
this can fail: a component shared by all responses at order one would correlate every pair
at order one, and no sampling fraction would restore network-free inference. In a classical
superpopulation with i.i.d.\ values, the pairwise covariance vanishes and the sampling
fraction does not affect the limiting variance; the interpolation in \eqref{eq:sigma} is
non-trivial precisely because network dependence makes the covariance of order exactly
$1/N$ and its sign can go either way, so that network dependence can reduce the variance as well as increase it. 
\end{remark}

\begin{remark}[Where hubs matter]\label{rem:hubs}
Hubs enter through $H_{\mathrm{rest}}$, which grows with the degree profile. If $H_{\mathrm{rest}}(\xi_1)^2$ fails to be uniformly integrable across $N$, a vanishing fraction of hubs carries a non-vanishing share of the variance and the limit need not be Gaussian; the Lyapunov conditions (B2), respectively (U2)--(U3), are what exclude this. In setting (U), (U1) additionally controls the influence of a single hub in the remainder terms.
\end{remark}

\begin{remark}[Random subsamples]\label{rem:random-S}
All statements are for deterministic $S$ and hold uniformly in $S$. If $S$ is drawn uniformly at random from the $n$-subsets of $[N]$, independently of everything else, the conclusions hold conditionally on $S$ and therefore unconditionally. Other sampling designs that depend on the network are not covered by our result.
\end{remark}

\subsection{Variance estimation and network-free inference}\label{sec:studentized}

Let $s_S^2=(n-1)^{-1}\sum_{i\in S}(Y_i-\bar Y_S)^2$ be the ordinary sample variance of the observed responses.

\begin{proposition}[Sample variance]\label{prop:svar}
Assume (M1)--(M4) and either setting (B) or setting (U). Then for every sequence of deterministic $S$ with $|S|=n\to\infty$, $n\le N$,
\[
s_S^2-\Var H_{\mathrm{root}}(\xi_1)\to0\quad\text{in probability}.
\]
More precisely,
\[
\E\left[\frac{1}{|S|}\sum_{i\in S}\big(Y_i-y^\circ(\xi_i)\big)^2 \right]\le C\Big(\frac{\E[gZ^2]}{N\rho_N}+\frac{\|\phi\|^2_{L^2}}{N}+\frac{r_N^2}{N}\Big)+CN^{-13}\,\E Z^2,
\]
where $r_N\to0$ is the remainder rate of the signal component.
\end{proposition}

Above, the term $CN^{-13}\E Z^2$ accounts for the event, of small probability, on which the
degrees fail to concentrate. Its order is not intrinsic. The tail bounds behind the
concentration events are exponential in the average degree, which under (M4) may be
of order only $\log N$, so they give probability $O(N^{-K})$ for any fixed $K$, at the cost of a larger constant $C_*$, but not exponentially small probability in $N$. Here
any $K>1$ would suffice; we fix $K=14$ throughout for such events for definiteness,
which also covers the network analysis, where some steps require bad events of
probability $o(N^{-6})$.

The proposition says that each observed response is close to its mean-field value, $Y_i=y^\circ(\xi_i)+e_i$ with $\frac{1}{|S|}\sum_{i \in S}e_i^2\to0$ at the rate $1/(N\rho_N)$, so that the empirical variance of the $Y_i$ is that of i.i.d.\ draws of $y^\circ(\xi)$. It holds for every $n\le N$, including the full population, and does not require $n/N\to0$. What requires $n/N\to0$ is only the identification of $\Var H_{\mathrm{root}}$ with the variance in Theorem~\ref{thm:clt}.

A corollary of  the above result is that in the $n/N \rightarrow 0$ regime, the sample variance is consistent and the usual Normal confidence interval provides asymptotically valid coverage.  It should be emphasized that the sample variance estimator uses neither the network, nor the graphon, nor $\gamma$, nor the sampling fraction.
\begin{corollary}[Studentized CLT]\label{cor:t}
Assume the conditions of Theorem~\ref{thm:clt} and $n/N\to0$. Then
\[
\frac{\sqrt n\,(\bar Y_S-\mu_N)}{s_S}\Rightarrow N(0,1).
\]
\end{corollary}

\begin{remark}[When the sample is not a vanishing fraction]\label{rem:not-vanishing}
If $n/N\to c\in(0,1]$, Proposition~\ref{prop:svar} still gives $s_S^2\to v_0$, while the variance of $\sqrt n\,\bar Y_S$ tends to $(1-c)v_0+cv_1$. Estimating $v_1$ requires information about the network: by Remark~\ref{rem:two-stage}, the
excess $c(v_1-v_0)$ comes from the fluctuation of the population mean, a component shared by
all sampled nodes, of which the sample observes a single realization. As for the variance
of a random effect observed in a single cluster, the responses within one population carry
no consistent information about it, because the shared shift is confounded with the unknown
mean. This is the precise sense in which network-free inference is available for vanishing subsamples but not in general, and the precise sense in which network dependence is a nuisance rather than an obstruction: its effect on the variance is an explicit, interpretable term that enters only through the sampling fraction.
\end{remark}

\begin{remark}[Rate of variance estimation]\label{rem:svar-rate}
The error in $s_S^2$ is $O_P(n^{-1/2})$ from the oracle term, under four moments of $y^\circ(\xi)$, plus a network term of order $(N\rho_N)^{-1/2}$. The network therefore does not slow variance estimation beyond the sparsity condition already needed for the CLT.
\end{remark}

Corollary~\ref{cor:t} and Remark~\ref{rem:not-vanishing} together suggest that whether network dependence can be ignored is not a property of the data but of the relation between the observed nodes and the population through which the dependence propagates. The decisive question is at what level the observed network is representative. If the $n$ sampled individuals interact mainly with a much larger population of which the observed network is a small part, as in a survey of individuals drawn from a city, a cohort drawn from a platform, or any design in which the units are sampled and the network is whatever connects them, then $n/N\to0$ is the right description, the response of each sampled node is driven by nodes that are not in the sample, and the ordinary $t$-interval is valid with no model for the network and no knowledge of it. If instead the observed network is the population, a complete classroom, firm or village, then $n=N$, the sampled nodes influence each other, and the variance carries the influence term, which cannot be estimated from the responses alone. In many applications the honest answer lies between, or is uncertain: the observed network is a window onto something larger, but how much larger is not known.

In that case there is a simple remedy. Drawing a random downsample of size $n_0=o(N)$ from the observed nodes and applying Corollary~\ref{cor:t} to it restores valid inference, since the subsample is a vanishing fraction of the observed network whether or not the observed network is itself a fraction of something larger. The cost is the rate: $n_0^{-1/2}$ in place of $N^{-1/2}$, which is the price of the influence term one has chosen not to model. When the network is unobserved or incompletely observed, or when one is unwilling to commit to a model for how the response propagates along it, this is essentially the only route to a valid interval, and the results above show that it costs nothing beyond the rate. 

\subsection{The finite-population target and the mean-field target}\label{sec:bias}

The centering in Theorem~\ref{thm:clt} is $\mu_N=\E Y_1$, the mean response in a population of size $N$. The mean-field value $\mu^\circ_N=\int_0^1y^\circ(u)\,du=\E Y^\circ_1$ is the mean response in the infinite-population limit with type distribution $U$ and connection profile $w_N$; when $w_N\equiv w$ it does not depend on $N$, and it is the quantity that a fixed-graphon analysis would call the parameter. The two differ by a bias that we now quantify, because it determines when the sample mean can be centered at $\mu^\circ_N$ instead.

\begin{proposition}[Finite-population bias]\label{prop:bias}
Assume (M1)--(M4). In setting (B),
\[
\big|\mu_N-\mu^\circ_N\big|\le CB\Big(\frac{\log N}{N\rho_N}+\frac1N\Big(\frac{\mathsf M_N}{\rho_N}\Big)^{1/2}\Big),
\]
and in setting (U), $\big|\mu_N-\mu^\circ_N\big|\le C\big(\log N/(N\rho_N)+\|w_N\|_\infty/N\big)$.
\end{proposition}

The bias is of the order of one over the average degree, up to a logarithm and a lower-order term, and it comes mainly from the edge noise: the expected-weight response $Y^W$ has mean within $O(r_N/\sqrt N)$ of $\mu^\circ_N$ by the argument of Section~\ref{sec:reduction} applied to $S=[N]$, so the signal component contributes at most $O(r_N/\sqrt N)$, and the rest of $\mu_N-\mu^\circ_N$ is $\E V_1=\E\big[\E[Y_1\mid\xi]-Y^W_1\big]$, the mean effect of replacing edge probabilities by edges.

\begin{corollary}[Centering at the mean-field target]\label{cor:fixed-centering}
Under the conditions of Theorem~\ref{thm:clt}, if $\sqrt n\,\log N/(N\rho_N)\to0$ then $\sqrt n\,(\bar Y_S-\mu^\circ_N)/\sigma_{n,N}\Rightarrow N(0,1)$, and when $n/N\to0$ the sample standard deviation $s_S$ may replace $\sigma_{n,N}$.
\end{corollary}

The second term in Proposition~\ref{prop:bias} is negligible at the $\sqrt n$ scale in both settings, since $\sqrt n(\mathsf M_N/\rho_N)^{1/2}/N\le(\mathsf M_N/(N\rho_N))^{1/2}\to0$ under (B3) and $\sqrt n\|w_N\|_\infty/N\le\|w_N\|_\infty/\sqrt N\to0$ under (U1). So the mean-field target is available exactly when the average degree is of larger order than $\sqrt n\log N$. Below that threshold the bias is visible at the $\sqrt n$ scale, and $\mu_N$ is the only target that the sample mean estimates without a model.

\subsection{Immediate consequences}\label{sec:extensions}
Central limit theorems for the mean play a foundational role in statistics not only because the mean is of intrinsic interest but also because many functionals of interest are smooth functions of means. In the random-design and random-network setting, least-squares coefficients are smooth functions of averages of products of node-level variables \citep{li2025assumption}. We give sufficient conditions under which vector means of the form $\bar U_S= \frac{1}{|S|} \sum_{i \in S}X_iY_i$, with $X_i=X(\xi_i)$ a node-level weight or covariate, satisfy a central limit theorem, and conclude asymptotic normality for regression coefficients by the delta method.

\subsubsection{Smooth functions}
Let $Z^{(1)},\dots,Z^{(d)}$ be inputs each satisfying the conditions of one of the two settings, with a common $\gamma$, let $X^{(1)},\dots,X^{(d)}$ be $\xi$-measurable node-level variables, and let $\bar U^{(j)}_S=\frac{1}{|S|}\sum_{i \in S}X^{(j)}_iY^{(j)}_i$ with $Y^{(j)}=(I-\Gamma P)^{-1}Z^{(j)}$; $X^{(j)}\equiv1$ recovers $\bar Y^{(j)}_S$. Write $\Phi^{(j)}_2(x):=\int_0^1\phi^{(j)}(x,v)^2\,dv$ for the influence profile of the $j$th response. For simplicity, we assume $\sup_x |X^{(j)}(x)| \leq K$ for $j \in [d]$; see Remark \ref{rem:bounded-cov} for further discussion.
The influence objects become
\[
H^{(j)}_{\mathrm{root}}:=X^{(j)}y^{\circ(j)},\qquad H^{(j)}_{\mathrm{rest}}(v):=\int_0^1X^{(j)}(x)\,\phi^{(j)}(x,v)\,dx ,
\]
the influence of a node of type $v$ on each node, weighted by that node's covariate; $H^{(j)}_{\mathrm{rest}}\in L^2$ by boundedness of $X^{(j)}$, since $\|\phi^{(j)}\|_{L^2}< \infty$. We also assume that the smallest eigenvalue of $\Sigma_{n,N}$ below is bounded away from zero.

\begin{corollary}[Vector means; smooth functions]\label{cor:delta}
Under the conditions above, $\sqrt n\,(\bar{\mathbf U}_S-\boldsymbol\mu_N)\Rightarrow N(0,\Sigma_{n,N})$ with $\boldsymbol\mu_N=\E\,\bar{\mathbf U}_S$ and $\Sigma_{n,N}=(1-n/N)\,\mathrm{Cov}\,\mathbf H_{\mathrm{root}}+(n/N)\,\mathrm{Cov}\,\mathbf H_{\mathrm{full}}$, read row-wise if the limits do not exist. If $F:\mathbb R^d\to\mathbb R$ is differentiable near $\boldsymbol\mu_N$ with gradient bounded away from zero and uniformly continuous there, then $\sqrt n\,(F(\bar{\mathbf U}_S)-F(\boldsymbol\mu_N))/(\nabla F^\top\Sigma_{n,N}\nabla F)^{1/2}\Rightarrow N(0,1)$. When $n/N\to0$, the sample covariance of the observed vectors $(X^{(j)}_iY^{(j)}_i)_j$ is consistent for $\Sigma_{n,N}$, and the plug-in studentization is valid without the network.
\end{corollary}

\begin{remark}[Bounded covariates]\label{rem:bounded-cov}
We state Corollary~\ref{cor:delta} for bounded weights $X^{(j)}$, since the purpose of
this section is to illustrate the reach of the main theorem rather than to find the
weakest regularity conditions. The signal part requires only the moment conditions $\E|X^{(j)}|^{2+\delta} <\infty$ and $\int_0^1 X^{(j)}(x)^2\Phi_2^{(j)}(x) \ dx < \infty$, and only the network part uses boundedness. We expect it can be relaxed:
by duality, unbounded weights on the sampled nodes play the role on the dual side
that unbounded inputs play on the primal side, which setting (U) already handles
under moment conditions. Our current proof uses bounded weights in two places,
through the bounds $\|b\|_\infty\le C/(N\rho_N)$ and $\sum_md_m|b_m|\le Cn$ on the dual vector,
which become random when the weights are unbounded. Replacing them would require
moment computations analogous to those of setting (U), which we do not pursue here.
\end{remark}

\subsubsection{Regression on network-derived covariates}\label{sec:network-covariates}

In network-linked regression the covariates are often computed from the network itself: degrees and other local subgraph counts, or coordinates of a spectral embedding. In \citet{li2025assumption} we studied such regressions when the responses are independent given the latent positions; here the responses are generated through the network. We show that in the regime $n/N\to0$ the conclusions of that paper carry over, and we recover the same distinction between two targets, with the same type of sparsity gap between them. Throughout this subsection $n/N\to0$, and the covariates are computed from the full network on $N$ nodes.

Let $\hat \ell_i\in\mathbb R^d$ be computed from $(A,\xi)$ and suppose it decomposes as
\begin{equation}\label{eq:cov-decomp}
\hat \ell_i=\ell(\xi_i)+\zeta_i ,
\end{equation}
where $\ell:[0,1]\to\mathbb R^d$ is measurable and $\zeta_i$ is network noise; for local subgraph counts, $\ell(\xi_i)=\E[\widetilde Q_i(\mathcal R)\mid\xi_i]$ with $\widetilde Q_i$ the latent count in which each edge indicator is replaced by $w_N(\xi_j,\xi_k)$. Write $Y_i=y^\circ(\xi_i)+e_i$ as in Proposition~\ref{prop:svar}, and let $X_i=X(\xi_i)\in\mathbb R^p$ be conventional $\xi$-measurable covariates with $X_{i,1}=1$. Set $\hat L_i=(X_i^\top,\hat \ell_i^\top)^\top$ and $L_i=(X_i^\top,\ell(\xi_i)^\top)^\top$. The two targets are
\[
\widetilde\beta_N:=\big(\E\,\hat L_1\hat L_1^\top\big)^{-1}\E\,\hat L_1Y_1,\qquad
\beta^*_N:=\big(\E\,L_1L_1^\top\big)^{-1}\E\,L_1Y_1 ,
\]
the projection of $Y$ onto the noisy and onto the latent covariates respectively, both at population size $N$, and the estimator is the least-squares coefficient $\hat\beta=\hat\Lambda^{-1}\hat\gamma$ with $\hat\Lambda=\frac{1}{|S|}\sum_{i \in S}\hat L_i\hat L_i^\top$, $\hat\gamma=\frac{1}{|S|}\sum_{i \in S}\hat L_iY_i$.

\paragraph*{Conditions on the covariates}
We assume that $X$ and $\ell$ are bounded, so that the weighted-means corollary applies to $\frac{1}{|S|}\sum_{i \in S}X_iY_i$, $\frac{1}{|S|}\sum_{i \in S}\ell(\xi_i)Y_i$, and the entries of $\frac{1}{|S|}\sum_{i \in S}L_iL_i^\top$; that $\E\,L_1L_1^\top$ has smallest eigenvalue bounded away from zero; and that the network noise satisfies
\begin{itemize}
\item[(C1)] $\E\|\zeta_1\|^2\to0$;
\item[(C2)] $\sqrt n\,\frac{1}{|S|}\sum_{i \in S}\big(y^\circ(\xi_i)\zeta_i-\E[y^\circ(\xi_i)\zeta_i]\big)\to0$ and $\sqrt n\,\frac{1}{|S|}\sum_{i \in S}\big(L_i\zeta_i^\top-\E[L_i\zeta_i^\top]\big) \xrightarrow{P}0$;
\item[(C3)] $\sqrt n\,\frac{1}{|S|}\sum_{i \in S}\big(e_i\zeta_i-\E[e_i\zeta_i]\big)\to0$ and $\sqrt n\,\frac{1}{|S|}\sum_{i \in S}\big(\zeta_i\zeta_i^\top-\E[\zeta_i\zeta_i^\top]\big)\xrightarrow{P}0$ .
\end{itemize}
(C1) says the covariate concentrates around its latent version; (C2) says its noise, weighted by the latent quantities, averages out at the $\sqrt n$ scale; (C3) says the same for the two products of noises, that of the response with that of the covariate and that of the covariate with itself. By exchangeability the centerings in (C2)--(C3) do not depend on $i$.

\begin{corollary}[Network-derived covariates]\label{cor:netcov}
Assume the conditions of Theorem~\ref{thm:clt}, $n/N\to0$, and the conditions above. Then
\[
\sqrt n\,(\hat\beta-\widetilde\beta_N)\Rightarrow N(0,\Sigma_\beta),\qquad
\Sigma_\beta=\Lambda^{-1}\,\mathrm{Cov}\big(L_1(Y_1-L_1^\top\beta^*_N)\big)\,\Lambda^{-1},\quad\Lambda=\E\,L_1L_1^\top,
\]
read row-wise if the limit does not exist, and the sandwich estimator built from the observed $(Y_i,\hat L_i)_{i\in S}$ is consistent for $\Sigma_\beta$. If in addition
\begin{equation}\label{eq:bias-cond}
\sqrt n\,\big\|\E[\hat L_1\hat L_1^\top]-\E[L_1L_1^\top]\big\|\to0\quad\text{and}\quad\sqrt n\,\big\|\E[\hat L_1Y_1]-\E[L_1Y_1]\big\|\to0 ,
\end{equation}
then the same holds with $\beta^*_N$ in place of $\widetilde\beta_N$.
\end{corollary}

The distinction between the two targets is the one found in \citet{li2025assumption}. The noisy target $\widetilde\beta_N$ absorbs every centering, including the mean of $\zeta_i\zeta_i^\top$, which is the second moment of the covariate's network noise and does not vanish at the $\sqrt n$ scale without further sparsity; the latent target $\beta^*_N$ is the coefficient one would compute from the noiseless covariates, coincides with the true coefficients when the model is linear in $(X,\ell)$, and requires \eqref{eq:bias-cond}. In the examples below, \eqref{eq:bias-cond} is where the sparsity gap between the two targets appears, and it is the subsampled version of the gap between Theorems~1 and~2 of that paper; the bias correction proposed there, replacing $\frac{1}{|S|} \sum_{i \in S}\hat \ell_i\hat \ell_i^\top$ by the leading term of its structural decomposition, applies here unchanged and we do not repeat it.

\begin{example}[Normalized degree]\label{ex:degree}
For simplicity, suppose that $w_N \leq K$ and let $\hat \ell_i=d_i/((N-1)\rho_N)$, so $\ell=g$ and $\zeta_i=\hat \ell_i-g(\xi_i)$. Given $\xi$, $\zeta_i$ is a centered sum of $N-1$ independent Bernoulli terms with $\Var(\zeta_i\mid\xi)\le g(\xi_i)/((N-1)\rho_N)$, so $\E\zeta_1^2\le\E g/(N\rho_N)\to0$, which is (C1). For (C2), the terms $y^\circ_i\zeta_i$ are conditionally uncorrelated given $\xi$ except through the edges within $S$, so
\[
\Var\left(\sqrt n\,\frac{1}{|S|}\sum_{i \in S}y^\circ_i\zeta_i\right)\le\frac{C\,\E[(y^\circ)^2g]}{N\rho_N}+\frac{Cn\,\E[(y^\circ)^2w_N]}{N^2\rho_N}\to0 ,
\]
and likewise for $L_i\zeta_i^\top$ with $\E[\|L\|^2g]$, which is finite since $X$ and $\ell=g$ are bounded. For the first part of (C3), Cauchy--Schwarz and Proposition~\ref{prop:svar} give
\[
\sqrt n\,\E\left|\frac{1}{|S|}\sum_{i \in S}e_i\zeta_i \right|\le\sqrt n\,(\E e_1^2)^{1/2}(\E\zeta_1^2)^{1/2}\le C\Big(\frac{n\log N}{N\rho_N}\Big)^{1/2}\Big(\frac{\E g}{N\rho_N}\Big)^{1/2}=C\,\frac{\sqrt{n\log N}}{N\rho_N},
\]
and for the second part, $\zeta_i^2-\E[\zeta_i^2\mid\xi_i]$ has conditional variance of order $g(\xi_i)^2/(N\rho_N)^2$ and the same uncorrelatedness, so $\sqrt n\,\frac{1}{|S|}\sum_{i \in S}(\zeta_i^2-\E\zeta_1^2)=O_P\big((N\rho_N)^{-1}\big)+O_P\big(\sqrt n/(N\rho_N)\cdot N^{-1/2}\big)$ under $\E g^2<\infty$. Hence (C1)--(C3) hold whenever
\begin{equation}\label{eq:degree-sparsity}
N\rho_N\gg\sqrt{n\log N},
\end{equation}
and Corollary~\ref{cor:netcov} applies to $\widetilde\beta_N$. For $\beta^*_N$, condition \eqref{eq:bias-cond} fails unless the mean of $\zeta_1^2$, which equals $\E g/((N-1)\rho_N)$ up to a vanishing correction, is $o(n^{-1/2})$: the latent target requires
\begin{equation}\label{eq:degree-sparsity-clean}
N\rho_N\gg\sqrt n ,
\end{equation}
together with the corresponding condition on $\E[e_1\zeta_1]$, which is implied by \eqref{eq:degree-sparsity}. When $n=N$, \eqref{eq:degree-sparsity-clean} is $\lambda_N=\omega(\sqrt N)$, the condition of Theorem~1 of \citet{li2025assumption} for acyclic subgraphs, and \eqref{eq:degree-sparsity} is, up to the logarithm, the $\lambda_N=\omega(1)$ of its Theorem~2 strengthened by the present sparsity regime.
\end{example}

\begin{example}[Spectral embedding]\label{ex:spectral}
Suppose the graphon is fixed, bounded, and of finite rank $d$ as in the sparse generalized random dot product graph (GRDPG) model \citep{rubin2022statistical}, and let $\hat \ell_i$ be the $i$th row of the sparsity-normalized adjacency spectral embedding. The embedding is identified up to an orthogonal matrix $Q_N$, and the decomposition \eqref{eq:cov-decomp} holds for $Q_N\hat \ell_i$ with $\ell$ the latent position, which is bounded because $w$ is. Under $N\rho_N\gg\log^{2c}N$, the two-to-infinity bound and row-wise expansion of \citet{athreya2018statistical,rubin2022statistical} give $\max_i\|\zeta_i\|=O_P\big(\sqrt{\log N/(N\rho_N)}\big)$ and $\zeta_i=(A-\E[A\mid\xi])_{i\cdot}U\Lambda^{-1/2}+O_P\big(\log^cN/(N\rho_N)\big)$ uniformly in $i$, so that given $\xi$ the leading part of $\zeta_i$ is a centered linear function of row $i$ of $A-\E[A\mid\xi]$ and rows for distinct sampled nodes share only the edges within $S$. The verification of (C1)--(C3) is then as in Example~\ref{ex:degree}, with the second-order term handled by $\sqrt n\log^cN/(N\rho_N)\to0$; all three hold under
\begin{equation}\label{eq:spectral-sparsity}
N\rho_N\gg\sqrt n\,\log^cN ,
\end{equation}
and Corollary~\ref{cor:netcov} applies to the rotated covariates. Since the least-squares map is equivariant under invertible linear transformations of the covariates, this gives $\sqrt n\,(M_N\hat\beta-\widetilde\beta_N)\Rightarrow N(0,\Sigma_\beta)$ with $M_N=\mathrm{diag}(I_p,Q_N)$, as in Theorem~5 of \citet{li2025assumption}; coefficients on conventional covariates, and rotation-invariant functionals of the embedding coefficients such as $\|\hat\beta_\ell\|^2$ under $\beta^*_\ell=0$, are unaffected by $Q_N$. The condition \eqref{eq:spectral-sparsity} is the downsampled version of the $\lambda_N=\omega(\sqrt N\log^{2c}N)$ of that theorem. For the latent target $\beta^*_N$ the bias condition \eqref{eq:bias-cond} requires in addition $\sqrt n\,\E\|\zeta_1\|^2\to0$, which the row-wise expansion gives under $N\rho_N\gg\sqrt n$, so for embeddings the two targets separate only by the logarithm.
\end{example}

\section{Proof outline}\label{sec:proof-outline}
The main difficulty is to average over the graph rather than condition on it, for a statistic in which every node depends on every other through the resolvent of a random matrix with random degrees in its denominators, in a regime where those degrees are of order $\log N$. By Remark~\ref{rem:two-stage}, what has to be established is the limit of the marginal variance and of the variance of the population mean. Weak-dependence conditions in network distance are not the right tool here. Once the
graph is averaged over, the nodes are exchangeable: every pair $(Y_i,Y_j)$ has the same joint law, so the dependence between two nodes is weak but the same for every pair,
and network distance, a feature of the realized graph, plays no role.

The argument we give is organized around one decomposition and two ideas. The decomposition replaces each edge by its conditional probability given the latent positions, splitting $\bar Y_S-\mu_N$ into a signal component driven by the positions alone and a network component driven by the edges given the positions. The first idea is that the signal component is a smooth functional of the empirical distribution of latent positions, and that its first-order behavior can be captured by an exact identity rather than an expansion, from which the influence function and the limiting variance are read off. The second idea is that the network component is best studied not through walks but through the effect of a single node's edges on a linear functional of the response: a leave-one-node-out identity expresses that effect through a dual vector, an Efron--Stein inequality over vertex blocks sums it over nodes, and reversibility of the walk turns the resulting bounds into integrated moments of the graphon. In both halves, the only terms that cost more than fixed-order moments are walks that return to a node, and these are isolated and paid for once; this is why the conditions do not grow with walk length and why hubs can be allowed.

\subsection{The signal: three operators and the influence function}\label{sec:operators}

This section proves the signal half of Theorem~\ref{thm:clt}. It identifies the influence function of Section~\ref{sec:variance} from an exact identity, with no expansion in walk length and no differentiability, and reduces the signal component of $\bar Y_S$ to a sum of independent terms plus a remainder. The remainder bounds are stated here and proved in Appendix~\ref{app:remainder}.

\subsubsection{Splitting the response}\label{sec:split}

Given the latent positions, each edge indicator $A_{ij}$ has conditional mean $p_{ij}:=\rho_Nw_N(\xi_i,\xi_j)$. Let $\hat P$ be the transition matrix of the walk on the expected-weight graph, $\hat P_{ij}=p_{ij}/\sum_{l\ne i}p_{il}$, and let
\begin{equation}\label{eq:yW}
Y^W:=(I-\Gamma\hat P)^{-1}Z
\end{equation}
be the response that \eqref{eq:SAR} would produce if every edge were replaced by its conditional probability. Then
\begin{equation}\label{eq:split}
\bar Y_S-\mu_N=\underbrace{\big( \bar{Y}_S^W-\E\bar{Y}_S^W\big)}_{\text{signal}}+\underbrace{\big(\bar V_S-\E\bar V_S\big)}_{\text{network}},\quad\ \bar{Y}_S^W:=\frac{1}{|S|}\sum_{i \in S}Y^W_i,\quad V:=Y-Y^W .
\end{equation}
The signal component is a function of the latent positions alone; the network component is the effect of the edge noise given the positions. This section treats the first, Section~\ref{sec:fluctuation} the second. The decomposition is exact, and both components are centered, so $\mu_N=\E\bar Y^W_S+\E\bar V_S$ with no approximation. For statistics such as subgraph counts, with the latent positions treated as random, the edge noise is of lower order than the fluctuation of the positions under mild sparsity, and one might expect the same here. Two things are different. The edge noise enters \eqref{eq:SAR} through the degree normalization, nonlinearly and in every denominator, so its negligibility is not automatic and is the subject of Section~\ref{sec:fluctuation}. And $\bar Y^W_S$ is not an average of independent terms: each $Y^W_i$ is the response of node $i$ in a population whose type distribution is the empirical measure of the positions, so $\bar Y^W_S$ is a perturbation of the mean-field quantity $\mu^\circ_N$. Identifying the first-order part of that perturbation is the subject of this section.

\subsubsection{Three operators}\label{sec:three-ops}

Write $\hat g_i:=(N-1)^{-1}\sum_{j\ne i}w_N(\xi_i,\xi_j)$ for the expected-degree ratio of node $i$, so that $\hat P_{ij}=w_N(\xi_i,\xi_j)/((N-1)\hat g_i)$. The entries of $Y^W$ are then the values, at the sample points, of a functional of the empirical distribution of latent positions. For a CDF $F$ define
\begin{gather*}
g_F(x):=\int w_N(x,v)\,dF(v),\\
(T_Ff)(x):=\frac{1}{g_F(x)}\int w_N(x,v)f(v)\,dF(v),\\
y_F^\circ:=(I-\gamma T_F)^{-1} z .
\end{gather*}
Under (M2), $g_F\ge c_0$ and $T_F$ is a Markov operator, so $y^\circ_F$ is well defined under (M1) by the same Neumann argument as \eqref{eq:neumann}; $y^\circ_F(x)$ is the response of a node of type $x$ in a population with type distribution $F$. At $F=F^\circ$, the uniform law, these are the population objects of Section~\ref{sec:variance}: $g_{F^\circ}=g$, $T_{F^\circ}=T$, and $y^\circ_{F^\circ}=y^\circ$ is the mean-field response \eqref{eq:mean-field}. At $F=F_N:=N^{-1}\sum_l\mathbf{1}(\xi_l \leq \cdot )$, the empirical distribution of latent positions, $(T_{F_N}f)(\xi_i)$ is the expected-weight transition from node $i$ including the self-term $l=i$, which is $0$ since $w(x,x) = 0$; thus,
\begin{equation}\label{eq:link}
Y^W_i=y^\circ_{F_N}(\xi_i).
\end{equation}

Before proceeding, note that viewing a matrix built from a kernel evaluated at random positions as a perturbation of an integral operator goes back to at least \citet{koltchinskii2000random}. \citet{vonLuxburg2008consistency} treat the degree-normalized case for weighted similarity graphs, which corresponds to the expected-weight operator $\hat P$ here, and \citet{tang2013universally,rohe2011spectral} the binary case for dense latent position graphs, in each case proving convergence of spectra and eigenvectors, whereas \citet{lunde2023validity} addresses convergence of eigenvalues for the random walk integral operator for sparse graphons. Our use of the viewpoint is different in aim: we need the first-order fluctuation of the resolvent applied to a fixed function, as a functional of the empirical distribution of latent positions, at the $\sqrt n$ scale, in a sparse regime, and with an exact remainder.

The argument rests on three operators on functions of type, all built from $w_N$ and $\gamma$ alone.
\begin{enumerate}
\item The Markov operator $T$, with transition density $p(v\mid x)=w_N(x,v)/g(x)$ and stationary measure $\pi(dx)=g(x)\,dx$. It is reversible with respect to $\pi$ and contracts every $L^r(\pi)$, $r\ge1$.
\item The resolvent $\cR:=(I-T_\gamma)^{-1}$, $T_\gamma:=\gamma T$, with $\|\cR\|_{L^r(\pi)\to L^r(\pi)}\le(1-\bar\gamma)^{-1}$, together with its $L^2(dx)$-adjoint $\cR^*=(I-T_\gamma^*)^{-1}$, where $(T_\gamma^*h)(u)=\int h(x)\gamma(x)p(u\mid x)\,dx$. The mean-field response is $y^\circ=\cR z$, and the occupation density $q:=\cR^*\one$ gives the duality $\int \cR f\,dx=\int qf\,dx$ for integrable $f$ (Lemma~\ref{lem:duality}).
\item The perturbation operator, for a signed measure $\nu$ on $[0,1]$,
\[
(\dot T_\nu f)(x):=\int p(v\mid x)\big(f(v)-(Tf)(x)\big)\,\nu(dv),
\]
linear in $\nu$, and vanishing when $\nu$ is Lebesgue measure because $\int p(v\mid x)\,dv=1$.
\end{enumerate}
The first is the walk in the latent space, and is what carries the intuition of Section~\ref{sec:model}. The second is the walk summed over all lengths. The third is new, and it is the operator that a first-order Taylor expansion of $1/\hat g_i$ around its conditional expectation would try to approximate; the next lemma shows that no expansion is needed.

\subsubsection{An exact first-order identity}\label{sec:identity}

\begin{lemma}[Exact expansion]\label{lem:identity}
Assume (M1)--(M2) and $y^\circ\in L^1(\pi)$. For every distribution function $F$ on $[0,1]$ for which $y^\circ_F$ is defined,
\begin{equation}\label{eq:Texact}
T_Ff-Tf=\frac{g}{g_F}\,\dot T_{F-F^\circ}f ,
\end{equation}
and, exactly,
\begin{equation}\label{eq:yexact}
y^\circ_F-y^\circ=\cR\,\gamma\,\dot T_{F-F^\circ}\,y^\circ+\mathcal Q_F ,
\end{equation}
where
\[
\mathcal Q_F=\cR\,\gamma\Big(\frac{g}{g_F}-1\Big)\dot T_{F-F^\circ}\,y^\circ+(I-\gamma T_F)^{-1}\gamma\,(T_F-T)\,\cR\,\gamma\,\frac{g}{g_F}\,\dot T_{F-F^\circ}\,y^\circ .
\]
\end{lemma}

For each fixed $x$ and $f$, both $g/g_F-1=(g-g_F)/g_F$ and $(T_F-T)f=(g/g_F)\,\dot T_{F-F^\circ}f$
are linear in $F-F^\circ$ up to the factor $1/g_F$, which satisfies $g_F\ge c_0$ under
(M2); pointwise,
\[
\Big|\frac{g}{g_F}-1\Big|\le\frac{|g_F-g|}{c_0},\qquad
\big|(T_F-T)f\big|\le\frac{g}{c_0}\,\big|\dot T_{F-F^\circ}f\big| .
\]
Since each term of $\mathcal Q_F$ contains two such factors, $\mathcal Q_F$ may be
thought of as a quadratic remainder in $F-F^\circ$, and the first term of
\eqref{eq:yexact} as the linear part of $F\mapsto y^\circ_F$ at $F^\circ$, obtained
without differentiating anything. Making this quantitative would require operator
norms, since in $\mathcal Q_F$ the operator $T_F-T$ acts on a function that itself
depends on $F$; we instead bound the remainder through its matrix form below:
\begin{equation}\label{eq:matrix-identity}
Y^W-Y^\circ-\Phi=(I-\Gamma\hat P)^{-1}\Gamma\eta ,
\end{equation}
where $Y^\circ=(y^\circ(\xi_1), \ldots y^\circ(\xi_N) )$, $\Phi$ is defined in \eqref{eq:per-node} below, and $\eta$ is an explicit per-node error vector that is quadratic in the empirical fluctuations (Lemma~\ref{lem:A-matrix}); its three parts are described after
Proposition~\ref{prop:remainder}.

\begin{remark}[Walk-level reading of the first-order term]\label{rem:taylor}
In the length-one term of \eqref{eq:neumann} on the expected-weight graph, the
first-order change has two parts: the empirical neighborhood average in the
numerator fluctuates around $(Ty^\circ)(\xi_i)$, and the normalization $1/\hat g_i$
fluctuates around $1/g(\xi_i)$. Expanding the latter to first order and projecting
both onto a single other node $v$ gives exactly
$\gamma(x)p(v\mid x)\big(y^\circ(v)-(Ty^\circ)(x)\big)/(N-1)$: the direct term from the
numerator and the dilution term from the denominator. Collecting the first-order terms over all walk
lengths reproduces $\cR\gamma\dot T_{F_N-F^\circ}y^\circ$ term by term.
Lemma~\ref{lem:identity} is this computation done exactly, with the remainder
$\mathcal Q_F$ in closed form.
\end{remark}

\subsubsection{The influence kernel and the influence function}\label{sec:influence}

Taking $F-F^\circ=\mathbf{1}(v \leq \cdot)-F^\circ$ in $\dot T$ gives $(\dot T_{\mathbf{1}(v \leq \cdot)-F^\circ}f)(x)=p(v\mid x)\big(f(v)-(Tf)(x)\big)$, so the influence kernel \eqref{eq:influence} is $f_v=\gamma\,\dot T_{\mathbf{1}(v \leq \cdot)-F^\circ}\,y^\circ$, and $\phi(\cdot,v)=\cR f_v$ is the first-order change in the response of every type when an infinitesimal mass of type $v$ is added to the population. By linearity of $\dot T_\nu$ in $\nu$,
\begin{equation}\label{eq:linear-term}
\cR \,\gamma\,\dot T_{F_N-F^\circ}\,y^\circ=\frac1N\sum_{l=1}^N\phi(\cdot,\xi_l) ,
\end{equation}
which is the linear term of \eqref{eq:yexact} at the empirical law.

The definition $H_{\rm rest}(v)=\int_0^1\phi(x,v)\,dx=\int_0^1(\cR f_v)(x)\,dx$ suggests that
evaluating the influence function requires a separate resolvent solve for every type
$v$. It does not: moving the integral across the resolvent, with $\cR^*$ the adjoint
with respect to Lebesgue measure,
\[
H_{\rm rest}(v)=\int q(x)\,f_v(x)\,dx,\qquad q:=\cR^*\one ,
\]
so a single solve, for $q$, suffices, after which $H_{\rm rest}(v)$ is one integral of
the explicit kernel $f_v$ for each $v$. This identity is not used in the proofs. The
same duality appears in the sample: $q=\cR^*\one$ plays the role in the population that
$\lambda=R_P^\top\one_{\tilde S}$ plays in the network analysis of
Section~\ref{sec:fluctuation}. We include a short proof of the identity below.

\begin{lemma}[Duality]\label{lem:duality}
Under (M1)--(M2), for integrable $f$, $\int \cR f\,dx=\int q\,f\,dx$ with $q=\cR^*\one$.
\end{lemma}

\begin{proof}
Expand both resolvents in Neumann series and match terms: $\int T_\gamma^kf\,dx=\int(T_\gamma^*)^k\one\cdot f\,dx$. Absolute convergence follows from contraction in $L^1(\pi)$ and $dx\le d\pi/c_0$.
\end{proof}

\begin{proposition}[Influence function]\label{prop:H}
Under (M1)--(M2),
\begin{gather*}
H_{\mathrm{rest}}(v)=\int_0^1\phi(x,v)\,dx=\int_0^1q(u)\,\gamma(u)\,p(v\mid u)\big(y^\circ(v)-(Ty^\circ)(u)\big)\,du,\\
\int_0^1H_{\mathrm{rest}}(v)\,dv=0 .
\end{gather*}
\end{proposition}

\begin{proof}
Lemma~\ref{lem:duality} with $f=f_v$; the mean vanishes because $\int p(v\mid u)(y^\circ(v)-Ty^\circ(u))\,dv=0$ for every $u$.
\end{proof}

\subsection{Reduction to a sum of independent terms}\label{sec:reduction}

Evaluating \eqref{eq:yexact} at the sample points and inserting \eqref{eq:link},
\begin{equation}\label{eq:per-node}
Y^W_i=y^\circ(\xi_i)+\Phi_i+R_i,\quad\Phi_i:=\frac1{N-1}\sum_{l\ne i}\phi(\xi_i,\xi_l),
\end{equation}
where $R_i:=Y^W_i-y^\circ(\xi_i)-\Phi_i$ is the remainder, controlled through \eqref{eq:matrix-identity}. Averaging over $S$,
\[
\bar{Y}_S^W-\mu^\circ_N=\frac{1}{|S|}\sum_{i \in S}\big(y^\circ(\xi_i)-\mu^\circ_N\big)+\frac{1}{|S|}\sum_{i \in S}\Phi_i+\frac{1}{|S|}\sum_{i \in S}R_i .
\]
The middle term is a generalized $U$-statistic of order two in the latent positions. Since $\int\phi(x,v)\,dv=0$ for almost every $x$ (because $\int f_v\,dv=0$), its Hoeffding projection onto single coordinates involves only the second argument, and the projection of the pair $(i,l)$ onto $\xi_l$ is $H_{\mathrm{rest}}(\xi_l)$. Carrying out the projection gives the following.

\begin{proposition}[Linear statistic and exact variance]\label{prop:linear}
Let
\[
L_S:=\frac1n\sum_{i\in S}\big(H_{\mathrm{root}}(\xi_i)-\mu^\circ_N\big)+\frac1N\sum_{l=1}^NH_{\mathrm{rest}}(\xi_l),\qquad H_{\mathrm{root}}=y^\circ .
\]
Then $L_S$ is a sum of independent centered terms, and for every $n\le N$, $n\Var(L_S)=\sigma^2_{n,N}$ as defined in \eqref{eq:sigma}. Moreover,
\[
\frac{1}{|S|}\sum_{i \in S}\Phi_i=L_S-\frac{1}{|S|}\sum_{i \in S}\big(y^\circ(\xi_i)-\mu^\circ_N\big)+D_S+E_S ,
\]
where $D_S$ is a degenerate $U$-statistic with $\E D_S^2\le2\|\phi\|^2_{L^2}/(n(N-1))$,
and $E_S$ accounts for the fact that the exact projection of
$\avg_{i\in S}\Phi_i$ onto $\xi_l$ has weight $1/(N-1)$ or $(n-1)/(n(N-1))$, according
to whether $l\notin S$ or $l\in S$, while $L_S$ uses $1/N$ for every $l$; it satisfies
$\E E_S^2\le C(N^{-3}+n^{-1}N^{-2})\Var H_{\rm rest}$.
\end{proposition}

What remains is the remainder. The matrix identity \eqref{eq:matrix-identity} makes it a resolvent applied to the per-node error $\eta$, and the resolvent is a contraction in the right norm: with $\hat\pi_i=\hat g_i/\sum_l\hat g_l$, the matrix $\hat P$ is $\hat\pi$-reversible, so $\|(I-\Gamma\hat P)^{-1}\Gamma\|_{\ell^2(\hat\pi)}\le\bar\gamma/(1-\bar\gamma)$\footnote{Here $\|v\|^2_{\ell^2(\hat\pi)}=\sum_i\hat\pi_iv_i^2$, and
$\|Q\|_{\ell^2(\hat\pi)}$ is the operator norm induced by it. Weighted norms of this
kind are used throughout the appendix; see Appendix~\ref{app:notation} for further details.}, and $\hat g_i\ge c_0$ transfers the bound from the whole population to any subset $S$ at no cost in $n$:
\begin{equation}\label{eq:transfer}
\sqrt n\,\E\Big|\frac{1}{|S|}\sum_{i \in S}R_i\Big|\le\frac{\bar\gamma}{(1-\bar\gamma)\sqrt{c_0}}\Big(N\,\E\big[\hat g_1\eta_1^2\big]\Big)^{1/2}\qquad\text{uniformly in }S .
\end{equation}
The whole of the remainder analysis is therefore a bound on the single weighted second moment $\E[\hat g_1\eta_1^2]$, which must be $o(1/N)$.

\begin{proposition}[Signal remainder]\label{prop:remainder}
Assume (M1)--(M4) and either setting (B) or setting (U). Then $\sqrt n\,\E|\frac{1}{|S|}\sum_{i \in S}R_i|\le Cr_N$ uniformly in $S$, with $r_N^2:=N\,\E[\hat g_1\eta_1^2]$ and
\[
r_N^2\le\begin{cases}
CB^2\Big(\dfrac{\mathsf M_N(1+\rho_N^{-1})}{N}+\dfrac{1+\mathsf M_N}{(N\rho_N)^2}\Big)&\text{setting (B)},\\[10pt]
C\Big(\dfrac{\|w_N\|_\infty^2}{N}+\dfrac{\|w_N\|_\infty^3}{N^2}\Big)&\text{setting (U)},
\end{cases}
\]
so that $r_N\to0$ under (B3), respectively (U1). Moreover $\|\phi\|_{L^2}\le C$ and $\|H_{\mathrm{rest}}\|_{L^{2+\delta}}\le C$ in both settings, and the Lyapunov condition holds for $L_S$ under (B1)--(B2), respectively (U2)--(U3).
\end{proposition}

The proof is in Appendix~\ref{app:remainder}. The error $\eta_i$ has three parts. The first is $(g(\xi_i)/\hat g_i-1)$ times a centered neighborhood average: a product of two mean-zero averages over the other nodes, whose second moment is of order $1/N^2$ times integrated second moments of $p(\cdot\mid x)$, which are controlled by $\mathsf M_N$. The second is a degenerate $U$-statistic of order two, whose kernel involves $T\phi$, and its second moment is again $1/N^2$ times a norm of $\phi$. The third is a diagonal term involving $\phi(\cdot,x)$ evaluated back at the node itself, the return of a node's influence to its origin, which is where the sup norm of the degree profile enters. In setting (B) this costs one factor $\rho_N^{-1}$ inside $r_N^2$; in setting (U), $\|w_N\|_\infty$ enters every one of these terms, which is how (U1) is used in this section.

\begin{proof}[Conclusion of the signal half of Theorem~\ref{thm:clt}]
By \eqref{eq:per-node} and Proposition~\ref{prop:linear},
\[
\sqrt n\,(\bar{Y}_S^W-\E\bar{Y}_S^W)=\sqrt n\,L_S+\sqrt n\,(D_S+E_S)+\sqrt n\,\left(\frac{1}{|S|}\sum_{i \in S}R_i-\E\frac{1}{|S|}\sum_{i \in S}R_i\right) .
\]
The second group is $o_P(1)$ by the variance bounds of Proposition~\ref{prop:linear} and $\|\phi\|_{L^2}\le C$; the third is $o_P(1)$ by Proposition~\ref{prop:remainder}. Since $L_S$ is a sum of independent terms with $n\Var L_S=\sigma^2_{n,N}$, bounded below by (M5), and the Lyapunov condition holds, $\sqrt n\,L_S/\sigma_{n,N}\Rightarrow N(0,1)$ by the Lindeberg--Feller theorem for triangular arrays.
\end{proof}

Centering at $\mu^\circ_N$ or at $\E\bar{Y}_S^W$ is equivalent at the $\sqrt n$ scale, since $\E\bar{Y}_S^W-\mu^\circ_N=\E\frac{1}{|S|}\sum_{i \in S}R_i=O(r_N/\sqrt n)$; the bias $\mu_N-\mu^\circ_N$, of order $\log N/(N\rho_N)$ up to lower-order terms, comes mainly from $\E\bar V_S$ and is the subject of Section~\ref{sec:bias}.

\subsection{The network fluctuation}\label{sec:fluctuation}

This section proves the network half of Theorem~\ref{thm:clt}: given the latent positions, replacing each edge by its conditional probability changes the mean response by $o_P(n^{-1/2})$.

\begin{theorem}[Network fluctuation]\label{thm:fluct}
Assume (M1)--(M4). In setting (B), for every $S$ with $|S|=n\le N$,
\[
n\Var(\bar V_S)\le C\,B^2\,(1+\mathsf M_N)\,\frac{\log N}{N\rho_N},\qquad C=C(\bar\gamma,c_0),
\]
so that $n\Var(\bar V_S)\to0$ under (B3). In setting (U), $\sup_{|S|\le N}n\Var(\bar V_S)\le C\big((N\rho_N)^{-1}+\|w_N\|_\infty^4/N^2\big)\to0$.
\end{theorem}

The rate in setting (B) is the one that fixes the sparsity regime: it is $o(1)$ exactly when the average degree is of larger order than $\log N$ (times $\mathsf M_N$), which (B3) guarantees. The proof has three ingredients, which we describe in turn; details are in Appendix~\ref{app:fluctuation}.

\subsubsection{Vertex-block Efron--Stein}
The randomness of the graph given the positions comes from the edge uniforms
$(U_{ij})_{i<j}$, which enter $V$ in a complicated way through the resolvent.
Rather than treat edges one at a time, we group them by vertex: the block of
node $j$ is $\mathcal B_j=(\xi_j,(U_{jl})_{l<j})$, the position of $j$ together with
the uniforms of its edges to earlier nodes. The blocks are independent, and $(\xi,A)$,
hence $V$, is a measurable function of them, so the Efron--Stein inequality applies
with one term per node. For independent random elements $\mathcal B_1,\dots,\mathcal B_N$ and square-integrable
$f=f(\mathcal B_1,\dots,\mathcal B_N)$, the Efron--Stein inequality
\citep{efron1981jackknife,steele1986efron} states that
\[
\Var f\le\frac12\sum_{j=1}^N\E\big(f-f^{(j)}\big)^2,
\]
where $f^{(j)}$ is computed with $\mathcal B_j$ replaced by an independent copy. Redrawing $\mathcal B_j$ changes only edges incident to $j$:
those to earlier nodes through $U_{jl}$, and those to later nodes through $\xi_j$.
A $\sigma$-field based on these random variables was first used by \citet{Borgs2008Convergent} and the network Efron--Stein inequality was proved by \citet{pmlr-v119-lin20c} using the same construction.

Bounding the effect of one row and column of a random matrix on a functional of its
resolvent is also a classical strategy in random matrix theory, for example in the
martingale arguments for linear spectral statistics of \citet{BaiSilverstein2004}; the
difference here is that the matrix is a degree-normalized walk on a sparse graph,
which is what the duality and reversibility arguments below address.

By exchangeability it suffices to resample a single node, which we take to be node
$1$, and
\begin{equation}\label{eq:ES}
n\Var(\bar V_S)\le4\,\E V_1^2+\frac{4N}{n}\sup_{\tilde S}\E\,G^2 ,
\end{equation}
where $G:=\sum_{i\in\tilde S}\big(V_i-V^{(0)}_i\big)$ is the effect of node $1$ on the
other sampled nodes: the difference between their network components on the full
graph and on the graph with node $1$ removed. Here $\tilde S$ ranges over the possible
positions of the sample relative to the resampled node, as explained below. The first
term is the variance of a single node's fluctuation; the second is the influence of
one node's edges on the other sampled nodes, which must be $o(n/N)$.

\subsubsection{Leave-one-node-out identity and dual vectors}

Let $P^{(0)}$ be the transition matrix of the graph with node $1$ removed (with row
and column $1$ set to zero), $R^{(0)}=(I-\Gamma P^{(0)})^{-1}$, and $Y^{(0)}=R^{(0)}Z$
the corresponding response. Subtracting $Y^{(0)}=Z+\Gamma P^{(0)}Y^{(0)}$ from
$Y=Z+\Gamma PY$, and adding and subtracting $\Gamma P^{(0)}Y$, gives the exact
identity
\[
Y-Y^{(0)}=R^{(0)}\Gamma\big(P-P^{(0)}\big)Y .
\]
The perturbation $P-P^{(0)}$ is supported on node $1$ and its neighbors: adding
node $1$ back changes the degrees of its neighbors and adds a new term to their
neighborhood averages. The reduced resolvent appears on the left and the full
response on the right, so the effect of node $1$'s presence on everyone else is
carried exactly by the full responses; no correction is needed. The effect of this
perturbation on a linear functional of the response, such as a sum over the sample,
is read off by duality.

\emph{Duality.} For a fixed graph, the response is linear in the inputs, $Y=R_PZ$
with $R_P=(I-\Gamma P)^{-1}$, so a sum over a set $\tilde S$ of nodes is also linear in
the inputs:
\[
\sum_{i\in\tilde S}Y_i=\sum_m\lambda_mZ_m,\qquad\lambda_m:=\sum_{i\in\tilde S}(R_P)_{im}.
\]
The weight $\lambda_m$ is the total effect of node $m$'s input on the sum: the change
in $\sum_{i\in\tilde S}Y_i$ when $Z_m$ increases by one and all other inputs are held
fixed. In walk terms, it collects the discounted walks that start in $\tilde S$ and
end at $m$. The vector $\lambda=R_P^\top\one_{\tilde S}$ is a column sum of the
resolvent, the solution of the dual, or adjoint, problem. A familiar statistical
analogue is the hat matrix of a linear regression, whose columns give the influence
of each observation on the fitted values.

\emph{The dual vector.} Because the walk is reversible with respect to the degrees,
$P^\top=DPD^{-1}$; equivalently $d_i(P^t)_{im}=d_m(P^t)_{mi}$, so that, weighted by the
degree of the starting node, a $t$-step path from $i$ to $m$ is as likely as the same
path from $m$ to $i$. The discount factors do not affect this: since $\Gamma$ and $D$
are both diagonal, they commute, and $\big((\Gamma P)^t\big)^\top=D(P\Gamma)^tD^{-1}$.
Summing over $t$,
\[
\lambda=Db,\qquad b:=(I-P\Gamma)^{-1}D^{-1}\one_{\tilde S}=\sum_{t\ge0}(P\Gamma)^tD^{-1}\one_{\tilde S}.
\]
The entry $b_m$ is a walk expectation: over walks started at $m$, the product of the
coefficients $\gamma$ along the walk times $\one\{\text{endpoint}\in\tilde S\}/d_{\text{endpoint}}$.
We call $b$ the \emph{dual vector} and work with it rather than with $\lambda$,
because $\lambda_m=d_mb_m$ grows with the degree of $m$, whereas $b$ is uniformly of
order $1/(N\rho_N)$ and contracts under the walk in every degree-weighted norm. The
degree factor is restored only at the end, where degree-biased averaging accounts for
it.

\emph{Which graph.} Efron--Stein asks how much the sample sum changes when node
$1$'s block is redrawn. The graph with node $1$ removed is the same before and after
redrawing, since it is a function of everything except node $1$'s block, so it serves
as a common reference point. Splitting off node $1$, and writing
$\tilde S:=S\setminus\{1\}$,
\[
\sum_{i\in S}V_i=\underbrace{\sum_{i\in\tilde S}V^{(0)}_i}_{\text{reduced graph}}
+\underbrace{\one\{1\in S\}\,V_1}_{\text{node $1$'s own term}}
+\underbrace{G}_{\text{node $1$'s effect on others}} ,
\]
with $G$ as in \eqref{eq:ES}. Redrawing node $1$'s block leaves the first term
unchanged, so if primes denote values after redrawing,
\[
\sum_{i\in S}V_i-\sum_{i\in S}V_i'=\one\{1\in S\}\,(V_1-V_1')+(G-G') .
\]
Node $1$'s own term is bounded separately, through $\E V_1^2$. Since $G'$ has the same
law as $G$, $\E(G-G')^2\le4\E G^2$, so it suffices to bound the size of $G$ for a
single draw: the effect of adding node $1$ back to the reduced graph. (For a general
resampled node $j$, relabeling $j$ as node $1$ turns $S\setminus\{j\}$ into such a set
$\tilde S$, which is why \eqref{eq:ES} takes a supremum over $\tilde S$.) For the same
reason, the dual vector is computed on the reduced graph, with $P^{(0)}$ and
$D^{(0)}$; this is the role of the reduced resolvent in the identity above. It is then
fixed when the block is redrawn, so node $1$'s edges enter $G$ only linearly, weighted
by the dual vector at node $1$'s neighbors. Since these edges are independent
Bernoulli variables given the types, the fluctuations of $G$ can be computed
directly.

\emph{The network component.} The quantity whose variance we bound is not the
response itself but the network component $V=Y-Y^W$ of \eqref{eq:split}. Both
responses are linear in the inputs, so the argument applies to each: the actual
response has dual vector $b$, and the expected-weight response has its own dual
vector $\hat b$, built in the same way from $\hat P^{(0)}$. The contribution $G$ of
node $1$ is therefore a sum over node $1$'s neighbors $m$ of dual weights times what
$m$ receives from node $1$. We write $G=G_1+G_2+G_3$, where $G_1$ is the fluctuation
of node $1$'s edges around their conditional probabilities, $G_2$ is the mismatch
between the two dual vectors $b$ and $\hat b$, and $G_3$ is the feedback of node $1$'s
presence through the rest of the network.

With heavy-tailed inputs, one further step is needed. The factors multiplying the
dual weights now include the responses and the network noise at node $1$'s neighbors,
which are unbounded and depend on the same graph. Bounding them uniformly would cost
the largest input in the population. Instead, we compute moments. Conditional on
everything except node $1$'s edges, each quantity is a polynomial in independent coin
flips $A_{m1}$ with known success probabilities. Its second or fourth moment is
computed by expanding the polynomial and sorting the terms by which neighbors
coincide, exactly as one computes $\E\big(\sum_iX_i\big)^4$ for independent summands,
or the variance of a $U$-statistic. Each coincidence pattern produces a population
moment of the graphon and the input, after averaging over node $1$'s type, and the
functionals in (U4) are precisely the moments that appear. In setting (B),
these factors are bounded, and this computation is not needed.

\section{Simulation study}
\label{sec:simulations}

We examine three questions numerically. First, does the ordinary sample variance
$s_S^2$ estimate the marginal variance of the responses across sampling fractions,
autoregressive configurations and sparsity levels? Second, when is the marginal
variance the right variance for the subsample mean, so that the $t$-interval with
standard error $s_S/\sqrt n$ is calibrated? Third, how does the discrepancy at larger
sampling fractions depend on the sign of the autoregressive coefficients? The theory
predicts that for negative coefficients the influence term can reduce the variance of
the mean below the marginal variance, so the naive interval should then be
conservative rather than anticonservative.

\subsection{Design}

We construct a model in which node-level inputs, autoregressive coefficients and the network are generated jointly
through shared latent positions. For
each node, we draw independent $\xi_i,V_{X,i},V_{\gamma,i},V_{\varepsilon,i}\sim
\mathrm{Unif}(0,1)$, independently across nodes, and independent edge uniforms
$U_{ij}\sim\mathrm{Unif}(0,1)$, $i<j$. Writing $t(u)=2u-1$, we set
\[
X_i=\beta_0+\beta_Ut(\xi_i)+s_Xt(V_{X,i}),\qquad
\gamma_i=a+b_Ut(\xi_i)+b_Vt(V_{\gamma,i}),\qquad
\varepsilon_i=s_\varepsilon(\xi_i)\,t(V_{\varepsilon,i}),
\]
with $\beta_0=\beta_U=1$, $s_X=0.5$ and $s_\varepsilon(u)=0.5\{1+0.5\,t(u)\}$, and
$Z_i=X_i+\varepsilon_i$. Thus $\E[\varepsilon_i\mid\xi_i,X_i,\gamma_i]=0$, the
residual variance $s_\varepsilon(\xi_i)^2/3$ varies with the latent position, and
$\Cov(X_i,\gamma_i)=\beta_Ub_U/3$, so the input and the susceptibility to peers are
correlated through the latent position. All configurations below
satisfy $|a|+|b_U|+|b_V|\le0.8$, so $|\gamma_i|\le0.8$.

The graphon is $w(u,v)=(0.5+uv)/0.75$, which integrates to one and is bounded below
by $2/3$. Responses are generated on the
full network of $N=5000$ nodes, and the subsample $S$ is drawn uniformly without
replacement, independently of the population.

We consider five autoregressive configurations:
\begin{itemize}
\item zero: $(a,b_U,b_V)=(0,0,0)$, so $\gamma_i=0$;
\item positive low: $(a,b_U,b_V)=(0.2,0.05,0.05)$, so $\gamma_i\in[0.1,0.3]$;
\item negative low: $(a,b_U,b_V)=(-0.2,-0.05,-0.05)$, so $\gamma_i\in[-0.3,-0.1]$;
\item positive high: $(a,b_U,b_V)=(0.6,0.1,0.1)$, so $\gamma_i\in[0.4,0.8]$;
\item negative high: $(a,b_U,b_V)=(-0.6,-0.1,-0.1)$, so $\gamma_i\in[-0.8,-0.4]$.
\end{itemize}
Subsample sizes are $n=\lfloor N^\alpha\rfloor$ with $\alpha\in\{0.5,0.7,0.9\}$, giving
$n=70$, $388$ and $2133$, and $n=cN$ with $c\in\{0.1,0.5,1\}$, where $c=1$ is full
observation.

For each configuration, we estimate the finite-population reference quantities
$\widetilde v_N:=\Var(Y^{(N)}_1)$, $\widetilde V_{n,N}:=n\Var(\bar Y_S)$ and
$\widetilde V_{N,N}:=N\Var(\bar Y_N)$ from $B_{\rm ref}=1000$ independent populations.
These are Monte Carlo estimates of finite-$N$ variances: $\widetilde v_N$ is the
marginal variance of a response, which the sample variance estimates, and
$\widetilde V_{n,N}$ is the variance that appears in the central limit theorem at
sampling fraction $n/N$. The sample variance $s_S^2$ and the confidence intervals are
computed from a separate set of $B=1000$ independent populations, with one subsample
drawn from each.

\subsection{Variance across sampling fractions}

Table~\ref{tab:sim-fraction} reports the results at $\rho_N=\log N/\sqrt N$, an average
degree of about $600$. The sample variance tracks the marginal variance closely in every
configuration: its Monte Carlo mean is within $0.012$ of $\widetilde v_N$ at every
subsample size, and its standard deviation decreases from $0.066$--$0.112$ at
$n=70$ to $0.007$--$0.015$ at $n=5000$. The table shows the two extreme subsample sizes; the
others are similar.

\begin{table}[ht]
\centering
\small
\setlength{\tabcolsep}{3.5pt}
\renewcommand{\arraystretch}{1.2}
\begin{tabular}{lccccccccc}
\toprule
& & \multicolumn{2}{c}{$s_S^2$, mean (SD)} & \multicolumn{6}{c}{$\widetilde V_{n,N}$ at $n$ (and $n/N$)}\\
\cmidrule(lr){3-4}\cmidrule(lr){5-10}
$\gamma$ setting & $\widetilde v_N$ & $n=70$ & $n=5000$
& \shortstack{$70$\\$(0.014)$} & \shortstack{$388$\\$(0.078)$} & \shortstack{$500$\\$(0.1)$}
& \shortstack{$2133$\\$(0.43)$} & \shortstack{$2500$\\$(0.5)$} & \shortstack{$5000$\\$(1)$}\\
\midrule
Zero          & 0.507 & 0.515 (0.074) & 0.507 (0.008) & 0.507 & 0.506 & 0.505 & 0.499 & 0.498 & 0.488 \\
Positive low  & 0.569 & 0.577 (0.082) & 0.569 (0.010) & 0.573 & 0.591 & 0.597 & 0.688 & 0.709 & 0.849 \\
Negative low  & 0.469 & 0.476 (0.070) & 0.469 (0.008) & 0.467 & 0.457 & 0.454 & 0.404 & 0.393 & 0.317 \\
Positive high & 0.809 & 0.820 (0.112) & 0.808 (0.015) & 0.865 & 1.123 & 1.214 & 2.540 & 2.839 & 4.869 \\
Negative high & 0.440 & 0.447 (0.066) & 0.440 (0.007) & 0.436 & 0.419 & 0.413 & 0.325 & 0.306 & 0.171 \\
\bottomrule
\end{tabular}
\caption{Marginal variance $\widetilde v_N$, sample variance $s_S^2$, and the
variance $\widetilde V_{n,N}$ of the scaled subsample mean, at $N=5000$ and
$\rho_N=\log N/\sqrt N$. The Monte Carlo standard error of $\widetilde V_{N,N}$ ranges
from $0.008$ (negative high) to $0.21$ (positive high).}
\label{tab:sim-fraction}
\end{table}

Whether the marginal variance is also the right variance for the subsample mean
depends on the sampling fraction. For the smallest subsample, $n/N\approx0.014$,
$\widetilde V_{n,N}$ is within about $7\%$ of $\widetilde v_N$ in every configuration.
As $n/N$ grows, the two separate, and the separation depends on the strength and sign
of the autoregressive coefficients. Without peer effects, $\widetilde V_{n,N}$ stays
close to $\widetilde v_N$ throughout, as it should. With positive coefficients it
grows: in the positive-high setting it is already about $40\%$ above the marginal
variance at $n/N\approx0.08$, and six times larger at $n=N$. With negative
coefficients it falls below the marginal variance, to less than $40\%$ of it at $n=N$
in the negative-high setting. These finite-sample results are consistent with the
theory, in which the network enters the variance only through the influence term
weighted by $n/N$. They also show that how small $n/N$ must be for the network to be
negligible depends on how strong the peer effects are.

\subsection{Sparsity and coverage}

We next vary the sparsity, for the positive-high and negative-high configurations and
for a small and a substantial subsample, $n=70$ and $n=2500$. We consider
$\rho_N=\log N/N$ and $2\log N/N$, average degrees of about $8.5$ and $17$, at the
logarithmic boundary of the theory, and $\rho_N=\log N/\sqrt N$ from the previous
experiment for comparison.  For each
configuration, we also compute the coverage of the nominal $95\%$ $t$-interval
$\bar Y_S\pm t_{n-1,0.975}\,s_S/\sqrt n$ for $\mu_N=\E Y^{(N)}_1$, estimated by the
average of $\bar Y_N$ over the reference populations. The Monte Carlo standard error
of each coverage estimate is at most $0.015$.

\begin{table}[ht]
\centering
\small
\setlength{\tabcolsep}{4pt}
\renewcommand{\arraystretch}{1.2}
\begin{tabular}{lcccccccc}
\toprule
& & & \multicolumn{3}{c}{$n=70$} & \multicolumn{2}{c}{$n=2500$} & \\
\cmidrule(lr){4-6}\cmidrule(lr){7-8}
$\rho_N$ & \shortstack{avg.\\degree} & $\widetilde v_N$ & \shortstack{$s_S^2$\\mean (SD)} & $\widetilde V_{n,N}$ & coverage & $\widetilde V_{n,N}$ & coverage & $\widetilde V_{N,N}$\\
\midrule
\multicolumn{9}{l}{\emph{Positive high, $\gamma_i\in[0.4,0.8]$}}\\
$\log N/N$       & 8.5 & 0.939 & 0.938 (0.136) & 1.004 & 0.950 & 3.290 & 0.712 & 5.641\\
$2\log N/N$      & 17  & 0.865 & 0.864 (0.124) & 0.927 & 0.942 & 3.129 & 0.667 & 5.394\\
$\log N/\sqrt N$ & 602 & 0.809 & 0.820 (0.112) & 0.865 & 0.939 & 2.839 & 0.675 & 4.869\\
\midrule
\multicolumn{9}{l}{\emph{Negative high, $\gamma_i\in[-0.8,-0.4]$}}\\
$\log N/N$       & 8.5 & 0.507 & 0.507 (0.075) & 0.503 & 0.953 & 0.357 & 0.985 & 0.207\\
$2\log N/N$      & 17  & 0.470 & 0.473 (0.073) & 0.466 & 0.951 & 0.328 & 0.980 & 0.186\\
$\log N/\sqrt N$ & 602 & 0.440 & 0.447 (0.066) & 0.436 & 0.946 & 0.306 & 0.980 & 0.171\\
\bottomrule
\end{tabular}
\caption{Variances and coverage of nominal $95\%$ $t$-intervals across sparsity levels,
at $N=5000$. At $n=2500$, the mean of $s_S^2$ is within $0.003$ of $\widetilde v_N$ in
every row.}
\label{tab:sim-sparsity}
\end{table}

Table~\ref{tab:sim-sparsity} shows three things. First, the sample variance
estimates the marginal variance accurately at every sparsity level, including the
two near the boundary. The marginal variance itself increases as the graph becomes
sparser, by about $0.13$ at average degree $8.5$ in the positive-high setting and
$0.07$ in the negative-high setting, relative to the densest level; the increase roughly
halves when the average degree doubles. This is consistent with the network-noise
contribution of order $1/(N\rho_N)$ in Proposition~\ref{prop:svar}.

Second, for the small subsample, the $t$-interval is close to nominal at every
sparsity level, with coverage between $0.939$ and $0.953$ across all levels,
using only the observed responses. Third, for $n=2500$, the interval is
miscalibrated in the direction the theory predicts. In the positive-high setting, the
variance of the mean is about $3.5$ times the marginal variance, the interval is too
narrow, and the coverage is $0.667$--$0.712$. In the negative-high setting, the
variance of the mean is about $0.7$ times the marginal variance, the interval is too
wide, and the coverage is $0.980$--$0.985$. The pattern is the same at every sparsity
level.

Taken together, the simulations support the use of the ordinary $t$-interval when the
subsample is a small fraction of the population, without knowledge of the network or of the propagation mechanism, across the sparsity levels considered. At substantial
sampling fractions, the variance of the mean includes the influence term, which the
sample variance does not estimate, and valid inference then requires information about
the network.

\section{Conclusion}\label{sec:conclusion}
In this work, we have proved a central limit theorem for sample means, and smooth functions of vector means, of responses generated by a heterogeneous spatial autoregression on a sparse random graph, a class broad enough to include random coefficients, nonparametric inputs, unbounded graphons with hubs, and heavy-tailed inputs, in the sparsity regime in which degrees concentrate. The limiting variance interpolates exactly between the marginal variance of the mean-field response and the variance of a node's total influence on the population. The practical consequence is the one we regard as most important: when the observed nodes are a vanishing fraction of the population through which dependence propagates, the ordinary sample variance is consistent and the ordinary $t$-interval is valid, with no knowledge of the network, the graphon, or the strength of peer influence; and when they are not, subsampling restores validity at the cost of the rate. We do not believe this phenomenon is specific to this model. It arises because the randomness of the graph and of the latent positions is integrated over rather than conditioned on, so that the dependence a sampled node inherits from the unsampled population is of smaller order than its own variation, and we expect similar behaviour for other common network regression models.

The two-stage structure of Remark~\ref{rem:two-stage} identifies what drives the phenomenon,
and hence what must be verified in other models. Whenever the responses are exchangeable
once the network is averaged over, their marginal variance converges to a limit that does
not depend on the realized network, and the population mean satisfies a central limit
theorem at the $N^{-1/2}$ scale, inference for the superpopulation mean is network-free
when $n/N\to0$, and the network enters the variance only through the population term
otherwise. The spatial autoregression is a demanding test of these conditions, and we
expect them to be considerably easier to verify in models with more limited propagation,
such as local spillovers or exposure-based outcome models.

Several questions remain open. We do not know whether the logarithm in the sparsity condition (M4) can be removed, or whether the fourth-moment-type conditions on the degree profile can be reduced to the $2+\delta$ moments that a Gaussian limit would seem to require. The sampling designs covered are those independent of the network; we suspect that network dependence also diminishes even in certain cases involving non-representative samples if the sample is small relative to the superpopulation. A related question is when sample splitting, and inference after selection based on it,
remains valid for network-dependent data: our results suggest that two halves that are each
small fractions of the population behave approximately as independent samples. Each of these
would require ideas beyond those of this paper.

\section*{Acknowledgements}
Robert Lunde acknowledges support from NSF grant DMS-2515923. During the preparation of this work the authors used Claude Opus 5.5 and Claude Fable 5.1 (Anthropic). Building on the authors’ earlier analysis, Claude Opus 5.5 suggested the exact first-order identity and the vertex-block Efron--Stein approach, and the models drafted proofs and text from specifications provided by the authors and checked proofs. The authors reviewed, verified and revised all content, and take full responsibility for the results.

\bibliographystyle{plainnat}
\bibliography{references}

@article{BaiSilverstein2004,
  author    = {Zhidong Bai and Jack W. Silverstein},
  title     = {CLT for linear spectral statistics of large-dimensional sample covariance matrices},
  journal   = {The Annals of Probability},
  volume    = {32},
  number    = {1A},
  pages     = {553--605},
  year      = {2004},
  publisher = {Institute of Mathematical Statistics},
}

@article{steele1986efron,
  author    = {J. Michael Steele},
  title     = {An {Efron}-{Stein} Inequality for Nonsymmetric Statistics},
  journal   = {The Annals of Statistics},
  volume    = {14},
  number    = {2},
  pages     = {753--758},
  year      = {1986},
  publisher = {Institute of Mathematical Statistics},
}

@article{efron1981jackknife,
  title={The Jackknife Estimate of Variance},
  author={Efron, Bradley and Stein, Charles},
  journal={The Annals of Statistics},
  volume={9},
  number={3},
  pages={586--596},
  year={1981},
}

@article{tropp2012user,
  title={User-friendly tail bounds for sums of random matrices},
  author={Tropp, Joel A.},
  journal={Foundations of Computational Mathematics},
  volume={12},
  number={4},
  pages={389--434},
  year={2012},
  publisher={Springer}
}

@article{maurer2003bound,
  author  = {Maurer, Andreas},
  title   = {A Bound on the Deviation Probability for Sums of Non-Negative Random Variables},
  journal = {Journal of Inequalities in Pure and Applied Mathematics (JIPAM)},
  volume  = {4},
  number  = {1},
  year    = {2003},
}

@book{kallenberg2002foundations,
  author    = {Kallenberg, Olav},
  title     = {Foundations of modern probability},
  edition   = {Second},
  series    = {Probability and its Applications (New York)},
  publisher = {Springer-Verlag},
  address   = {New York},
  year      = {2002},
}

@InProceedings{pmlr-v119-lin20c,
  title = 	 {On the Theoretical Properties of the Network Jackknife},
  author =       {Lin, Qiaohui and Lunde, Robert and Sarkar, Purnamrita},
  booktitle = 	 {Proceedings of the 37th International Conference on Machine Learning},
  pages = 	 {6105--6115},
  year = 	 {2020},
  volume = 	 {119},
  series = 	 {Proceedings of Machine Learning Research},
}

@article{Borgs2008Convergent,
  title     = {Convergent sequences of dense graphs I: Subgraph frequencies, metric properties and testing},
  author    = {Borgs, Christian and Chayes, Jennifer T. and Lov{\'a}sz, L{\'a}szl{\'o} and S{\'o}s, Vera T. and Vesztergombi, Katalin},
  journal   = {Advances in Mathematics},
  volume    = {219},
  number    = {6},
  pages     = {1801--1851},
  year      = {2008},
  publisher = {Elsevier},
}

@article{lunde2023validity,
  title={On the Validity of Conformal Prediction for Network Data Under Non-Uniform Sampling},
  author={Lunde, Robert},
  journal={Bernoulli (forthcoming)},
  year={2026}
}

@article{rohe2011spectral,
  author    = {Karl Rohe and Sourav Chatterjee and Bin Yu},
  title     = {Spectral clustering and the high-dimensional stochastic blockmodel},
  journal   = {The Annals of Statistics},
  volume    = {39},
  number    = {4},
  pages     = {1878--1915},
  year      = {2011},
}

@article{tang2013universally,
  author    = {Tang, Minh and Sussman, Daniel L. and Priebe, Carey E.},
  title     = {Universally consistent vertex classification for latent positions graphs},
  journal   = {The Annals of Statistics},
  volume    = {41},
  number    = {3},
  pages     = {1406--1430},
  year      = {2013},
}

@article{vonLuxburg2008consistency,
  author    = {Ulrike von Luxburg and Mikhail Belkin and Olivier Bousquet},
  title     = {Consistency of spectral clustering},
  journal   = {The Annals of Statistics},
  volume    = {36},
  number    = {2},
  pages     = {555--586},
  year      = {2008},
  publisher = {Institute of Mathematical Statistics},
}

@article{koltchinskii2000random,
  author    = {Koltchinskii, Vladimir and Gin{\'e}, Evarist},
  title     = {Random matrix approximation of spectra of integral operators},
  journal   = {Bernoulli},
  volume    = {6},
  number    = {1},
  pages     = {113--167},
  year      = {2000},
}

@article{li2022random,
  title={Random graph asymptotics for treatment effect estimation under network interference},
  author={Li, Shuangning and Wager, Stefan},
  journal={The Annals of Statistics},
  volume={50},
  number={4},
  pages={2334--2358},
  year={2022},
}

@article{rubin2022statistical,
  title={A statistical interpretation of spectral embedding: the generalised random dot product graph},
  author={Rubin-Delanchy, Patrick and Cape, Joshua and Tang, Minh and Priebe, Carey E},
  journal={Journal of the Royal Statistical Society Series B: Statistical Methodology},
  volume={84},
  number={4},
  pages={1446--1473},
  year={2022},
  publisher={Oxford University Press}
}

@article{athreya2018statistical,
  title={Statistical inference on random dot product graphs: a survey},
  author = {Athreya, Avanti and Fishkind, Donniell E. and Tang, Minh and Priebe, Carey E. and Park, Youngser and Vogelstein, Joshua T. and Levin, Keith and Lyzinski, Vince and Qin, Yichen},
  journal={Journal of Machine Learning Research},
  volume={18},
  number={226},
  pages={1--92},
  year={2018}
}

@article{Bickel-Chen-Levina-method-of-moments,
  author = "Peter J. Bickel and Aiyou Chen and Elizaveta Levina",
  title = "The method of moments and degree distributions for network models",
  journal = "Annals of Statistics",
  volume = 39,
  year = 2011,
  pages = "38-59",
  }

@article{hastie1993varying,
  author    = {Hastie, Trevor and Tibshirani, Robert},
  title     = {Varying-Coefficient Models},
  journal   = {Journal of the Royal Statistical Society: Series B (Methodological)},
  volume    = {55},
  number    = {4},
  pages     = {757--779},
  year      = {1993},
  publisher = {Wiley},
}

@article{aldous-representation-array,
  author= "David J. Aldous",
  title = "Representations for partially exchangeable arrays of random variables",
  journal= "Journal of Multivariate Analysis",
  volume = 11,
  year= 1981,
  pages= "581-598"
  }

@misc{hoover-exchangeability,
  author = {Hoover, Douglas N.},
  title  = {Relations on probability spaces and arrays of random variables},
  note   = {Preprint, Institute for Advanced Study, Princeton, NJ},
  year   = {1979},
}

@article{friedkin1990social,
  title={Social influence and opinions},
  author={Friedkin, Noah E. and Johnsen, Eugene C.},
  journal={Journal of Mathematical Sociology},
  volume={15},
  number={3-4},
  pages={193--206},
  year={1990},
  publisher={Taylor \& Francis},
  }

@article{bollobas2011sparse,
  title={Sparse graphs: Metrics and random models},
  author={Bollob{\'a}s, B{\'e}la and Riordan, Oliver},
  journal={Random Structures \& Algorithms},
  volume={39},
  number={1},
  pages={1--38},
  year={2011},
  publisher={Wiley Online Library},
  }

@article{borgs-lp-part-one,
   author = {Christian Borgs and Jennifer T. Chayes and Henry Cohn and Yufei Zhao},
    title = {An {$L^p$} theory of sparse graph convergence {I}: Limits, sparse random graph models, and power law distributions},
  journal = {Transactions of the American Mathematical Society},
  volume= {372},
  year= {2019},
  pages={3019-3062}
     }

@article{Bickel-Chen-on-modularity,
  author  = {Peter J. Bickel and Aiyou Chen},
  title   = {A Nonparametric View of Network Models and {Newman}-{Girvan} and Other Modularities},
  journal = {Proceedings of the National Academy of Sciences (USA)},
  volume  = {106},
  pages   = {21068--21073},
  year    = {2009}
}

@article{diaconis1980finite,
  author = {Diaconis, P. and Freedman, D.},
  title = {Finite Exchangeable Sequences},
  journal = {The Annals of Probability},
  volume = {8},
  number = {4},
  pages = {745--764},
  year = {1980},
  publisher = {Institute of Mathematical Statistics},
 }

@article{shalizi2011homophily,
  title={Homophily and contagion are generically confounded in observational social network studies},
  author={Shalizi, Cosma Rohilla and Thomas, Andrew C},
  journal={Sociological Methods \& Research},
  volume={40},
  number={2},
  pages={211--239},
  year={2011},
}

@misc{lei2024leastsquares,
  title={Least Squares Inference for Data with Network Dependency}, 
  author={Jing Lei and Kehui Chen and Haeun Moon},
  year={2024},
  eprint={2404.01977},
  archivePrefix={arXiv},
  primaryClass={stat.ME},
  url={https://arxiv.org/abs/2404.01977}
}

@article{Lee03072021,
author = {Youjin Lee and Elizabeth L. Ogburn},
title = {Network Dependence Can Lead to Spurious Associations and Invalid Inference},
journal = {Journal of the American Statistical Association},
volume = {116},
number = {535},
pages = {1060--1074},
year = {2021},

}

@article{bramoulle2009identification,
  title   = {Identification of Peer Effects through Social Networks},
  author  = {Bramoull{\'e}, Yann and Djebbari, Habiba and Fortin, Bernard},
  journal = {Journal of Econometrics},
  volume  = {150},
  number  = {1},
  pages   = {41--55},
  year    = {2009}
}

@book{lovasz2012large,
  title     = {Large Networks and Graph Limits},
  author    = {Lov{\'a}sz, L{\'a}szl{\'o}},
  publisher = {American Mathematical Society},
  address   = {Providence, RI},
  year      = {2012}
}

@article{kojevnikov2021limit,
  title   = {Limit Theorems for Network Dependent Random Variables},
  author  = {Kojevnikov, Denis and Marmer, Vadim and Song, Kyungchul},
  journal = {Journal of Econometrics},
  volume  = {222},
  number  = {2},
  pages   = {882--908},
  year    = {2021}
}

@article{li2025assumption,
  title={Assumption-lean Inference for Network-linked Data},
  author={Li, Wei and Chakraborty, Nilanjan and Lunde, Robert},
  journal={arXiv preprint arXiv:2510.00287},
  year={2025}
}

@article{baldi1989normal,
  author  = {Baldi, Pierre and Rinott, Yosef},
  title   = {On normal approximations of distributions in terms of dependency graphs},
  journal = {The Annals of Probability},
  year    = {1989},
  volume  = {17},
  number  = {4},
  pages   = {1646--1650}
}

@article{chen2004normal,
  author  = {Chen, Louis H. Y. and Shao, Qi-Man},
  title   = {Normal approximation under local dependence},
  journal = {The Annals of Probability},
  year    = {2004},
  volume  = {32},
  number  = {3A},
  pages   = {1985--2028}
}

@article{jenish2009central,
  author  = {Jenish, Nazgul and Prucha, Ingmar R.},
  title   = {Central limit theorems and uniform laws of large numbers for arrays of random fields},
  journal = {Journal of Econometrics},
  year    = {2009},
  volume  = {150},
  number  = {1},
  pages   = {86--98}
}

@article{wu2005nonlinear,
  author  = {Wu, Wei Biao},
  title   = {Nonlinear system theory: Another look at dependence},
  journal = {Proceedings of the National Academy of Sciences},
  year    = {2005},
  volume  = {102},
  number  = {40},
  pages   = {14150--14154}
}

@article{xu2024functional,
  author  = {Wu, Zeqi and Jiang, Wen and Xu, Xingbai},
  title   = {Applications of functional dependence to spatial econometrics},
  journal = {Econometric Theory},
  volume = {41},
  number = {5},
  pages = {1044--1079},
  year    = {2025},
}

@unpublished{jiang2025limit,
  author = {Jiang, Wen and Wang, Yachen and Wu, Zeqi and Xu, Xingbai},
  title  = {Limit theorems for network data without metric structure},
  year   = {2025},
  note   = {arXiv:2511.17928}
}

@article{leung2026normal,
  author  = {Leung, Michael P. and Moon, Hyungsik Roger},
  title   = {Normal approximation in large network models},
  journal = {The Review of Economic Studies},
  year    = {2026},
  note    = {Advance article},
}

@article{ogburn2024causal,
  author  = {Ogburn, Elizabeth L. and Sofrygin, Oleg and D{\'\i}az, Iv{\'a}n and van der Laan, Mark J.},
  title   = {Causal inference for social network data},
  journal = {Journal of the American Statistical Association},
  year    = {2024},
  volume  = {119},
  number  = {545},
  pages   = {597--611},
}

@article{silverman1976limit,
  author  = {Silverman, Bernard W.},
  title   = {Limit theorems for dissociated random variables},
  journal = {Advances in Applied Probability},
  year    = {1976},
  volume  = {8},
  number  = {4},
  pages   = {806--819}
}

@article{eagleson1978limit,
  author  = {Eagleson, G. K. and Weber, N. C.},
  title   = {Limit theorems for weakly exchangeable arrays},
  journal = {Mathematical Proceedings of the Cambridge Philosophical Society},
  year    = {1978},
  volume  = {84},
  number  = {1},
  pages   = {123--130}
}

@article{davezies2021empirical,
  author  = {Davezies, Laurent and D'Haultf{\oe}uille, Xavier and Guyonvarch, Yannick},
  title   = {Empirical process results for exchangeable arrays},
  journal = {The Annals of Statistics},
  year    = {2021},
  volume  = {49},
  number  = {2},
  pages   = {845--862},
}

@article{menzel2021bootstrap,
  author  = {Menzel, Konrad},
  title   = {Bootstrap with cluster-dependence in two or more dimensions},
  journal = {Econometrica},
  year    = {2021},
  volume  = {89},
  number  = {5},
  pages   = {2143--2188}
}

@article{bhattacharya2015subsampling,
  author  = {Bhattacharyya, Sharmodeep and Bickel, Peter J.},
  title   = {Subsampling bootstrap of count features of networks},
  journal = {The Annals of Statistics},
  year    = {2015},
  volume  = {43},
  number  = {6},
  pages   = {2384--2411}
}

@article{lunde2023subsampling,
  author  = {Lunde, Robert and Sarkar, Purnamrita},
  title   = {Subsampling sparse graphons under minimal assumptions},
  journal = {Biometrika},
  year    = {2023},
  volume  = {110},
  number  = {1},
  pages   = {15--32}
}

@book{cliff1973spatial,
  author    = {Cliff, Andrew D. and Ord, J. Keith},
  title     = {Spatial Autocorrelation},
  publisher = {Pion},
  address   = {London},
  year      = {1973}
}

@article{kelejian1998generalized,
  author  = {Kelejian, Harry H. and Prucha, Ingmar R.},
  title   = {A generalized spatial two-stage least squares procedure for estimating a spatial autoregressive model with autoregressive disturbances},
  journal = {The Journal of Real Estate Finance and Economics},
  year    = {1998},
  volume  = {17},
  number  = {1},
  pages   = {99--121}
}

@article{kelejian1999generalized,
  author  = {Kelejian, Harry H. and Prucha, Ingmar R.},
  title   = {A generalized moments estimator for the autoregressive parameter in a spatial model},
  journal = {International Economic Review},
  year    = {1999},
  volume  = {40},
  number  = {2},
  pages   = {509--533}
}

@article{kelejian2010specification,
  author  = {Kelejian, Harry H. and Prucha, Ingmar R.},
  title   = {Specification and estimation of spatial autoregressive models with autoregressive and heteroskedastic disturbances},
  journal = {Journal of Econometrics},
  year    = {2010},
  volume  = {157},
  number  = {1},
  pages   = {53--67}
}

@article{lee2004asymptotic,
  author  = {Lee, Lung-Fei},
  title   = {Asymptotic distributions of quasi-maximum likelihood estimators for spatial autoregressive models},
  journal = {Econometrica},
  year    = {2004},
  volume  = {72},
  number  = {6},
  pages   = {1899--1925},
}

@article{manski1993identification,
  author  = {Manski, Charles F.},
  title   = {Identification of endogenous social effects: The reflection problem},
  journal = {The Review of Economic Studies},
  year    = {1993},
  volume  = {60},
  number  = {3},
  pages   = {531--542}
}

@article{lin2010identifying,
  author  = {Lin, Xu},
  title   = {Identifying peer effects in student academic achievement by spatial autoregressive models with group unobservables},
  journal = {Journal of Labor Economics},
  year    = {2010},
  volume  = {28},
  number  = {4},
  pages   = {825--860}
}

@article{aquaro2021estimation,
  author  = {Aquaro, Michele and Bailey, Natalia and Pesaran, M. Hashem},
  title   = {Estimation and inference for spatial models with heterogeneous coefficients: An application to {US} house prices},
  journal = {Journal of Applied Econometrics},
  year    = {2021},
  volume  = {36},
  number  = {1},
  pages   = {18--44},
  }

@article{johnsson2021estimation,
  author  = {Johnsson, Ida and Moon, Hyungsik Roger},
  title   = {Estimation of peer effects in endogenous social networks: Control function approach},
  journal = {The Review of Economics and Statistics},
  year    = {2021},
  volume  = {103},
  number  = {2},
  pages   = {328--345},
}

@article{auerbach2022identification,
  author  = {Auerbach, Eric},
  title   = {Identification and estimation of a partially linear regression model using network data},
  journal = {Econometrica},
  year    = {2022},
  volume  = {90},
  number  = {1},
  pages   = {347--365}
}

@article{boucher2022estimating,
title     = {Estimating peer effects using partial network data},
author    = {Vincent Boucher and Aristide Houndetoungan},
journal   = {Review of Economics and Statistics},
year      = {2025},
pages     = {1--48},
}

@article{lewbel2023social,
  author  = {Lewbel, Arthur and Qu, Xi and Tang, Xun},
  title   = {Social networks with unobserved links},
  journal = {Journal of Political Economy},
  year    = {2023},
  volume  = {131},
  number  = {4},
  pages   = {898--946}
}

@article{zhu2017network,
  author  = {Zhu, Xuening and Pan, Rui and Li, Guodong and Liu, Yuewen and Wang, Hansheng},
  title   = {Network vector autoregression},
  journal = {The Annals of Statistics},
  year    = {2017},
  volume  = {45},
  number  = {3},
  pages   = {1096--1123},
}

@article{zhu2019network,
  author  = {Zhu, Xuening and Wang, Weining and Wang, Hansheng and H{\"a}rdle, Wolfgang Karl},
  title   = {Network quantile autoregression},
  journal = {Journal of Econometrics},
  year    = {2019},
  volume  = {212},
  number  = {1},
  pages   = {345--358}
}

@article{zhu2020grouped,
  author  = {Zhu, Xuening and Pan, Rui},
  title   = {Grouped network vector autoregression},
  journal = {Statistica Sinica},
  year    = {2020},
  volume  = {30},
  number  = {3},
  pages   = {1437--1462}
}

@article{chen2023community,
  author  = {Chen, Elynn Y. and Fan, Jianqing and Zhu, Xuening},
  title   = {Community network auto-regression for high-dimensional time series},
  journal = {Journal of Econometrics},
  year    = {2023},
  volume  = {235},
  number  = {2},
  pages   = {1239--1256}
}

@article{armillotta2023nonlinear,
  author  = {Armillotta, Mirko and Fokianos, Konstantinos},
  title   = {Nonlinear network autoregression},
  journal = {The Annals of Statistics},
  year    = {2023},
  volume  = {51},
  number  = {6},
  pages   = {2526--2552},
}

@article{armillotta2024count,
  author  = {Armillotta, Mirko and Fokianos, Konstantinos},
  title   = {Count network autoregression},
  journal = {Journal of Time Series Analysis},
  year    = {2024},
  volume  = {45},
  number  = {4},
  pages   = {584--612}
}

@article{knight2020generalized,
  author  = {Knight, Marina and Leeming, Kathryn and Nason, Guy and Nunes, Matthew},
  title   = {Generalized network autoregressive processes and the {GNAR} package},
  journal = {Journal of Statistical Software},
  year    = {2020},
  volume  = {96},
  number  = {5},
  pages   = {1--36}
}

@unpublished{wu2023random,
  author = {Wu, Weichi and Leng, Chenlei},
  title  = {A random graph-based autoregressive model for networked time series},
  year   = {2023},
  note   = {arXiv:2309.08488}
}

@book{kish1965survey,
  author    = {Kish, Leslie},
  title     = {Survey Sampling},
  publisher = {Wiley},
  address   = {New York},
  year      = {1965},
}

@article{abadie2020sampling,
  author  = {Abadie, Alberto and Athey, Susan and Imbens, Guido W. and Wooldridge, Jeffrey M.},
  title   = {Sampling-based versus design-based uncertainty in regression analysis},
  journal = {Econometrica},
  volume  = {88},
  number  = {1},
  pages   = {265--296},
  year    = {2020},
}

@article{blom1976some,
  author  = {Blom, Gunnar},
  title   = {Some properties of incomplete {U}-statistics},
  journal = {Biometrika},
  volume  = {63},
  number  = {3},
  pages   = {573--580},
  year    = {1976},
}

@article{janson1984asymptotic,
  author  = {Janson, Svante},
  title   = {The asymptotic distributions of incomplete {U}-statistics},
  journal = {Zeitschrift f{\"u}r Wahrscheinlichkeitstheorie und verwandte Gebiete},
  volume  = {66},
  pages   = {495--505},
  year    = {1984},
}

@article{chen2019randomized,
  author  = {Chen, Xiaohui and Kato, Kengo},
  title   = {Randomized incomplete {U}-statistics in high dimensions},
  journal = {The Annals of Statistics},
  volume  = {47},
  number  = {6},
  pages   = {3127--3156},
  year    = {2019},
}

\section*{Appendix}
The appendix is organized as follows. Appendix~\ref{app:notation} collects notation
and elementary facts used throughout. Appendix~\ref{app:identity} proves the exact
first-order identity, its matrix form, and the reduction of the signal component to a
sum of independent terms. Appendix~\ref{app:remainder} establishes the population
bounds, the Lyapunov condition and the properties of the limiting variance, and bounds
the signal remainder in settings (B) and (U). Appendix~\ref{app:tools} contains the
tools for the network component: concentration and operator bounds, moments of the
graph noise, products of Bernoulli sums, the leave-one-node-out identity, dual vectors
and degree-biased averaging. Appendix~\ref{app:fluctuation} proves
Theorem~\ref{thm:fluct}, first for bounded inputs and then for heavy-tailed inputs.
Appendix~\ref{app:svar} treats variance estimation, Appendix~\ref{app:weighted}
weighted means, and Appendix~\ref{app:bias} the finite-population bias. A reader
interested only in the network-free $t$-interval for $n/N\to0$ needs
Appendices~\ref{app:identity}, \ref{app:remainder} and~\ref{app:svar}, the setting (B)
part of Appendix~\ref{app:fluctuation}, and the parts of Appendix~\ref{app:tools} it
cites.

\appendix

\numberwithin{equation}{section}
\section{Notation and conventions}\label{app:notation}

Table~\ref{tab:not-global} lists the notation used throughout; later sections have their own tables for local notation.

\begin{table}[htbp]
\centering\small
\begin{tabular}{p{0.22\linewidth}p{0.55\linewidth}p{0.15\linewidth}}
\toprule
Symbol & Meaning & Defined in\\
\midrule
\multicolumn{3}{l}{\emph{Population (functions of type)}}\\
$w=w_N$, $\rho=\rho_N$ & graphon and sparsity level & Def.~\ref{def:graphon}\\
$W_N$ & largest edge weight $\|w_N\|_\infty$ & App.~\ref{app:notation}\\
$g(x)$ & expected-degree profile: a node of type $x$ has expected degree $(N-1)\rho g(x)$ & \S\ref{sec:conditions}\\
$p(v\mid x)$ & transition density of the walk on types & \S\ref{sec:conditions}\\
$T$, $\pi$ & one-step average of the walk; its stationary law $\pi(dx)=g(x)\,dx$ & App.~\ref{app:notation}\\
$T_\gamma$, $\cR$ & discounted step $\gamma T$; resolvent $(I-T_\gamma)^{-1}$, summing all walk lengths & App.~\ref{app:notation}\\
$y^\circ$ & mean-field response $\cR z$: response of a type in an infinite population & \eqref{eq:mean-field}\\
$\phi(x,v)$, $\phi_v$ & total influence of a type-$v$ node on a type-$x$ node; $\phi_v=\phi(\cdot,v)$ is this influence as a function of the receiving type, so $T\phi_v$ averages it over one step of the walk & \S\ref{sec:variance}\\
$\Phi_s(u)$ & size $\|\phi(u,\cdot)\|_{L^s}$ of the influence received by type $u$ from a random influencer (compare $\phi_v$, which fixes the influencer) & App.~\ref{app:notation}\\
$H_{\rm root}$, $H_{\rm rest}$ & a node's contribution through its own response; through everyone else's & \eqref{eq:H}\\
$\kappa_2$, $\mathsf M_N$ & $\kappa_2(x)=\int w(x,v)^2dv$; mixed norm $\int\kappa_2^2$ (always $\ge1$) & App.~\ref{app:notation}\\
$\Theta_s$ & $s$-th moment of the one-step influence kernel & App.~\ref{app:notation}\\
$e_x(v)$ & deviation $y^\circ(v)-(Ty^\circ)(x)$ of a type-$v$ neighbor's response from the average around $x$ & App.~\ref{app:notation}\\
$f_v$ & one-step influence kernel $\gamma(x)p(v\mid x)e_x(v)$: direct effect plus dilution & App.~\ref{app:notation}\\

\midrule
\multicolumn{3}{l}{\emph{Sample (functions of the types $\xi_1,\dots,\xi_N$)}}\\
$\hat g_i$ & empirical expected-degree ratio $(N-1)^{-1}\sum_{j\ne i}w(\xi_i,\xi_j)$ & App.~\ref{app:notation}\\
$\hat P$, $Y^W$ & walk on the expected-weight graph; response with edges replaced by their probabilities & \eqref{eq:yW}\\
$\Phi_i$, $R_i$ & linear term of node $i$; remainder $Y^W_i-y^\circ(\xi_i)-\Phi_i$ & \eqref{eq:per-node}\\
$\eta$ & per-node error: degree-ratio error times a centered neighborhood average, minus a degenerate $U$-statistic and a diagonal term; it generates the remainder, $R=(I-\Gamma\hat P)^{-1}\Gamma\eta$ & Lem.~\ref{lem:A-matrix}\\
$r_N$ & remainder rate, $r_N^2=N\,\E[\hat g_1\eta_1^2]$ & Lem.~\ref{lem:A-transfer}\\
\midrule
\multicolumn{3}{l}{\emph{Network (given the types)}}\\
$p_{ij}$, $\mu_i$, $d_i$ & edge probability $\rho w(\xi_i,\xi_j)$; expected degree; degree & App.~\ref{app:notation}\\
$P$, $X$, $\cE$ & actual walk $D^{-1}A$; centered edges $A-\E[A\mid\F]$; perturbation $P-\hat P$ & App.~\ref{app:notation}\\
$V$ & network component $Y-Y^W$ & \eqref{eq:split}\\
$\dmin$, $\varepsilon$ & lower bound $c_0(N-2)\rho$ on all expected degrees; concentration level $C_\varepsilon\sqrt{\log N/\dmin}$ & App.~\ref{app:notation}\\
$\F$ & $\sigma$-field of the types & App.~\ref{app:notation} \\
\bottomrule
\end{tabular}
\caption{Notation used throughout the appendix.}
\label{tab:not-global}
\end{table}

\paragraph*{Array and constants} We write $w=w_N$, $\rho=\rho_N$ and $W_N:=\|w_N\|_\infty$; in setting (B), $W_N$ may be as large as $1/(2\rho_N)$. Unless stated otherwise, $C$ denotes a constant depending only on $(\bar\gamma,c_0,\delta,K)$, and may change from line to line. In setting (B) we keep the dependence on $B$ explicit and write $B_y:=B/(1-\bar\gamma)$. We set
\[
\dmin:=c_0(N-2)\rho_N,\qquad \varepsilon:=C_\varepsilon\sqrt{\log N/\dmin},
\]
with an absolute constant $C_\varepsilon$ chosen in Lemma~\ref{lem:A-conc}. Under (M4) with $C_*\ge16C_\varepsilon^2$ we have $\varepsilon\le1/4$, and $\dmin\to\infty$. Note that $\dmin$ is a lower bound on the expected degrees $\mu_i$, which is why it carries the letter $\mu$. Two further conventions: $\E[\,\cdot\mid x]$ abbreviates $\E[\,\cdot\mid\xi_1=x]$, and ``on $\Omega$'' for an event $\Omega$ means that the stated inequality holds pointwise on $\Omega$; indicators of named events are written $\one_\Omega$, and bare $\one\{\cdot\}$ is reserved for indicator functions of explicit conditions.

\paragraph*{Population objects} Let $g(x)=\int w(x,v)\,dv$ denote the expected-degree
profile, so that a node of type $x$ has expected degree $(N-1)\rho_Ng(x)$, and let
$p(v\mid x)=w(x,v)/g(x)$ denote the transition density of a random walk
$(\zeta_k)_{k\ge0}$ on the latent space.  We use the terms ``type'' and ``latent position'' interchangeably. The Markov operator
$(Tf)(x)=\int p(v\mid x)f(v)\,dv=\E[f(\zeta_1)\mid\zeta_0=x]$ averages $f$ over one
step of the walk. Since $g(x)p(v\mid x)=w(x,v)$ is symmetric in $(x,v)$, the walk is
reversible with respect to $\pi(dx)=g(x)\,dx$, which is therefore a stationary
probability measure ($\int g=1$). Let $T_\gamma=\gamma T$, that is,
$(T_\gamma f)(x)=\gamma(x)(Tf)(x)$. When $\gamma\ge0$, $T_\gamma$ is the kernel of the
walk terminated at each step with probability $1-\gamma(\zeta_k)$; in general it is a
signed kernel with $\|T_\gamma\|_{\infty\to\infty}\le\bar\gamma<1$. The resolvent
$\cR:=(I-T_\gamma)^{-1}=\sum_{k\ge0}T_\gamma^k$ sums these discounted one-step
averages over all walk lengths, so that
$(\cR f)(x)=\E\big[\sum_{k\ge0}\prod_{j<k}\gamma(\zeta_j)\,f(\zeta_k)\bigm|\zeta_0=x\big]$. Shifting the index in the Neumann series gives $\cR=I+T_\gamma\cR$: the walk either
stops at once or takes one step and then continues as before.

We also use the $\pi$-preserving Markov operator $\bar T:=(1-\bar\gamma)\sum_{k\ge0}\bar\gamma^kT^k$, the multi-step densities $p^{(k)}(v\mid x)$, and
\[
K_t(x):=\int p(v\mid x)^t\,dv,\qquad \kappa_2(x):=\int w(x,v)^2\,dv=g(x)^2K_2(x),\qquad \mathsf M_N=\int\kappa_2^2 .
\]
Note that, to lighten notation, all integrals are over $[0,1]$ (or $[0,1]^2$ for $\iint$) with
respect to Lebesgue measure unless another measure, such as $\pi$, is written
explicitly. When the variable of integration is clear from context, we omit the
differential, writing for example $\int g$ for $\int_0^1 g(x)\,dx$ and $\iint w^2$ for
$\int_0^1\int_0^1 w(x,v)^2\,dx\,dv$.

The mean-field response is $y^\circ=\cR z$. We write $e_x(v):=y^\circ(v)-(Ty^\circ)(x)$, $f_v(x)=\gamma(x)p(v\mid x)e_x(v)$,  $\phi_v:=\phi(\cdot,v)=\cR f_v$, $\Phi_s(u):=\|\phi(u,\cdot)\|_{L^s(dv)}$, and
\[
\Theta_s:=\iint w(x,v)\,p(v\mid x)^{s-1}\,|e_x(v)|^s\,dv\,dx,\qquad s\ge1 .
\]
By $\cR=I+T_\gamma\cR$, the mean-field response and the influence profiles satisfy
the fixed-point equations
\[
y^\circ=z+\gamma\,Ty^\circ,\qquad \phi_v=f_v+\gamma\,T\phi_v ,
\]
which are used repeatedly below; the second holds pointwise wherever the Neumann
series converges absolutely. The sample resolvents satisfy the analogous equations
$Y^W=Z+\Gamma\hat PY^W$ and $Y=Z+\Gamma PY$.

For setting (U) we use the functionals of (U4):
\begin{gather*}
W_2(x,z)=\int w(x,v)w(v,z)\,dv,\qquad \mathfrak D=\iint g(x)^2W_2(x,z)^2,\\
\mathfrak h(z)=\int\kappa_2(v)w(v,z)\,dv,\qquad \mathcal H=\int\mathfrak h^2,\qquad \mathfrak G=\int g^3(Tg)^2 .
\end{gather*}
The functionals in (U4) are degree moments and two- and three-point coincidence
functionals: weights of short paths in the latent space that backtrack
($\mathsf M_N$, $\mathcal H$), share a neighbor ($\mathfrak D$), or pass through
high-degree neighbors ($\mathfrak G$). They arise in the network fluctuation, when the
effect of a single node's edges is averaged over that node's latent type, and each
controls one kind of coincidence among that node's neighbors.

They are easiest to interpret in a product graphon. Suppose each node has a
\emph{sociability} $a(\xi_i)\ge a_0>0$, and that connection probabilities are
proportional to the product of sociabilities,
\[
w(x,v)=\frac{a(x)\,a(v)}{\big(\int a\big)^2},
\]
as in Chung--Lu and degree-corrected block models. Then a node's expected degree is
proportional to its sociability, $g(x)=a(x)/\int a$, and each functional in (U4) is
bounded by a constant multiple of a moment of $a(\xi)$ of order at most four. For
example, $\mathsf M_N=\E[a(\xi)^4]\,\big(\E[a(\xi)^2]\big)^2/\big(\E\,a(\xi)\big)^8$.
Thus (U4) holds whenever $\E[a(\xi)^4]<\infty$: hubs are allowed, provided the degree
distribution has four moments.

\paragraph*{Sample objects} $\hat g_i=(N-1)^{-1}\sum_{j\ne i}w(\xi_i,\xi_j)$, $g_i:=g(\xi_i)$, and
\[
\hat P_{ij}=\frac{w(\xi_i,\xi_j)}{\sum_{l\ne i}w(\xi_i,\xi_l)}=\frac{g_i}{\hat g_i}\,\frac{p(\xi_j \mid \xi_i)}{N-1}.
\]
The expected-weight response is $Y^W=(I-\Gamma\hat P)^{-1}Z$, and $\Phi_i=(N-1)^{-1}\sum_{l\ne i}\phi(\xi_i,\xi_l)$.

\paragraph*{Network objects} The edge probabilities are $p_{ij}=\rho w(\xi_i,\xi_j)$, and $\mu_i=\sum_{j\ne i}p_{ij}=\rho(N-1)\hat g_i$. Let $\F:=\sigma(\xi_1,\dots,\xi_N)$ be the $\sigma$-field generated by the types. Given
$\F$, the edges $A_{ij}$, $i<j$, are independent $\mathrm{Bern}(p_{ij})$, since they are functions of the independent edge uniforms.   We write $d_i$ for the degrees, $P=D^{-1}A$ (zero rows at isolated nodes), $X:=A-\E[A\mid\F]$, $\cE:=P-\hat P$, $R_P:=(I-\Gamma P)^{-1}$, $M:=(I-\Gamma\hat P)^{-1}$, $V:=Y-Y^W$, and $\bar V_S=\avg_{i\in S}V_i$, where $\avg_{i\in S}$ is shorthand for average over elements in $S$.

\paragraph*{Weighted norms} For positive weights $\nu=(\nu_i)$ on an index set and
$s\ge1$, we write $\|v\|_{\ell^s(\nu)}:=\big(\sum_i\nu_i|v_i|^s\big)^{1/s}$ and
$\|v\|_\infty:=\max_i|v_i|$. The weights used are degrees $d$, expected degrees $\mu$,
expected-degree ratios $\hat g$, and their normalized versions such as
$\hat\pi=\hat g/\sum_l\hat g_l$. For functions of type, $\|f\|_{L^s(\pi)}:=\big(\int|f|^s\,d\pi\big)^{1/s}$,
and $\|f\|_{L^s}$ without a measure refers to Lebesgue measure. For a matrix or
integral operator $Q$, $\|Q\|_{\ell^s(\nu)}$ and $\|Q\|_{L^s(\pi)}$ denote the induced
operator norms, $\sup_{v\ne0}\|Qv\|/\|v\|$ in the corresponding norm, and
$\|Q\|_{\infty\to\infty}$ the operator norm on $\ell^\infty$, the maximum absolute row
sum. Rescaling the weights by a constant does not change an induced operator norm, so
$\ell^2(\hat g)$ and $\ell^2(\hat\pi)$, for instance, give the same operator norms.

\paragraph*{Elementary facts used throughout}
\begin{enumerate}[label=(F\arabic*),leftmargin=*]
\item\label{F:M3} By (M3), $w\le1/(2\rho)$, so $g\le1/(2\rho)$, $p(v\mid x)\le(2\rho c_0)^{-1}$ and $K_2\le(2\rho c_0)^{-1}$.
\item\label{F:M2} By (M2), $g\ge c_0$ pointwise and $\hat g_i\ge c_0$ almost surely. Hence $\|f\|_{L^r(dx)}\le c_0^{-1/r}\|f\|_{L^r(\pi)}$.
\item\label{F:Mbig} $g^2\le\kappa_2$ by Jensen. Hence $\E g^2\le\int\kappa_2\le\mathsf M_N^{1/2}$, $\E g^4\le\mathsf M_N$, and $\mathsf M_N\ge(\int\kappa_2)^2\ge(\int g^2)^2\ge1$.
\item\label{F:contr} $T$ contracts every $L^r(\pi)$, $r\ge1$, by Jensen and stationarity. Hence $\|\cR\|_{L^r(\pi)}\le(1-\bar\gamma)^{-1}$ and $\|y^\circ\|_{L^r(\pi)}\le(1-\bar\gamma)^{-1}\|Z\|_{L^r(\pi)}$. In setting (B), $|y^\circ|,|Ty^\circ|\le B_y$ and $|e_x(v)|\le2B_y$.
\item\label{F:diag} With the convention $w(x,x)=0$, $(T_{F_N}f)(\xi_i)=(\hat Pf(\xi_\cdot))_i$ exactly, so $Y^W_i=y^\circ_{F_N}(\xi_i)$ (equation \eqref{eq:link}); no self-loop correction arises.
\item\label{F:ops} \emph{Norms of walk operators.} Let $Q$ be a row-stochastic matrix
that is reversible with respect to positive weights $\nu$, meaning $\nu_iQ_{ij}=\nu_jQ_{ji}$
(for example $P$ with $\nu=d$, or $\hat P$ with $\nu=\mu$ or $\nu=\hat g$), and let
$|\gamma|\le\bar\gamma$. Then, for every vector $v$:
\begin{itemize}
\item[(i)] $|(Q\Gamma v)_i|\le\bar\gamma\|v\|_\infty$, by H\"older with exponents
$(1,\infty)$ row by row, since each row of $Q$ is a probability vector;
\item[(ii)] $\sum_i\nu_i|(Q\Gamma v)_i|^s\le\bar\gamma^s\sum_i\nu_i|v_i|^s$ for $s\ge1$,
by Jensen row by row and stationarity, $\sum_i\nu_iQ_{ij}=\nu_j$;
\item[(iii)] consequently $(I-Q\Gamma)^{-1}=\sum_t(Q\Gamma)^t$ has norm at most
$(1-\bar\gamma)^{-1}$ on $\ell^\infty$ and on every $\ell^s(\nu)$, and the same holds for
$\Gamma Q$.
\end{itemize}
\end{enumerate}

\section{The exact identity and the linear statistic}\label{app:identity}

\begin{proof}[Proof of Lemma~\ref{lem:identity}]
For \eqref{eq:Texact}, note that $T_F$ is itself a Markov operator: by the definition
of $g_F$,
\[
\frac{1}{g_F(x)}\int w(x,v)\,dF(v)=1\qquad\text{for every }x .
\]
Hence $(Tf)(x)=g_F(x)^{-1}\int w(x,v)\,(Tf)(x)\,dF(v)$, and subtracting this from the
definition of $T_Ff$ gives
\[
(T_Ff)(x)-(Tf)(x)=\frac{1}{g_F(x)}\int w(x,v)\big(f(v)-(Tf)(x)\big)\,dF(v).
\]
Writing $w(x,v)=g(x)p(v\mid x)$, this is $(g/g_F)\int p(v\mid\cdot)\big(f(v)-(Tf)(\cdot)\big)\,dF(v)$. The same integral against $dF^\circ(v)=dv$ vanishes, because $\int p(v\mid x)\,dv=1$. So the integral against $dF$ equals the integral against $d(F-F^\circ)$, which is $\dot T_{F-F^\circ}f$. This proves \eqref{eq:Texact}.

For \eqref{eq:yexact}, subtract $y^\circ=z+\gamma Ty^\circ$ from $y^\circ_F=z+\gamma T_Fy^\circ_F$. This gives
\[
y^\circ_F-y^\circ=\gamma T_F(y^\circ_F-y^\circ)+\gamma(T_F-T)y^\circ,\qquad\text{so}\qquad
y^\circ_F-y^\circ=(I-\gamma T_F)^{-1}\gamma(T_F-T)y^\circ .
\]
By \eqref{eq:Texact}, $(T_F-T)y^\circ=(g/g_F)\dot T_{F-F^\circ}y^\circ$. Next, the resolvent identity
$(I-\gamma T_F)^{-1}=\cR+(I-\gamma T_F)^{-1}\gamma(T_F-T)\cR$, which follows by multiplying
$(I-\gamma T_F)-(I-T_\gamma)=-\gamma(T_F-T)$ on the left by $(I-\gamma T_F)^{-1}$ and on the
right by $\cR$, gives
\[
y^\circ_F-y^\circ=\cR\gamma\frac g{g_F}\dot T_{F-F^\circ}y^\circ+(I-\gamma T_F)^{-1}\gamma(T_F-T)\cR\gamma\frac g{g_F}\dot T_{F-F^\circ}y^\circ .
\]
Finally, writing $g/g_F=1+(g/g_F-1)$ in the first term isolates $\cR\gamma\dot T_{F-F^\circ}y^\circ$, and the two remaining terms are $\mathcal Q_F$.
\end{proof}

Lemma~\ref{lem:identity} identifies the first-order term, but its remainder
$\mathcal Q_{F_N}$ is awkward to bound directly. It involves the random operator
$(I-\gamma T_{F_N})^{-1}$, and evaluated at the sample points its linear term contains
the diagonal value $\phi(\xi_i,\xi_i)$ of an a.e.-defined kernel. We therefore work
with an equivalent matrix form. It expresses the remainder as the expected-weight
resolvent $M=(I-\Gamma\hat P)^{-1}$ applied to an explicit per-node error $\eta$, which
involves only pairs of distinct nodes. Bounding the remainder then reduces to two
steps: a contraction bound for $M$ (Lemma~\ref{lem:A-transfer}) and a second-moment
bound for $\eta$ (Appendix~\ref{app:remainder}).

\begin{lemma}[Matrix form of the identity]\label{lem:A-matrix}
Assume (M1)--(M2) and $y^\circ\in L^1(\pi)$. Define, for $i\in[N]$,
\begin{gather*}
\chi_i:=\frac1{N-1}\sum_{j\ne i}p(\xi_j \mid \xi_i)\,e_{\xi_i}(\xi_j),\qquad
\cU_i:=\frac1{(N-1)^2}\sum_{\substack{j\ne l\\ j,l\ne i}}k_i(\xi_j,\xi_l),
\\
\cD_i:=\frac1{(N-1)^2}\sum_{j\ne i}p(\xi_j \mid \xi_i)\big[(T\phi_{\xi_j})(\xi_i)-\phi(\xi_j,\xi_i)\big],
\end{gather*}
where $k_i(u,v):=p(u\mid\xi_i)\big[(T\phi_v)(\xi_i)-\phi(u,v)\big]$. Here
$\phi_v=\phi(\cdot,v)$ is the influence of a type-$v$ node as a function of the
receiving type, so that $(T\phi_{\xi_j})(\xi_i)$ is the one-step average, around node
$i$, of node $j$'s influence.

Let $R:=Y^W-Y^\circ-\Phi$. Then, almost surely,
\begin{equation}\label{eq:A-matrix}
R=M\Gamma\eta,\qquad \eta_i:=\Big(\frac{g_i}{\hat g_i}-1\Big)\chi_i-\frac{g_i}{\hat g_i}\big(\cU_i+\cD_i\big).
\end{equation}
\end{lemma}

This is \eqref{eq:matrix-identity} of the main text.

\begin{proof}
Recall $\hat P_{ij}=(g_i/\hat g_i)p(\xi_j \mid \xi_i)/(N-1)$ for $j\ne i$, and $\sum_j\hat P_{ij}=1$. Equivalently,
\begin{equation}\label{eq:A-ghatratio}
\frac{\hat g_i}{g_i}=\frac1{N-1}\sum_{j\ne i}p(\xi_j\mid\xi_i).
\end{equation}
The proof has four steps: an equation for $Y^W-y^\circ$, an equation for $\Phi$, their difference, and the identification of the difference with $\eta$.

\emph{Step 1: the equation for $Y^W-y^\circ$.} Write $y^\circ_i:=y^\circ(\xi_i)$. Evaluating $y^\circ=z+\gamma Ty^\circ$ at
$\xi_i$ gives $y^\circ_i=Z_i+\gamma_i(Ty^\circ)(\xi_i)$, while
$Y^W_i=Z_i+\gamma_i(\hat PY^W)_i$. Subtracting, and adding and subtracting
$\gamma_i(\hat Py^\circ)_i$,
\[
Y^W_i-y^\circ_i=\gamma_i\big(\hat P(Y^W-y^\circ)\big)_i+\gamma_i\upsilon_i,
\qquad \upsilon_i:=(\hat Py^\circ)_i-(Ty^\circ)(\xi_i),
\]
that is, $(I-\Gamma\hat P)(Y^W-y^\circ)=\Gamma\upsilon$. Since the rows of $\hat P$
sum to one, $\upsilon_i=\sum_{j\ne i}\hat P_{ij}\big(y^\circ_j-(Ty^\circ)(\xi_i)\big)
=\sum_{j\ne i}\hat P_{ij}\,e_{\xi_i}(\xi_j)$, and substituting
$\hat P_{ij}=(g_i/\hat g_i)\,p(\xi_j\mid\xi_i)/(N-1)$ gives
$\upsilon_i=(g_i/\hat g_i)\chi_i$.

\emph{Step 2: the equation for $\Phi$.} By the fixed-point equation $\phi_v=f_v+\gamma T\phi_v$ of Appendix~\ref{app:notation}, which holds almost surely at the sample points, evaluating at $\xi_i$ with $v=\xi_l$ and averaging over $l\ne i$ gives
\begin{gather*}
\Phi_i=\frac1{N-1}\sum_{l\ne i}f_{\xi_l}(\xi_i)+\gamma_i\frac1{N-1}\sum_{l\ne i}(T\phi_{\xi_l})(\xi_i)=\gamma_i(\chi_i+\mathfrak a_i),\\
\mathfrak a_i:=\frac1{N-1}\sum_{l\ne i}(T\phi_{\xi_l})(\xi_i),
\end{gather*}
since $f_{\xi_l}(\xi_i)=\gamma_ip(\xi_l\mid\xi_i)e_{\xi_i}(\xi_l)$ averages to $\gamma_i\chi_i$. Subtracting $\gamma_i(\hat P\Phi)_i$ from both sides,
\[
\big((I-\Gamma\hat P)\Phi\big)_i=\gamma_i(\chi_i+\varrho_i),\qquad \varrho_i:=\mathfrak a_i-(\hat P\Phi)_i .
\]

\emph{Step 3: the difference.} Since $R=(Y^W-y^\circ)-\Phi$, Steps 1 and 2 give
\[
(I-\Gamma\hat P)R=\Gamma\Big(\frac g{\hat g}\chi-\chi-\varrho\Big).
\]

\emph{Step 4: identification of $\varrho$.} By \eqref{eq:A-ghatratio}, $\mathfrak a_i=(g_i/\hat g_i)(N-1)^{-1}\sum_{j\ne i}p(\xi_j\mid\xi_i)\,\mathfrak a_i$, while $(\hat P\Phi)_i=(g_i/\hat g_i)(N-1)^{-1}\sum_{j\ne i}p(\xi_j\mid\xi_i)\,\Phi_j$. Hence
\[
\varrho_i=\frac{g_i}{\hat g_i}\cdot\frac1{N-1}\sum_{j\ne i}p(\xi_j \mid \xi_i)\,(\mathfrak a_i-\Phi_j).
\]
In the difference $\mathfrak a_i-\Phi_j=(N-1)^{-1}\big[\sum_{l\ne i}(T\phi_{\xi_l})(\xi_i)-\sum_{l\ne j}\phi(\xi_j,\xi_l)\big]$, the indices $l\notin\{i,j\}$ appear in both sums; the first sum has the extra index $l=j$ and the second the extra index $l=i$. So
\[
\mathfrak a_i-\Phi_j=\frac1{N-1}\Big[\sum_{l\notin\{i,j\}}\big((T\phi_{\xi_l})(\xi_i)-\phi(\xi_j,\xi_l)\big)+(T\phi_{\xi_j})(\xi_i)-\phi(\xi_j,\xi_i)\Big].
\]
Multiplying by $p(\xi_j\mid\xi_i)$ and averaging over $j$, the sum over $l\notin\{i,j\}$ gives $\cU_i$ (by the definition of $k_i$), and the last two terms give $\cD_i$. So $\varrho_i=(g_i/\hat g_i)(\cU_i+\cD_i)$. Substituting into Step 3,
\[
(I-\Gamma\hat P)R=\Gamma\Big(\Big(\frac g{\hat g}-1\Big)\chi-\frac g{\hat g}(\cU+\cD)\Big)=\Gamma\eta,
\]
and applying $M=(I-\Gamma\hat P)^{-1}$ gives $R=M\Gamma\eta$.
\end{proof}

\begin{lemma}[Contraction and transfer to subsets]\label{lem:A-transfer}
Let $\hat\pi_i=\hat g_i/\sum_l\hat g_l$. Then $\hat P$ is $\hat\pi$-reversible and $\|M\Gamma\|_{\ell^2(\hat\pi)}\le\bar\gamma/(1-\bar\gamma)$. Let $r_N^2:=N\,\E[\hat g_1\eta_1^2]$ and $C_R:=\bar\gamma^2/((1-\bar\gamma)^2c_0)$. Then for every $S\subseteq[N]$ with $|S|=n$,
\[
\avg_{i\in S}R_i^2\le\frac{C_R}{n}\sum_{i=1}^N\hat g_i\eta_i^2,\qquad
\E\,\avg_{i\in S}R_i^2\le\frac{C_Rr_N^2}{n},\qquad
\sqrt n\,\E\Big|\avg_{i\in S}R_i\Big|\le C_R^{1/2}r_N .
\]
In particular $\E R_1^2\le C_Rr_N^2/N$, and $|\E\bar Y^W_S-\mu^\circ_N|\le C_R^{1/2}r_N/\sqrt n$.
\end{lemma}

\begin{proof}
\emph{Contraction.} $\hat g_i\hat P_{ij}=w(\xi_i,\xi_j)/(N-1)$ is symmetric, so $\hat P$ is reversible with respect to $\hat g$, and hence with respect to $\hat\pi$. By \ref{F:ops}(ii)--(iii) with $\nu=\hat g$, $\|\Gamma\hat P\|_{\ell^2(\hat g)}\le\bar\gamma$, so $\|M\Gamma\|_{\ell^2(\hat g)}\le\sum_{t\ge0}\bar\gamma^t\cdot\bar\gamma=\bar\gamma/(1-\bar\gamma)$; the same bound holds in $\ell^2(\hat\pi)$, which differs by a constant rescaling of the weights.

\emph{First bound.} Using $\hat g_i\ge c_0$ to pass from the average over $S$ to a $\hat g$-weighted sum over all nodes, and then the contraction with $R=M\Gamma\eta$,
\[
\avg_{i \in S}R_i^2\le\frac{1}{nc_0}\sum_{i=1}^N\hat g_iR_i^2
\le\frac{\bar\gamma^2}{(1-\bar\gamma)^2nc_0}\sum_{i=1}^N\hat g_i\eta_i^2
=\frac{C_R}{n}\sum_{i=1}^N\hat g_i\eta_i^2 .
\]

\emph{Second and third bounds.} Taking expectations and using exchangeability, $\E\sum_i\hat g_i\eta_i^2=N\E[\hat g_1\eta_1^2]=r_N^2$, which gives the second bound. The third follows from the second by Jensen, $\E|\avg_SR_i|\le(\E\avg_SR_i^2)^{1/2}$.

\emph{Final steps.} Taking $S=[N]$ and using exchangeability, $\E R_1^2=\E\avg_{[N]}R_i^2\le C_Rr_N^2/N$. For the mean, $\bar Y^W_S-\mu^\circ_N=\avg_S(y^\circ_i-\mu^\circ_N)+\avg_S\Phi_i+\avg_SR_i$; the first two terms have mean zero (the second by Lemma~\ref{lem:A-pop}(d)), so $|\E\bar Y^W_S-\mu^\circ_N|=|\E\avg_SR_i|\le C_R^{1/2}r_N/\sqrt n$.
\end{proof}

\begin{proof}[Proof of Proposition~\ref{prop:linear}]
\emph{The variance of $L_S$.} Write $a=H_{\rm root}-\mu^\circ_N$ and $b=H_{\rm rest}$, both centered. A node $i\in S$ contributes $a(\xi_i)/n+b(\xi_i)/N$ to $L_S$, and a node $l\notin S$ contributes $b(\xi_l)/N$. These contributions are independent across nodes, so
\[
\Var L_S=n\Var\Big(\frac an+\frac bN\Big)+(N-n)\frac{\Var b}{N^2}
=\frac{\Var a}n+\frac{2\Cov(a,b)}N+\frac{\Var b}N ,
\]
and therefore
\[
n\Var L_S=\Var a+\frac{2n}N\Cov(a,b)+\frac nN\Var b=\Big(1-\frac nN\Big)\Var a+\frac nN\Var(a+b).
\]

\emph{The centered kernel.} Let $\psi(x,v):=\phi(x,v)-H_{\rm rest}(v)$. Then $\int\psi(x,v)\,dx=0$ by the definition of $H_{\rm rest}$, and $\int\psi(x,v)\,dv=0$ for a.e.\ $x$ by Lemma~\ref{lem:A-pop}(d), so $\psi$ is degenerate. Since $H_{\rm rest}(v)=\int\phi(x,v)\,dx$ is the $L^2(dx\,dv)$-projection of $\phi$ onto functions of $v$, we have $\|\psi\|_{L^2}^2=\|\phi\|_{L^2}^2-\|H_{\rm rest}\|_{L^2}^2\le\|\phi\|^2_{L^2}$.

\emph{The decomposition.} Write $\phi(\xi_i,\xi_l)=H_{\rm rest}(\xi_l)+\psi(\xi_i,\xi_l)$ in
\[
\avg_{i\in S}\Phi_i=\frac1{n(N-1)}\sum_{i\in S}\sum_{l\ne i}\phi(\xi_i,\xi_l).
\]
For each $l$, the number of $i\in S$ with $i\ne l$ is $n-\one\{l\in S\}$, so
\[
\avg_{i\in S}\Phi_i=\sum_{l=1}^Nc_lH_{\rm rest}(\xi_l)+D_S,\qquad D_S:=\frac1{n(N-1)}\sum_{i\in S}\sum_{l\ne i}\psi(\xi_i,\xi_l),
\]
with $c_l=1/(N-1)$ for $l\notin S$ and $c_l=(n-1)/(n(N-1))$ for $l\in S$. The statistic $L_S$ uses the weight $1/N$ for every $l$; the difference is $E_S:=\sum_l(c_l-N^{-1})H_{\rm rest}(\xi_l)$.

\emph{The term $E_S$.} One checks $|c_l-N^{-1}|=1/(N(N-1))$ for $l\notin S$ and $|c_l-N^{-1}|\le1/(n(N-1))$ for $l\in S$. By independence,
\begin{align*}
\E E_S^2=\sum_l(c_l-N^{-1})^2\Var H_{\rm rest}&\le\Big(\frac{N}{N^2(N-1)^2}+\frac{n}{n^2(N-1)^2}\Big)\Var H_{\rm rest}\\
&\le C\big(N^{-3}+n^{-1}N^{-2}\big)\Var H_{\rm rest}.
\end{align*}

\emph{The term $D_S$.} Expand $\E D_S^2$ as a sum over quadruples $(i,l,i',l')$ with $i\ne l$, $i'\ne l'$, of $\E[\psi(\xi_i,\xi_l)\psi(\xi_{i'},\xi_{l'})]$. If some index occurs exactly once among $(i,l,i',l')$, condition on the other variables: by independence, and since $\psi$ has mean zero in each argument for every value of the other, the conditional expectation vanishes. So the surviving terms have $\{i,l\}=\{i',l'\}$: either $(i',l')=(i,l)$, at most $n(N-1)$ terms, or $(i',l')=(l,i)$, which requires $l\in S$, at most $n(n-1)$ terms. Each is at most $\|\psi\|^2_{L^2}$ by Cauchy--Schwarz. Hence
\[
\E D_S^2\le\frac{n(N-1)+n(n-1)}{n^2(N-1)^2}\|\psi\|^2_{L^2}\le\frac{2\|\phi\|^2_{L^2}}{n(N-1)} .
\]
\end{proof}

\begin{proof}[Proof of Theorem~\ref{thm:clt}]
By \eqref{eq:split}, $\sqrt n(\bar Y_S-\mu_N)=\sqrt n(\bar Y^W_S-\E\bar Y^W_S)+\sqrt n(\bar V_S-\E\bar V_S)$. The second term tends to $0$ in $L^2$ by Theorem~\ref{thm:fluct}, proved in Appendix~\ref{app:fluctuation}. The first equals $\sqrt nL_S+o_P(1)$ by Propositions~\ref{prop:linear} and~\ref{prop:remainder}. $L_S$ is a sum of independent centered terms with $n\Var L_S=\sigma^2_{n,N}$, which is bounded below by (M5). The Lyapunov condition holds by Lemma~\ref{lem:A-lyap}. The Lindeberg--Feller theorem for triangular arrays finishes the proof.
\end{proof}

\section{Population bounds and the signal remainder}\label{app:remainder}
\begin{table}[htbp]
\centering\small
\begin{tabular}{p{0.22\linewidth}p{0.55\linewidth}p{0.15\linewidth}}
\toprule
Symbol & Meaning & Defined in\\
\midrule
$K_t(x)$ & $\int p(v\mid x)^tdv$; $K_2$ measures how concentrated the neighborhood of type $x$ is & App.~\ref{app:notation}\\
$\bar T$ & walk averaged over a geometric number of steps, $(1-\bar\gamma)\sum_k\bar\gamma^kT^k$ & App.~\ref{app:notation}\\
$\Phi_s(u)$ & size $\|\phi(u,\cdot)\|_{L^s}$ of the influence received by type $u$ & App.~\ref{app:notation}\\
$\chi_i$, $\cU_i$, $\cD_i$ & centered neighborhood average; degenerate $U$-statistic; diagonal term & Lem.~\ref{lem:A-matrix}\\
$\Ohgone$, $\omega_1$ & event $\hat g_1\ge g(\xi_1)/2$; relative error $\hat g_1/g(\xi_1)-1$ & Lem.~\ref{lem:A-cond}\\
$\tilde F(x)^2$ & second moment of one neighbor's contribution to $x$'s centered average & Lem.~\ref{lem:A-cond}\\
$\mathfrak c_1$, $\mathfrak c_2$ & return of a node's influence to itself (hub term); its second-moment analogue & Lem.~\ref{lem:A-cond}\\
$\Lambda_N$, $\Lambda'_N$ & pair and coincidence functionals controlling $r_N^2$ & Lem.~\ref{lem:A-cond}\\
$\psi_u$ & discounted multi-step density $\sum_k\bar\gamma^kp^{(k+1)}(\cdot\mid u)$ (setting (B)) & Lem.~\ref{lem:A-rowcol}\\
$\mathfrak y$ & $\bar T|Ty^\circ|$ (setting (U)) & Prop.~\ref{prop:A-remU}\\
\bottomrule
\end{tabular}
\caption{Notation for the signal remainder.}
\label{tab:not-remainder}
\end{table}

\subsection{Population bounds}
We now derive bounds for various quantities defined in the population objects paragraph of Appendix \ref{app:notation}.

\begin{lemma}[Population bounds]\label{lem:A-pop}
Assume (M1)--(M2) and $y^\circ\in L^1(\pi)$. For $s\ge1$:
\begin{itemize}
\item[(a)] $\|y^\circ\|_{L^s(\pi)}\le(1-\bar\gamma)^{-1}\|Z\|_{L^s(\pi)}$ and $\|y^\circ\|_{L^s(dx)}\le c_0^{-1/s}\|y^\circ\|_{L^s(\pi)}$;
\item[(b)] $\int\Phi_s^s\,d\pi\le(1-\bar\gamma)^{-s}\,\Theta_s$;
\item[(c)] $\|H_{\rm rest}\|^s_{L^s}\le c_0^{-1}(1-\bar\gamma)^{-s}\Theta_s$ and $\|\phi\|^2_{L^2(dx\,dv)}\le c_0^{-1}(1-\bar\gamma)^{-2}\Theta_2$;
\item[(d)] $\Theta_1\le2\|y^\circ\|_{L^1(\pi)}<\infty$, and $\int\phi(x,v)\,dv=0$ for a.e.\ $x$. In particular $\E\Phi_i=0$ and $\E H_{\rm rest}(\xi_1)=0$.
\end{itemize}
\end{lemma}

\begin{proof}
(a) is \ref{F:contr} together with \ref{F:M2}.

(b) Let $F_s(x):=\|f_\cdot(x)\|_{L^s(dv)}=\big(\int|f_v(x)|^s\,dv\big)^{1/s}$ be the size of the one-step influence received by type $x$. Since $|\gamma|\le1$ and $w\,p^{s-1}=g\,p^s$,
\[
\int F_s^s\,d\pi\le\iint g(x)\,p(v\mid x)^s|e_x(v)|^s\,dv\,dx=\Theta_s .
\]
We compare $\Phi_s$ with $F_s$ in three steps.
\begin{itemize}
\item \emph{Expand the resolvent.} $\phi(u,v)=\sum_{k\ge0}(T_\gamma^kf_v)(u)$ and $|T_\gamma^kf_v|\le\bar\gamma^kT^k|f_v|$ pointwise, so by Minkowski's inequality in $L^s(dv)$,
\[
\Phi_s(u)\le\sum_{k\ge0}\bar\gamma^k\big\|(T^k|f_\cdot|)(u)\big\|_{L^s(dv)}\le\sum_{k\ge0}\bar\gamma^k(T^kF_s)(u),
\]
where the second step is Minkowski's integral inequality, since $T^k$ averages over the type with a probability kernel.
\item \emph{Jensen over walk lengths.} The weights $(1-\bar\gamma)\bar\gamma^k$ sum to one, so by Jensen, and then by Jensen for the probability kernel $T^k$,
\[
\Phi_s^s\le(1-\bar\gamma)^{-s}\Big(\sum_k(1-\bar\gamma)\bar\gamma^kT^kF_s\Big)^s\le(1-\bar\gamma)^{1-s}\sum_k\bar\gamma^k\,T^k(F_s^s).
\]
\item \emph{Stationarity.} $\int T^k(F_s^s)\,d\pi=\int F_s^s\,d\pi$, so integrating against $\pi$ gives $\int\Phi_s^s\,d\pi\le(1-\bar\gamma)^{1-s}(1-\bar\gamma)^{-1}\Theta_s=(1-\bar\gamma)^{-s}\Theta_s$.
\end{itemize}

(c) Since $H_{\rm rest}(v)=\int\phi(x,v)\,dx$, Minkowski's integral inequality gives $\|H_{\rm rest}\|_{L^s(dv)}\le\int\|\phi(x,\cdot)\|_{L^s}\,dx=\int\Phi_s\,dx$. By Jensen on the probability space $([0,1],dx)$, $\int\Phi_s\,dx\le(\int\Phi_s^s\,dx)^{1/s}$, and $dx\le d\pi/c_0$ by \ref{F:M2}. Combining with (b) gives the first claim. For the second, $\|\phi\|^2_{L^2(dx\,dv)}=\int\Phi_2^2\,dx\le c_0^{-1}\int\Phi_2^2\,d\pi$, and (b) with $s=2$ applies.

(d) By the triangle inequality, $\Theta_1=\iint g\,p\,|e_x(v)|\le\int(T|y^\circ|)\,d\pi+\int|Ty^\circ|\,d\pi\le2\|y^\circ\|_{L^1(\pi)}$, using stationarity. In particular $Ty^\circ$ is finite a.e., and
\[
\int f_v(x)\,dv=\gamma(x)\int p(v\mid x)\big(y^\circ(v)-(Ty^\circ)(x)\big)\,dv=\gamma(x)\big((Ty^\circ)(x)-(Ty^\circ)(x)\big)=0 .
\]
To pass this through the resolvent, note that $\iint\sum_k\bar\gamma^kT^k|f_v|\,dv\,d\pi\le(1-\bar\gamma)^{-1}\Theta_1<\infty$ by stationarity, so Fubini allows exchanging $\int dv$ with the Neumann series. This gives $\int\phi(x,v)\,dv=\sum_kT_\gamma^k\big(\int f_v\,dv\big)(x)=0$ for $\pi$-a.e., hence a.e., $x$. The last two statements follow: $\E\Phi_i=\E\big[\E[\phi(\xi_i,\xi_l)\mid\xi_i]\big]=0$, and $\E H_{\rm rest}(\xi_1)=\iint\phi(x,v)\,dx\,dv=0$.
\end{proof}

We now provide a lemma that bounds $\Theta_s=\int\big(\int|p(v\mid x)e_x(v)|^s\,dv\big)\,\pi(dx)$, which is the
$s$-th moment of the one-step influence kernel: how much a single node of type $v$
moves the neighborhood average of a node of type $x$, averaged over receiving types
with the stationary weight $\pi$. Since the total influence $\phi$ is this one-step
influence propagated by the walk, and the walk contracts $L^s(\pi)$, the total
influence on the population mean, $H_{\rm rest}$, is controlled by $\Theta_s$ with no
dependence on walk length (Lemma~\ref{lem:A-pop}). The power $p^s$ is where hubs enter:
a type $v$ with large $p(v\mid x)$ has a large share of $x$'s neighborhood, and
$\Theta_s<\infty$ limits how much of the influence such types can carry.
\begin{lemma}[$\Theta_s$ in the two settings]\label{lem:A-theta}
\begin{itemize}
\item[(B)] If $|Z|\le B$, then $\Theta_{2+\delta}\le(2B_y)^{2+\delta}c_0^{-1-\delta}\iint w^{2+\delta}$ and
\[
\Theta_2\le4B_y^2c_0^{-1}\iint w^2\le4B_y^2c_0^{-1}\Big(\iint w^{2+\delta}\Big)^{2/(2+\delta)}.
\]
\item[(U)] Let $\mathsf N_\delta:=\int\big(\int w(x,v)^{2+\delta}dx\big)^{4/(2-\delta)}dv$. Then
\begin{gather*}
\Theta_{2+\delta}\le C\,\mathsf N_\delta^{(2-\delta)/4}\,\E[gZ^4]^{(2+\delta)/4},\\
\Theta_2\le C\,\mathsf M_N^{1/2}\,\E[gZ^4]^{1/2}.
\end{gather*}
\end{itemize}
In particular, $\Theta_2$ and $\Theta_{2+\delta}$ are bounded under (B1)--(B2), respectively (U2)--(U4).
\end{lemma}

\begin{proof}
Since $g\ge c_0$, $w\,p^{s-1}=w^s/g^{s-1}\le c_0^{1-s}w^s$, so in both settings
\[
\Theta_s\le c_0^{1-s}\iint w(x,v)^s\,|e_x(v)|^s\,dv\,dx .
\]

(B) Here $|e_x(v)|\le2B_y$ by \ref{F:contr}, which gives the bounds on $\Theta_{2+\delta}$ and on $\Theta_2$ in terms of $\iint w^2$. For the last inequality, apply Jensen on the probability space $([0,1]^2,dx\,dv)$: $\iint w^2\le(\iint w^{2+\delta})^{2/(2+\delta)}$.

(U) First, $|e_x(v)|^s\le2^{s-1}(|y^\circ(v)|^s+|Ty^\circ(x)|^s)$. For $s=2+\delta$, consider the term with $y^\circ(v)$; integrating over $x$ first and using symmetry of $w$, then H\"older in $v$ with exponents $\big(\tfrac{4}{2-\delta},\tfrac4{2+\delta}\big)$,
\[
\iint w(x,v)^{2+\delta}|y^\circ(v)|^{2+\delta}\,dx\,dv=\int\Big(\int w(x,v)^{2+\delta}dx\Big)|y^\circ(v)|^{2+\delta}dv\le\mathsf N_\delta^{(2-\delta)/4}\,\|y^\circ\|_{L^4}^{2+\delta},
\]
and the term with $Ty^\circ(x)$ is handled in the same way with the roles of $x$ and $v$ exchanged. By \ref{F:M2} and \ref{F:contr}, $\|Ty^\circ\|_{L^4}$ and $\|y^\circ\|_{L^4}$ are at most $C\|Z\|_{L^4(\pi)}=C\,\E[gZ^4]^{1/4}$. For $\Theta_2$, similarly, by Cauchy--Schwarz,
\[
\iint w^2\big(y^\circ(v)^2+(Ty^\circ)(x)^2\big)=\int\kappa_2\big((y^\circ)^2+(Ty^\circ)^2\big)\le\mathsf M_N^{1/2}\big(\|y^\circ\|^2_{L^4}+\|Ty^\circ\|^2_{L^4}\big).
\]
\end{proof}

We now verify Lyapunov's condition for the main projection term.
\begin{lemma}[Lyapunov condition]\label{lem:A-lyap}
Under (B1)--(B2), respectively (U2)--(U3), $\|H_{\rm root}-\mu^\circ_N\|_{L^{2+\delta}}+\|H_{\rm rest}\|_{L^{2+\delta}}\le C$. Hence, with $L_S=\sum_jX_j$ as in Proposition~\ref{prop:linear},
\[
\frac{\sum_j\E|X_j|^{2+\delta}}{(\Var L_S)^{1+\delta/2}}\le\frac{C\,n^{-\delta/2}}{\sigma_{n,N}^{2+\delta}}\to0\qquad\text{under (M5)}.
\]
\end{lemma}
\begin{proof}
\emph{Step 1: moment bounds.} In setting (B), $|y^\circ|\le B_y$ by \ref{F:contr}, so
$\|y^\circ-\mu^\circ_N\|_{L^{2+\delta}}\le2B_y$. In setting (U), since $2+\delta\le4$
and $dx$ is a probability measure on $[0,1]$, Lemma~\ref{lem:A-pop}(a) and (U2) give
\[
\|y^\circ-\mu^\circ_N\|_{L^{2+\delta}}\le2\|y^\circ\|_{L^{2+\delta}}\le2\|y^\circ\|_{L^4}
\le2c_0^{-1/4}(1-\bar\gamma)^{-1}\E[gZ^4]^{1/4}\le C .
\]
In both settings, Lemma~\ref{lem:A-pop}(c) gives
$\|H_{\rm rest}\|^{2+\delta}_{L^{2+\delta}}\le c_0^{-1}(1-\bar\gamma)^{-(2+\delta)}\Theta_{2+\delta}$,
and $\Theta_{2+\delta}\le C$ by Lemma~\ref{lem:A-theta}, under (B1)--(B2), respectively
(U2)--(U3).

\emph{Step 2: the Lyapunov ratio.} The summands of $L_S$ are
$X_j=\one\{j\in S\}\big(y^\circ(\xi_j)-\mu^\circ_N\big)/n+H_{\rm rest}(\xi_j)/N$. For
$j\in S$, the inequality $|a+b|^{2+\delta}\le2^{1+\delta}(|a|^{2+\delta}+|b|^{2+\delta})$
and Step~1 give $\E|X_j|^{2+\delta}\le C\big(n^{-(2+\delta)}+N^{-(2+\delta)}\big)$; for
$j\notin S$, $\E|X_j|^{2+\delta}\le CN^{-(2+\delta)}$. Summing over the $n$ sampled
and $N-n$ unsampled nodes,
\[
\sum_j\E|X_j|^{2+\delta}\le C\big(n^{-(1+\delta)}+N^{-(1+\delta)}\big)\le Cn^{-(1+\delta)} .
\]
Finally, $\Var L_S=\sigma^2_{n,N}/n$ by Proposition~\ref{prop:linear}, so
$(\Var L_S)^{1+\delta/2}=\sigma_{n,N}^{2+\delta}n^{-(1+\delta/2)}$, and the ratio is at
most $Cn^{-\delta/2}/\sigma_{n,N}^{2+\delta}$, which tends to zero under (M5).
\end{proof}

The following lemma provides insight into the components of $\sigma_{n,N}^2$.
\begin{lemma}[The limiting variance]\label{lem:variance-bounds}
For every $n\le N$:
\begin{itemize}
\item[(a)] $\sigma^2_{n,N}\ge(1-n/N)\Var H_{\rm root}(\xi_1)$;
\item[(b)] $\Var H_{\rm full}(\xi_1)\le2\Var H_{\rm root}(\xi_1)+2\|H_{\rm rest}\|^2_{L^2}\le2\Var H_{\rm root}(\xi_1)+C\Theta_2$.
\end{itemize}
Consequently, (M5) holds whenever $\liminf\Var H_{\rm root}(\xi_1)>0$ and $\limsup n/N<1$. If $n/N\to0$, (M5) is equivalent to $\liminf\Var H_{\rm root}(\xi_1)>0$.
\end{lemma}

\begin{proof}
(a) holds because $\sigma^2_{n,N}$ is a convex combination of two nonnegative terms. (b) is $(a+b)^2\le2a^2+2b^2$ together with Lemma~\ref{lem:A-pop}(c). For the last statement, $\sigma^2_{n,N}-\Var H_{\rm root}=(n/N)(\Var H_{\rm full}-\Var H_{\rm root})$, so $|\sigma^2_{n,N}-\Var H_{\rm root}|\le(n/N)(\Var H_{\rm full}+\Var H_{\rm root})\to0$ by (b).
\end{proof}

\subsection{Conditional second moments of the error}
By Lemma~\ref{lem:A-transfer}, the signal remainder is controlled by the single
quantity $r_N^2=N\,\E[\hat g_1\eta_1^2]$. This subsection computes that second moment
conditionally on the type of node $1$, reducing it to two population functionals:
$\Lambda_N$, which collects the contributions of pairs of distinct other nodes and
enters at order $N^{-2}$, and $\Lambda'_N$, which collects coincidences, where the same
node appears twice, and enters at order $N^{-3}$. The next two subsections bound
$\Lambda_N$ and $\Lambda'_N$ in settings (B) and (U).

Fix node $1$, write $x=\xi_1$, and let $\Ohgone:=\{\hat g_1\ge g(x)/2\}$ and $\omega_1:=\hat g_1/g(x)-1$. Define
\begin{gather*}
\tilde F(x)^2:=\int p(v\mid x)^2e_x(v)^2dv,\\
\mathfrak c_1(x):=\int p(v\mid x)|\phi(v,x)|\,dv,\\
\mathfrak c_2(x):=\int p(v\mid x)^2\phi(v,x)^2dv,
\end{gather*}
and the integrated quantities
\begin{align*}
\Lambda_N&:=\int g(x)\Big[K_2(x)\tilde F(x)^2+K_2(x)\,T(\Phi_2^2)(x)
+\int p(u\mid x)^2\Phi_2(u)^2\,du+\mathfrak c_1(x)^2\Big]dx,\\
\Lambda'_N&:=\int g(x)\Big[\int\big(p(v\mid x)^4+p(v\mid x)^2\big)e_x(v)^2\,dv
+\int p(v\mid x)^2(T\phi_v)(x)^2\,dv+\mathfrak c_2(x)\Big]dx .
\end{align*}

Each functional has a simple interpretation in terms of the walk in type space.
\begin{itemize}
\item $K_2(x)=\int p(v\mid x)^2dv$ measures how concentrated the neighborhood of a type
$x$ is: it is large when $x$'s neighbors come from a few types, and the average of
$N$ draws from $p(\cdot\mid x)$ then has variance of order $K_2(x)/N$. Its integrated
square, $\int gK_2^2=\int\kappa_2^2/g^3$, is controlled by the mixed norm $\mathsf M_N$.
\item $\tilde F(x)^2=\int p(v\mid x)^2e_x(v)^2dv$ is the second moment of the
contribution of one random neighbor to $x$'s centered neighborhood average.
\item $\Phi_2(u)^2=\int\phi(u,v)^2dv$ is the influence profile of $u$: the second moment
of the total influence that a node of type $u$ receives from one node of random type.
\item $\mathfrak c_1(x)=\int p(v\mid x)|\phi(v,x)|\,dv$ is the influence of $x$ on its
own neighbors, averaged over them: how much of $x$'s influence returns to $x$ through
the network. This is the hub term, since a hub's neighbors are numerous and
$\phi(\cdot,x)$ grows with the degree of $x$; $\mathfrak c_2$ is its second-moment
analogue.
\item The terms with $p^4$, $p^2(T\phi)^2$ and $\mathfrak c_2$ arise when the same
other node appears twice in an average, and therefore carry an extra factor $1/N$.
\end{itemize}

\begin{lemma}[Conditional moments of $\eta_1$]\label{lem:A-cond}
Assume (M1)--(M4) and $y^\circ\in L^1(\pi)$. Then:
\begin{itemize}
\item[(a)] $\Pr(\Ohgone^c\mid x)\le e^{-\dmin/4}$;
\item[(b)] on $\Ohgone$, $\hat g_1\eta_1^2\le4g(x)\omega_1^2\chi_1^2+8g(x)(\cU_1^2+\cD_1^2)$;
\item[(c)] $\E[\hat g_1\eta_1^2\one_{\Ohgone}]\le C\big(\Lambda_N/N^2+\Lambda'_N/N^3\big)$.
\end{itemize}
\end{lemma}

\begin{proof}
(a) Given $x$, $\hat g_1/g(x)=\avg_{j\ne1} p(\xi_j \mid x)$ by \eqref{eq:A-ghatratio}: an average of $N-1$ i.i.d.\ nonnegative terms with mean $1$ and second moment $K_2(x)\le(2\rho c_0)^{-1}$, by \ref{F:M3}. The event $\Ohgone^c$ is that this average falls below $1/2$, that is, that the sum falls below its mean by $t=(N-1)/2$. Lemma~\ref{lem:A-lt} gives probability at most
\[
\exp\Big\{-\frac{(N-1)^2/4}{2(N-1)K_2(x)}\Big\}=\exp\Big\{-\frac{N-1}{8K_2(x)}\Big\}\le\exp\Big\{-\frac{(N-1)\rho c_0}{4}\Big\}\le e^{-\dmin/4}.
\]

(b) On $\Ohgone$, $\hat g_1\ge g/2$. The first term of $\eta_1$ contributes $\hat g_1(g/\hat g_1-1)^2\chi_1^2=g^2\omega_1^2\chi_1^2/\hat g_1\le2g\omega_1^2\chi_1^2$, and the second contributes $\hat g_1(g/\hat g_1)^2(\cU_1+\cD_1)^2\le2g(\cU_1+\cD_1)^2$. The claim follows from $(a+b)^2\le2a^2+2b^2$, applied once to the two terms of $\eta_1$ and once to $\cU_1+\cD_1$.

(c) We bound the conditional expectations given $x$ of the three quantities in (b), then multiply by $g(x)$ and integrate.

\emph{The product $\omega_1^2\chi_1^2$.} Given $x$, let $a_j=p(\xi_j \mid x)-1$ and $b_j=p(\xi_j \mid x)e_x(\xi_j)$ for $j\ne1$. These are i.i.d.\ and centered, since $\int p(v\mid x)\,dv=1$ and $\int p(v\mid x)e_x(v)\,dv=0$, and $\omega_1=\avg_ja_j$, $\chi_1=\avg_jb_j$. Expanding,
\[
\E[\omega_1^2\chi_1^2\mid x]=\frac1{(N-1)^4}\sum_{j,k,l,m}\E[a_ja_kb_lb_m\mid x],
\]
and a term vanishes unless every index appears at least twice. The surviving patterns are $j=k\ne l=m$, $j=l\ne k=m$ and $j=m\ne k=l$, each with at most $(N-1)^2$ terms, and $j=k=l=m$, with $N-1$ terms. Using $(\E ab)^2\le\E a^2\,\E b^2$,
\[
\E[\omega_1^2\chi_1^2\mid x]\le\frac{3\E a^2\E b^2}{(N-1)^2}+\frac{\E a^2b^2}{(N-1)^3}\le\frac{3K_2\tilde F^2}{(N-1)^2}+\frac{2\int(p^4+p^2)e_x^2}{(N-1)^3},
\]
since $\E a^2\le\int p^2=K_2$, $\E b^2=\tilde F^2$, and $a^2\le2p^2+2$.

\emph{The degenerate $U$-statistic $\cU_1$.} The kernel $k_1$ is degenerate in both arguments given $x$: $\int k_1(u,v)\,du=\int p(u\mid x)\,du\,(T\phi_v)(x)-(T\phi_v)(x)=0$, and $\int k_1(u,v)\,dv=0$ by Lemma~\ref{lem:A-pop}(d). So the computation of $\E D_S^2$ in the proof of Proposition~\ref{prop:linear}, with $N-1$ variables and normalization $(N-1)^{-2}$, gives $\E[\cU_1^2\mid x]\le2(N-1)^{-2}\iint k_1^2$. By $(a+b)^2\le2a^2+2b^2$,
\begin{align*}
\iint k_1^2&\le2K_2(x)\int(T\phi_v)(x)^2\,dv+2\iint p(u\mid x)^2\phi(u,v)^2\,du\,dv\\
&\le2K_2(x)\,T(\Phi_2^2)(x)+2\int p(u\mid x)^2\Phi_2(u)^2du,
\end{align*}
using $(T\phi_v)^2\le T(\phi_v^2)$ (Jensen) and $\int T(\phi_v^2)\,dv=T(\Phi_2^2)$. Hence $\E[\cU_1^2\mid x]\le4(N-1)^{-2}\big[\int p(u\mid x)^2\Phi_2(u)^2du+K_2(x)T(\Phi_2^2)(x)\big]$.

\emph{The diagonal term $\cD_1$.} Here $\cD_1=(N-1)^{-1}\avg_jc_j$ with i.i.d.\ $c_j:=p(\xi_j\mid x)\big[(T\phi_{\xi_j})(x)-\phi(\xi_j,x)\big]$, so
\[
\E[\cD_1^2\mid x]=\frac{(\E c)^2+\Var(c)/(N-1)}{(N-1)^2}\le\frac{(\E c)^2+\E c^2/(N-1)}{(N-1)^2}.
\]
By Cauchy--Schwarz, $\big(\int p\,|T\phi_v(x)|\,dv\big)^2\le K_2(x)\int(T\phi_v)(x)^2dv\le K_2T(\Phi_2^2)$, so
\[
(\E c)^2\le2K_2T(\Phi_2^2)+2\mathfrak c_1^2,\qquad
\E c^2\le2\int p^2(T\phi_v)^2+2\mathfrak c_2 .
\]

\emph{Collecting.} Multiplying the three bounds by $g(x)$, as in (b), and integrating over $x$: the terms with $(N-1)^{-2}$ are the four integrands of $\Lambda_N$, and the terms with $(N-1)^{-3}$ are the three integrands of $\Lambda'_N$.
\end{proof}

\subsection{Setting (B)}

\begin{lemma}[Row and column bounds, bounded inputs]\label{lem:A-rowcol}
Let $|Z|\le B$ and $\psi_u:=\sum_{k\ge0}\bar\gamma^kp^{(k+1)}(\cdot\mid u)$. For a.e.\ $u,v,x$:
\begin{itemize}
\item[(a)] $\tilde F(x)^2\le4B_y^2K_2(x)$;
\item[(b)] $|\phi(u,v)|\le2\bar\gamma B_y\psi_u(v)$ and $\Phi_2(u)^2\le CB^2\,(\bar TK_2)(u)$;
\item[(c)] $|\phi(v,x)|\le(2\bar\gamma B_y/c_0)\,g(x)\psi_x(v)$, and $\mathfrak c_1(x)^2\le CB^2g(x)^2K_2(x)(\bar TK_2)(x)$;
\item[(d)] $\|\psi_u\|_\infty\le C/\rho$ and $|(T\phi_v)(x)|\le CB/\rho$.
\end{itemize}
\end{lemma}

\begin{proof}
(a) follows from $|e|\le2B_y$.

(b) \emph{Pointwise bound.} Since $|e|\le2B_y$, $|f_v(s)|\le2\bar\gamma B_yp(v\mid s)$. Applying $T_\gamma^k$ and using $T^kp(v\mid\cdot)=p^{(k+1)}(v\mid\cdot)$ (Chapman--Kolmogorov),
\[
|\phi(u,v)|\le\sum_{k\ge0}\bar\gamma^k\,T^k|f_v|(u)\le2\bar\gamma B_y\sum_{k\ge0}\bar\gamma^kp^{(k+1)}(v\mid u)=2\bar\gamma B_y\psi_u(v).
\]
\emph{Norm bound.} As in the proof of Lemma~\ref{lem:A-pop}(b), $\Phi_2\le\sum_k\bar\gamma^kT^kF_2$, and here $F_2\le2\bar\gamma B_yK_2^{1/2}$. So $\Phi_2\le CB\,\bar T(K_2^{1/2})$, and $(\bar TK_2^{1/2})^2\le\bar TK_2$ by Jensen, since $\bar T$ is a Markov operator.

(c) \emph{Pointwise bound.} Exchanging the roles of the arguments in (b), $|\phi(v,x)|\le2\bar\gamma B_y\sum_k\bar\gamma^kp^{(k+1)}(x\mid v)$. Multi-step reversibility, $g(v)p^{(k+1)}(x\mid v)=g(x)p^{(k+1)}(v\mid x)$, together with $g(v)\ge c_0$, turns this into $|\phi(v,x)|\le(2\bar\gamma B_y/c_0)g(x)\psi_x(v)$.

\emph{Bound on $\mathfrak c_1$.} Inserting the pointwise bound,
\[
\mathfrak c_1(x)\le\frac{CB}{c_0}g(x)\sum_k\bar\gamma^k\int p(v\mid x)p^{(k+1)}(v\mid x)\,dv\le CBg(x)K_2(x)^{1/2}\sum_k\bar\gamma^k\|p^{(k+1)}(\cdot\mid x)\|_{L^2},
\]
by Cauchy--Schwarz. Since $p^{(k+1)}(\cdot\mid x)=\int p^{(k)}(s\mid x)p(\cdot\mid s)\,ds$, Minkowski gives $\|p^{(k+1)}(\cdot\mid x)\|_{L^2}\le(T^kK_2^{1/2})(x)$. Hence $\mathfrak c_1(x)\le CBg(x)K_2(x)^{1/2}(\bar TK_2^{1/2})(x)$. Square and apply Jensen's inequality, $(\bar TK_2^{1/2})^2\le\bar TK_2$.

(d) By \ref{F:M3}, $p(v\mid s)\le(2\rho c_0)^{-1}$, so $p^{(k+1)}(v\mid u)=\int p^{(k)}(s\mid u)p(v\mid s)\,ds\le(2\rho c_0)^{-1}$, and summing over $k$ gives $\|\psi_u\|_\infty\le C/\rho$. Then $|(T\phi_v)(x)|\le\int p(u\mid x)|\phi(u,v)|\,du\le2\bar\gamma B_y\sup_u\|\psi_u\|_\infty$.
\end{proof}

The proof of the next proposition has the same structure in both settings. On the
event $\Ohgone$, Lemma~\ref{lem:A-cond} reduces $r_N^2$ to $\Lambda_N/N$ and
$\Lambda'_N/N^2$, whose terms are of four kinds: the \emph{product term}
$\E[gK_2\tilde F^2]$; the two \emph{degenerate-$U$ terms}
$\E[gK_2\,T(\Phi_2^2)]$ and $\E[g\int p(u\mid\cdot)^2\Phi_2(u)^2\,du]$, which come
from $\cU_1$ and involve the influence profile $\Phi_2$; the \emph{hub term}
$\E[g\mathfrak c_1^2]$; and the three \emph{coincidence terms} that make up
$\Lambda'_N$,
\[
\E\Big[g\int(p^4+p^2)\,e^2\Big],\qquad\E\Big[g\int p(v\mid\cdot)^2(T\phi_v)^2\,dv\Big],\qquad\E[g\mathfrak c_2],
\]
which arise when the same other node appears twice in an average: the first from the
product $\omega_1\chi_1$, the second and third from the diagonal term $\cD_1$. In
setting (B), the row and column bounds of Lemma~\ref{lem:A-rowcol} reduce every term
to the norms $\|K_2\|_{L^2(\pi)}$ and $\|gK_2\|_{L^2(\pi)}$, both controlled by
$\mathsf M_N$. The implicit bound $g\le1/(2\rho_N)$ of (M3) is used once for the hub
term, giving the factor $\rho_N^{-1}$, and for the coincidence terms, giving
$\rho_N^{-2}$. Off $\Ohgone$, every quantity is polynomial in $B/\rho_N$, while
$\Ohgone^c$ has probability at most $e^{-\dmin/4}$, so this part is negligible.

\begin{proposition}[Remainder, setting (B)]\label{prop:A-remB}
Under (M1)--(M4) and (B1)--(B2), with $C_*$ large enough,
\[
r_N^2=N\,\E[\hat g_1\eta_1^2]\le CB^2\Big(\frac{\mathsf M_N(1+\rho_N^{-1})}{N}+\frac{1+\mathsf M_N}{(N\rho_N)^2}\Big).
\]
\end{proposition}

\begin{proof}
\emph{Preliminaries.} Throughout, we use Lemma~\ref{lem:A-rowcol}, \ref{F:M3}, the fact that $T$ and $\bar T$ contract $L^2(\pi)$, and the norms
\[
\|K_2\|^2_{L^2(\pi)}=\int\kappa_2^2/g^3\le c_0^{-3}\mathsf M_N,\qquad\|gK_2\|^2_{L^2(\pi)}=\int\kappa_2^2/g\le c_0^{-1}\mathsf M_N,
\]
which follow from $K_2=\kappa_2/g^2$ and $g\ge c_0$.

\emph{The functional $\Lambda_N$.} We bound its four terms in turn.
\begin{itemize}
\item \emph{Product term.} By Lemma~\ref{lem:A-rowcol}(a), $\E[gK_2\tilde F^2]\le4B_y^2\int K_2^2\,d\pi=4B_y^2\|K_2\|^2_{L^2(\pi)}\le CB^2\mathsf M_N$.
\item \emph{First degenerate-$U$ term.} By Lemma~\ref{lem:A-rowcol}(b), $T(\Phi_2^2)\le CB^2\,T\bar TK_2$, so by Cauchy--Schwarz and the $L^2(\pi)$ contraction of $T\bar T$,
\[
\E[gK_2T(\Phi_2^2)]\le CB^2\int K_2\,T\bar TK_2\,d\pi\le CB^2\|K_2\|^2_{L^2(\pi)}\le CB^2\mathsf M_N .
\]
\item \emph{Second degenerate-$U$ term.} Exchanging the order of integration, $\E[g\int p(u\mid\cdot)^2\Phi_2(u)^2du]=\int\Phi_2(u)^2\big(\int w(x,u)^2/g(x)\,dx\big)du$, and $\int w(x,u)^2/g(x)\,dx\le c_0^{-1}\kappa_2(u)=c_0^{-1}(gK_2)(u)g(u)$. So, by Lemma~\ref{lem:A-rowcol}(b) and Cauchy--Schwarz,
\[
\E\Big[g\int p(u\mid\cdot)^2\Phi_2(u)^2du\Big]\le CB^2\int\bar TK_2\cdot gK_2\,d\pi\le CB^2\mathsf M_N ,
\]
using $\int\bar TK_2\cdot gK_2\,d\pi\le\|K_2\|_{L^2(\pi)}\|gK_2\|_{L^2(\pi)}$.
\item \emph{Hub term.} By Lemma~\ref{lem:A-rowcol}(c), $g\mathfrak c_1^2\le CB^2g^3K_2\,\bar TK_2$. Using $g\le1/(2\rho)$ once, $\E[g\mathfrak c_1^2]\le CB^2\rho^{-1}\int gK_2\cdot\bar TK_2\,d\pi$, and Cauchy--Schwarz as in the previous item gives at most $CB^2\rho^{-1}\mathsf M_N$.
\end{itemize}

\emph{The functional $\Lambda'_N$.} Again term by term.
\begin{itemize}
\item \emph{First coincidence term.} $gK_4=\int w^4/g^3\le c_0^{-3}(2\rho)^{-2}\kappa_2$ by \ref{F:M3}, and $gK_2=\kappa_2/g\le\kappa_2/c_0$. With $|e|\le2B_y$, $\E[g\int(p^4+p^2)e^2]\le CB^2\rho^{-2}\int\kappa_2\le CB^2\rho^{-2}\mathsf M_N^{1/2}$.
\item \emph{Second coincidence term.} By Lemma~\ref{lem:A-rowcol}(d), $|T\phi_v|\le CB/\rho$, so $\E[g\int p^2(T\phi_v)^2]\le CB^2\rho^{-2}\E[gK_2]\le CB^2\rho^{-2}\mathsf M_N^{1/2}$.
\item \emph{Third coincidence term.} Using $p(\cdot\mid x)\le(2\rho g(x))^{-1}$, the pointwise bound of Lemma~\ref{lem:A-rowcol}(c) and $\|\psi_x\|_\infty\le C/\rho$, $\mathfrak c_2(x)\le CB^2\rho^{-2}g(x)K_2(x)^{1/2}(\bar TK_2^{1/2})(x)$. Hence, by Cauchy--Schwarz in $L^2(\pi)$,
\[
\E[g\mathfrak c_2]\le CB^2\rho^{-2}\|gK_2^{1/2}\|_{L^2(\pi)}\|K_2^{1/2}\|_{L^2(\pi)}\le CB^2\rho^{-2}\mathsf M_N,
\]
using $\int g^3K_2=\int g\kappa_2\le\mathsf M_N^{3/4}$ and $\int gK_2\le C\mathsf M_N^{1/2}$.
\end{itemize}
So $\Lambda_N\le CB^2\mathsf M_N(1+\rho^{-1})$ and $\Lambda'_N\le CB^2\rho^{-2}(1+\mathsf M_N)$, and Lemma~\ref{lem:A-cond}(c) gives the stated bound on $\Ohgone$.

\emph{Off $\Ohgone$.} By \ref{F:M3} and Lemma~\ref{lem:A-rowcol}(d), every factor in $\eta_1$ is bounded by a power of $1/\rho$: $|\eta_1|\le CB\rho^{-3}$, and $\hat g_1\le(2\rho)^{-1}$. Since $\rho\ge1/N$ and by Lemma~\ref{lem:A-cond}(a), $\E[\hat g_1\eta_1^2\one_{\Ohgone^c}]\le CB^2N^7e^{-\dmin/4}\le CB^2N^{-5}$ for $C_*\ge48$. This is below the stated bound.
\end{proof}

\subsection{Setting (U)}

\begin{proposition}[Remainder, setting (U)]\label{prop:A-remU}
Under (M1)--(M4) and (U1)--(U4), with $C_*\ge60$,
\[
r_N^2\le C\Big(\frac{W_N^2}{N}+\frac{W_N^3}{N^2}\Big)\to0 .
\]
\end{proposition}

\begin{proof}
\emph{Preliminaries.} Since $\iint w=1$, $W_N\ge1$. Pointwise, $K_2\le W_N/g$, $p(\cdot\mid x)\le W_N/g(x)$ and $p\le W_N/c_0$. Let $\mathfrak y:=\bar T|Ty^\circ|$, so that $\|\mathfrak y\|_{L^r(\pi)}\le\|y^\circ\|_{L^r(\pi)}$. Since $|f_v(s)|\le(W_N/c_0)\big(|y^\circ(v)|+|Ty^\circ(s)|\big)$, summing the Neumann series as in Lemma~\ref{lem:A-rowcol}(b) gives
\begin{equation}\label{eq:A-phiU}
|\phi(u,v)|\le C W_N\big(|y^\circ(v)|+\mathfrak y(u)\big),\qquad |(T\phi_v)(x)|\le CW_N\big(|y^\circ(v)|+(T\mathfrak y)(x)\big).
\end{equation}
Write $\Sigma_y:=\int(y^\circ)^2d\pi\le C\E[gZ^2]\le C$. All of the following bounds are used in the proof; $W_N$ enters every one of them.

\emph{The functionals.}
\begin{itemize}
\item \emph{Product term.} $\tilde F^2\le2\int p^2(y^\circ)^2+2K_2(Ty^\circ)^2\le(4W_N/g)\,T((y^\circ)^2)$, using $p\le W_N/g$ and $(Ty^\circ)^2\le T((y^\circ)^2)$. With $K_2\le W_N/g$, $\E[gK_2\tilde F^2]\le CW_N^2\int T((y^\circ)^2)/g\,dx\le CW_N^2\Sigma_y$.
\item \emph{Degenerate-$U$ terms.} Using $K_2\le W_N/g$ and $p^2\le(W_N/g)\,p$, both terms are at most $W_N\E[T(\Phi_2^2)]\le CW_N\Theta_2$, by Lemma~\ref{lem:A-pop}(b).
\item \emph{Hub term.} By \eqref{eq:A-phiU}, $|\phi(v,x)|\le CW_N(|y^\circ(x)|+\mathfrak y(v))$, so $\mathfrak c_1(x)\le CW_N(|y^\circ(x)|+T\mathfrak y(x))$, and $\E[g\mathfrak c_1^2]\le CW_N^2\Sigma_y$.
\item \emph{Coincidence terms.} $p^4\le(W_N/g)^3p$ gives $\E[g\int(p^4+p^2)e^2]\le CW_N^3\Sigma_y$. Next, $\E[g\int p^2(T\phi_v)^2]\le W_N^2c_0^{-1}\E[T\Phi_2^2]\le CW_N^2\Theta_2$. Finally, $\mathfrak c_2\le(W_N/g)\int p\,\phi(v,x)^2dv$ and \eqref{eq:A-phiU} give $\E[g\mathfrak c_2]\le CW_N^3\Sigma_y$.
\end{itemize}
Thus $\Lambda_N\le CW_N^2$ and $\Lambda'_N\le CW_N^3$, using Lemma~\ref{lem:A-theta}, and Lemma~\ref{lem:A-cond}(c) gives the stated bound on $\Ohgone$.

\emph{Off $\Ohgone$.} Use $g/\hat g_1\le W_N/c_0$, $|Ty^\circ(x)|+T\mathfrak y(x)\le CW_N$ and \eqref{eq:A-phiU}. These give $|\eta_1|\le CW_N^4\big(1+\avg_{j\ne1}(|y^\circ(\xi_j)|+\mathfrak y(\xi_j))\big)$, so $\E\eta_1^4\le CW_N^{16}$, since $\|y^\circ\|_{L^4}+\|\mathfrak y\|_{L^4}\le C$. With $\hat g_1\le W_N\le N^{1/2}$, Cauchy--Schwarz and Lemma~\ref{lem:A-cond}(a) give
\[
\E[\hat g_1\eta_1^2\one_{\Ohgone^c}]\le W_N(\E\eta_1^4)^{1/2}\Pr(\Ohgone^c)^{1/2}\le CN^{9/2}e^{-\dmin/8}\le CN^{9/2}N^{-C_*/8}\le CN^{-3}
\]
for $C_*\ge60$.
\end{proof}

\begin{proof}[Proof of Proposition~\ref{prop:remainder}]
Combine Lemma~\ref{lem:A-transfer} with Propositions~\ref{prop:A-remB} and~\ref{prop:A-remU}. The bounds on $\|\phi\|_{L^2}$ and $\|H_{\rm rest}\|_{L^{2+\delta}}$, and the Lyapunov condition, are Lemmas~\ref{lem:A-pop}, \ref{lem:A-theta} and~\ref{lem:A-lyap}.
\end{proof}

\begin{remark}[Where (U1) and (M3) enter the signal half]\label{rem:A-where}
In setting (U), $W_N$ enters every term of $\Lambda_N$ and $\Lambda'_N$: the product, degenerate-$U$, hub and coincidence terms. It also enters the bad-event bound. In setting (B), the implicit bound $W_N\le1/(2\rho_N)$ of (M3) is used only once in $\Lambda_N$ (the hub term), in $\Lambda'_N$, and on $\Ohgone^c$. Each of these sits inside $r_N^2$, and none is applied at the $\sqrt n$ scale.
\end{remark}

\section{Tools for the network fluctuation}\label{app:tools}

\subsection{Concentration and operators}
The following tail inequality, which holds only for non-negative random variables and only handles one tail, will be useful in later arguments.
\begin{lemma}[Lower tail for nonnegative sums]\label{lem:A-lt}
If $X_1,\dots,X_k\ge0$ are independent with finite second moments, then for $t\ge0$,
\[
\Pr\left(\sum_{i=1}^k X_i\le\E\sum_{i=1}^k X_i-t\right)\le\exp\left\{-\frac{t^2}{2\sum_{i=1}^k\E X_i^2}\right\}.
\]
\end{lemma}

\begin{proof}
Write $S=\sum_iX_i$ and $V=\sum_i\E X_i^2$, and let $\lambda\ge0$.

\emph{Moment generating function.} For $u\ge0$, $e^{-u}\le1-u+u^2/2$, and $1+y\le e^y$. Applying both with $u=\lambda X_i\ge0$,
\[
\E e^{-\lambda X_i}\le1-\lambda\E X_i+\tfrac{\lambda^2}2\E X_i^2\le\exp\big\{-\lambda\E X_i+\tfrac{\lambda^2}2\E X_i^2\big\}.
\]

\emph{Chernoff bound.} By Markov's inequality applied to $e^{-\lambda S}$, and independence,
\[
\Pr(S\le\E S-t)=\Pr\big(e^{-\lambda S}\ge e^{-\lambda(\E S-t)}\big)\le e^{\lambda(\E S-t)}\prod_i\E e^{-\lambda X_i}\le\exp\big\{-\lambda t+\tfrac{\lambda^2}2V\big\}.
\]
Choosing $\lambda=t/V$ gives the claim. See also \citet{maurer2003bound}.
\end{proof}

Let $\mathcal V\subseteq[N]$ be a node set. Given $\F$, let the edges $A_{ik}$, $i<k$ in $\mathcal V$, be independent $\mathrm{Bern}(p_{ik})$. Let $\mu_i=\sum_{k\in\mathcal V\setminus i}p_{ik}$, and suppose $\mu_i\ge \dmin$ for all $i$. This holds for the full graph and for every graph with one node deleted. Write $\bar D=\diag(\mu_i)$.

\begin{lemma}[Concentration]\label{lem:A-conc}
There is an absolute constant $C_\varepsilon$ such that, if $\varepsilon=C_\varepsilon\sqrt{\log N/\dmin}\le1/4$, each of the following events (the \emph{degree}, \emph{spectral} and \emph{sup-norm} events) has conditional probability at least $1-N^{-14}$ given $\F$:
\begin{gather*}
\Odeg:=\Big\{\max_i\Big|\frac{d_i}{\mu_i}-1\Big|\le\varepsilon\Big\},\qquad
\Osp:=\Big\{\varepsilon_*:=\|\bar D^{-1/2}X\bar D^{-1/2}\|_{\rm op}\le\varepsilon\Big\},\\
\Osup:=\Big\{\max_i\Big|\sum_kX_{ik}\big(u_k-(\hat Pu)_i\big)\Big|\le\varepsilon\mu_i\|u\|_\infty\Big\},
\end{gather*}
where in $\Osup$ the vector $u$ is any fixed $\F$-measurable vector. We write $\Ods:=\Odeg\cap\Osp$. When these events are needed for a particular node set, the set is indicated by a superscript: ${}^{\rm full}$ for $[N]$ and ${}^{(0)}$ for $\mathcal V=\{2,\dots,N\}$.
\end{lemma}

\begin{proof}
\emph{Degrees.} Given $\F$, $d_i-\mu_i=\sum_kX_{ik}$ is a sum of independent centered terms with $|X_{ik}|\le1$ and total variance at most $\mu_i$. Scalar Bernstein with $t=\varepsilon\mu_i$ gives a tail exponent at least $\varepsilon^2\mu_i/3\ge\varepsilon^2\dmin/3=C_\varepsilon^2\log N/3$, using $\varepsilon\le1/4$. A union bound over the at most $N$ nodes gives failure probability at most $2N^{1-C_\varepsilon^2/3}$.

\emph{Spectral norm.} Write $\bar D^{-1/2}X\bar D^{-1/2}=\sum_{i<k}X_{ik}(e_ie_k^\top+e_ke_i^\top)/\sqrt{\mu_i\mu_k}$, a sum of independent centered random matrices. Each summand has operator norm at most $1/\sqrt{\mu_i\mu_k}\le1/\dmin$. The variance proxy $\sum_{i<k}\E X_{ik}^2(e_ie_k^\top+e_ke_i^\top)^2/(\mu_i\mu_k)$ is diagonal, with entries $\sum_k\Var(X_{ik})/(\mu_i\mu_k)\le\sum_kp_{ik}/(\mu_i\dmin)\le1/\dmin$. Matrix Bernstein \citep[Theorem~1.4]{tropp2012user}, applied to the sum and its negative, gives
\[
\Pr(\varepsilon_*\ge t\mid\F)\le2N\exp\Big\{-\frac{t^2/2}{(1+t/3)/\dmin}\Big\},
\]
and at $t=\varepsilon\le1/4$ this is at most $2N^{1-3C_\varepsilon^2/7}$.

\emph{Sup-norm event.} For fixed $i$, the summands $X_{ik}(u_k-(\hat Pu)_i)$ are independent and centered, bounded by $c=2\|u\|_\infty$, with total variance at most $\mu_ic^2$. Bernstein at $t=\varepsilon\mu_ic/2$ gives exponent at least $\varepsilon^2\mu_i/9$, and a union bound over $i$ finishes.

\emph{Constants.} Taking $C_\varepsilon^2\ge144$ makes each probability at most $N^{-14}$; then $\varepsilon\le1/4$ as soon as $C_*\ge16C_\varepsilon^2$.
\end{proof}

The network component is $V=Y-Y^W$, the difference between the responses of the
actual and the expected-weight walk. The next lemma records the facts about these
two walks that are used throughout. First, $V$ is the resolvent of the actual walk
applied to the perturbation $\cE=P-\hat P$ acting on $Y^W$, and $\cE$ can be written
in terms of the centered edges $X=A-\E[A\mid\F]$ alone. Second, on the degree event
$\Odeg$ of Lemma~\ref{lem:A-conc}, every degree is within a constant factor of its expected
value, so $P$ is reversible with respect to the degrees and contracts the same
weighted norms as $\hat P$. Third, on the spectral and sup-norm events, $\cE$ is
small: of order $\varepsilon\asymp(\log N/(N\rho_N))^{1/2}$ in $\ell^2(\mu)$, and
uniformly on fixed $\F$-measurable vectors.

\begin{lemma}[Algebra, and operators on the concentration events]\label{lem:A-ops}
Everywhere,
\begin{equation}\label{eq:A-Ealg}
(\cE v)_i=\frac1{d_i}\sum_jX_{ij}\big(v_j-(\hat Pv)_i\big)\quad(d_i\ge1),\qquad V=R_P\Gamma\cE Y^W,
\end{equation}
and $\|\Gamma P\|_{\infty\to\infty},\|\Gamma\hat P\|_{\infty\to\infty}\le\bar\gamma$, so that $\|Y\|_\infty,\|Y^W\|_\infty\le\max_i|Z_i|/(1-\bar\gamma)$. On $\Odeg$:
\begin{itemize}
\item[(a)] $\mu_i/2\le d_i\le3\mu_i/2$, and $P$ is stochastic and reversible with respect to $d$;
\item[(b)] $P\Gamma$ and $\Gamma P$ contract $\ell^\infty$ and every $\ell^s(d)$, $s\ge1$, by the factor $\bar\gamma$. Likewise $\hat P\Gamma$ and $\Gamma\hat P$ contract $\ell^\infty$ and $\ell^s(\mu)$.
\end{itemize}
On $\Ods=\Odeg\cap\Osp$, $\|\cE v\|_{\ell^2(\mu)}\le4\varepsilon\|v\|_{\ell^2(\mu)}$. On $\Odeg\cap\Osup$, $|(\cE u)_i|\le2\varepsilon\|u\|_\infty$.
\end{lemma}

\begin{proof}
\emph{The identities \eqref{eq:A-Ealg}.} For $d_i\ge1$, $(Pv)_i-(\hat Pv)_i=d_i^{-1}\sum_jA_{ij}v_j-(\hat Pv)_i=d_i^{-1}\sum_jA_{ij}\big(v_j-(\hat Pv)_i\big)$. Since $\sum_jp_{ij}\big(v_j-(\hat Pv)_i\big)=\mu_i(\hat Pv)_i-\mu_i(\hat Pv)_i=0$, $A$ may be replaced by $X=A-p$. Next, $Y=R_PZ$ and $Y^W=MZ$, so
\[
Y-Y^W=R_P\big[(I-\Gamma\hat P)-(I-\Gamma P)\big]Y^W=R_P\Gamma(P-\hat P)Y^W=R_P\Gamma\cE Y^W .
\]
The $\ell^\infty$ bounds are \ref{F:ops}(i)--(iii); for $P$ they use that rows of $P$ sum to at most one, which holds everywhere.

\emph{(a) and (b).} On $\Odeg$, $|d_i-\mu_i|\le\varepsilon\mu_i\le\mu_i/4$, which gives (a); in particular no node is isolated, so $P$ is stochastic, and $d_iP_{ik}=A_{ik}$ is symmetric. Part (b) is \ref{F:ops}(ii) with $\nu=d$ for $P$, and with $\nu=\mu$ for $\hat P$, since $\mu_i\hat P_{ik}=p_{ik}$ is symmetric.

\emph{The spectral bound.} By \eqref{eq:A-Ealg}, $\cE v=D^{-1}Xv-D^{-1}\diag(X\one)\hat Pv$. For the first term, using $d_i\ge\mu_i/2$,
\[
\|D^{-1}Xv\|^2_{\ell^2(\mu)}=\sum_i\frac{\mu_i}{d_i^2}(Xv)_i^2\le4\sum_i\frac{(Xv)_i^2}{\mu_i}=4\big\|(\bar D^{-1/2}X\bar D^{-1/2})\bar D^{1/2}v\big\|^2\le4\varepsilon_*^2\|v\|^2_{\ell^2(\mu)} .
\]
For the second term, $|(X\one)_i|/d_i=|d_i-\mu_i|/d_i\le\varepsilon\mu_i/d_i\le2\varepsilon$, and $\hat P$ contracts $\ell^2(\mu)$. So on $\Ods$, $\|\cE v\|_{\ell^2(\mu)}\le(2\varepsilon_*+2\varepsilon)\|v\|_{\ell^2(\mu)}\le4\varepsilon\|v\|_{\ell^2(\mu)}$.

\emph{The sup-norm bound.} On $\Odeg\cap\Osup$, by \eqref{eq:A-Ealg}, $|(\cE u)_i|\le d_i^{-1}\varepsilon\mu_i\|u\|_\infty\le2\varepsilon\|u\|_\infty$.
\end{proof}

The spectral and sup-norm events of Lemma~\ref{lem:A-conc} hold only with
$\varepsilon\asymp(\log N/(N\rho_N))^{1/2}$, which is where a logarithm enters.
When the perturbation acts on an $\F$-measurable vector, its moments can instead be
computed directly, conditionally on $\F$, using only the degree event. The next lemma
does this, and gives bounds of order $(N\rho_N)^{-s/2}$ with no logarithm.

\begin{lemma}[Graph-noise moments on the degree event alone]\label{lem:A-noise}
In the setting of Lemma~\ref{lem:A-conc}, let $n_{\mathcal V}=|\mathcal V|-1$ and $\hat g^{\mathcal V}_i=\mu_i/(\rho n_{\mathcal V})\ge c_0$. Let $Y^{\mathcal V}$ and $y^{\mathcal V}$ be the responses on $\mathcal V$, set $V^{\mathcal V}=Y^{\mathcal V}-y^{\mathcal V}$, and set $\gn^{\mathcal V}=PV^{\mathcal V}+\cE y^{\mathcal V}$. For $s\in\{2,4\}$,
\[
\E\Big[\one_{\Odeg}\sum_{i\in\mathcal V}\hat g^{\mathcal V}_i\big(|V^{\mathcal V}_i|^s+|(\cE y^{\mathcal V})_i|^s+|\gn^{\mathcal V}_i|^s\big)\Bigm|\F\Big]\le C\,\dmin^{-s/2}\sum_{i\in\mathcal V}\hat g^{\mathcal V}_i|Z_i|^s .
\]
Also, on $\Odeg$, $\sum_i\hat g^{\mathcal V}_i|Y^{\mathcal V}_i|^4\le C\sum_i\hat g^{\mathcal V}_iZ_i^4$.
\end{lemma}

\begin{proof}
Drop the superscript $\mathcal V$, and write $n:=n_{\mathcal V}$ within this proof. On $\Odeg$, $d_i\ge\mu_i/2$ and $d_i\le\frac32\mu_i$, so
\begin{equation}\label{eq:A-dg}
\hat g_i\le\frac{2d_i}{\rho n},\qquad d_i\le\tfrac32\rho n\,\hat g_i .
\end{equation}

\emph{Step 1: moments of $\cE y$ given $\F$.} By \eqref{eq:A-Ealg}, $(\cE y)_k=d_k^{-1}\sum_jX_{kj}c_{kj}$ with $c_{kj}=y_j-(\hat Py)_k$, which is $\F$-measurable. Given $\F$ the $X_{kj}$ are independent, centered, with $|X_{kj}|\le1$ and $\Var(X_{kj})\le p_{kj}=\mu_k\hat P_{kj}$, so
\begin{gather*}
\E\Big[\Big(\sum_jX_{kj}c_{kj}\Big)^2\Bigm|\F\Big]\le\mu_k(\hat Pc_k^2)_k,\\
\E\Big[\Big(\sum_jX_{kj}c_{kj}\Big)^4\Bigm|\F\Big]\le\mu_k(\hat Pc_k^4)_k+3\mu_k^2(\hat Pc_k^2)_k^2\le4\mu_k^2(\hat Pc_k^4)_k ,
\end{gather*}
using Jensen and $\mu_k\ge1$ in the last step. On $\Odeg$, $d_k\ge\mu_k/2$, so for $s\in\{2,4\}$,
\[
\E\big[\one_{\Odeg}d_k|(\cE y)_k|^s\bigm|\F\big]\le C\mu_k^{1-s/2}(\hat P|c_k|^s)_k\le C\dmin^{1-s/2}(\hat P|y|^s)_k ,
\]
where the last step uses $|c_{kj}|^s\le2^{s-1}(|y_j|^s+|(\hat Py)_k|^s)$ and $|(\hat Py)_k|^s\le(\hat P|y|^s)_k$.

\emph{Step 2: from $d$-weights to $\hat g$-weights.} By \eqref{eq:A-dg}, $\hat g_k|(\cE y)_k|^s\le\frac2{\rho n}d_k|(\cE y)_k|^s$. Summing Step~1 over $k$,
\[
\E\Big[\one_{\Odeg}\sum_k\hat g_k|(\cE y)_k|^s\Bigm|\F\Big]\le\frac{C\dmin^{1-s/2}}{\rho n}\sum_k(\hat P|y|^s)_k .
\]
Since $\hat g\ge c_0$ and $\hat g$ is stationary for $\hat P$ ($\hat g_k\hat P_{kj}$ is symmetric), $\sum_k(\hat P|y|^s)_k\le c_0^{-1}\sum_k\hat g_k(\hat P|y|^s)_k=c_0^{-1}\sum_j\hat g_j|y_j|^s$. By \ref{F:ops}(ii)--(iii) with $\nu=\hat g$, $\sum\hat g|y|^s\le(1-\bar\gamma)^{-s}\sum\hat g|Z|^s$. Finally $1/(\rho n)\le c_0/\dmin$, which gives the bound for $\cE y$.

\emph{Step 3: $V$ and $\gn$.} Since $V=R_P\Gamma\cE y$ and, on $\Odeg$, $R_P\Gamma$ is bounded on $\ell^s(d)$ by $\bar\gamma/(1-\bar\gamma)$ (\ref{F:ops}), $\sum d|V|^s\le C\sum d|\cE y|^s$. Converting both sides with \eqref{eq:A-dg}, $\sum\hat g|V|^s\le\frac2{\rho n}\sum d|V|^s\le\frac{C}{\rho n}\sum d|\cE y|^s\le C\sum\hat g|\cE y|^s$, and Step~2 applies. For $\gn=PV+\cE y$, the $\ell^s(d)$ contraction of $P$ gives the same bound.

\emph{The bound on $Y$.} On $\Odeg$, $Y=R_PZ$, and the $\ell^4(d)$ contraction of $\Gamma P$ gives $\sum d|Y|^4\le C\sum dZ^4$; \eqref{eq:A-dg} converts both sides to $\hat g$-weights.
\end{proof}

\begin{corollary}[Single-node graph noise]\label{cor:A-V1}
Under (M1)--(M4),
\[
\E V_1^2\le\frac{C\,\E[gZ^2]}{\dmin}+CN^{-13}\E Z^2,\qquad\E V_1^4\le\frac{C\,\E[gZ^4]}{\dmin^2}+CN^{-13}\E Z^4 .
\]
\end{corollary}

\begin{proof}
Fix $s\in\{2,4\}$. By exchangeability, $\E V_1^s=\E[N^{-1}\sum_iV_i^s]$, and we split according to $\Odeg$ for the full graph.

\emph{On $\Odeg$.} Since $\hat g_i\ge c_0$, $N^{-1}\sum_i|V_i|^s\le(c_0N)^{-1}\sum_i\hat g_i|V_i|^s$. Lemma~\ref{lem:A-noise} with $\mathcal V=[N]$ bounds its conditional expectation by $C\dmin^{-s/2}N^{-1}\sum_i\hat g_i|Z_i|^s$, and $\E[\hat g_1|Z_1|^s]=\E[g|Z|^s]$ since $\E[\hat g_1\mid\xi_1]=g(\xi_1)$.

\emph{Off $\Odeg$.} Everywhere, $|V_i|\le|Y_i|+|Y^W_i|\le2\max_j|Z_j|/(1-\bar\gamma)$, and $\max_j|Z_j|^s\le\sum_j|Z_j|^s$. Since $\Pr(\Odeg^c\mid\F)\le N^{-14}$ and $\sum_j|Z_j|^s$ is $\F$-measurable, this part contributes at most $CN^{-14}\cdot N\E|Z|^s=CN^{-13}\E|Z|^s$.
\end{proof}

In setting (U), the change caused by resampling node $1$ involves products of sums
over node $1$'s edges, such as $\big(\sum_mA_{m1}x_m\big)^2\big(\sum_mA_{m1}z_m\big)^2$,
where $x$ carries the dual weights and $z$ the graph noise at the neighbors. Given
everything except node $1$'s block, and given $\xi_1$, the edges $A_{m1}$ are
independent Bernoulli variables, but the weights are unbounded, so these moments
cannot be controlled through sup norms. The next lemma expands them exactly. Each
term corresponds to a pattern of coincidences among node $1$'s neighbors (all
distinct, two equal, and so on) and is a product of sums
$s_{ij}=\sum_mp_{m1}x_m^iz_m^j$. Each pattern is then averaged over the type of
node $1$ and bounded by one of the functionals in (U4); see Table~\ref{tab:A-patterns}.
\begin{lemma}[Products of Bernoulli sums]\label{lem:A-bern}
Let $A_m\sim\mathrm{Bern}(p_m)$ be independent, $X_m=A_m-p_m$, $x,z\ge0$ and $a$ real. Write $s_{ij}=\sum p_mx_m^iz_m^j$ and $t_{ij}=\sum p_mx_m^i|a_m|^j$. Then
\begin{align*}
\E\Big[\Big(\sum Ax\Big)^2\Big(\sum Az\Big)^2\Big]&\le s_{10}^2s_{01}^2+s_{20}s_{01}^2+s_{10}^2s_{02}+s_{20}s_{02}+4s_{11}s_{10}s_{01}\\
&\qquad+2s_{11}^2+2s_{21}s_{01}+2s_{12}s_{10}+s_{22},\\
\E\Big[\Big(\sum Ax\Big)^2\Big(\sum Xa\Big)^2\Big]&\le s_{10}^2t_{02}+s_{20}t_{02}+2t_{11}^2+2t_{12}s_{10}+t_{22},\\
\E\Big[\Big(\sum Ax\Big)^2\sum Az\Big]&\le s_{10}^2s_{01}+s_{20}s_{01}+2s_{11}s_{10}+s_{21},\\
\E\Big[\sum Ax\sum Az\Big]&\le s_{10}s_{01}+s_{11}.
\end{align*}
\end{lemma}

\begin{proof}
Expand into index tuples and group them by the set partition of the positions induced by equal indices. Bound the sum over assignments of distinct indices to blocks by the unrestricted sum. For blocks consisting only of $A$'s, $\E\prod A=p$. For blocks containing an $X$, $|\E[A^iX^j]|\le p$ when $i+j\ge2$, and $\E X=0$. The first display has the $15$ partitions of four positions, with the multiplicities shown: for instance, $s_{10}^2s_{01}^2$ is the pattern in which all four indices are distinct, $4s_{11}s_{10}s_{01}$ collects the four ways in which exactly one $x$-index coincides with one $z$-index, and $s_{22}$ is the pattern in which all four indices coincide. Each coincidence replaces a product of two first-moment sums by a single sum with a higher power, and it is these higher powers that, after averaging over the type of node $1$, produce the functionals of (U4). In the second, partitions in which an $X$-position is a singleton vanish, which leaves the five listed types. The remaining displays are the partitions of three and two positions.
\end{proof}

\subsection{Leave-one-node-out structure}
\begin{table}[htbp]
\centering\small
\begin{tabular}{p{0.22\linewidth}p{0.55\linewidth}p{0.15\linewidth}}
\toprule
Symbol & Meaning & Defined in\\
\midrule
$\mathcal V$, $n_0$ & nodes other than node $1$; $N-2$ & App.~\ref{app:tools}\\
superscript ${}^{(0)}$ & object of the graph with node $1$ deleted ($d^{(0)}$, $P^{(0)}$, $Y^{(0)}$, \dots) & App.~\ref{app:tools}\\
$\hat g^{(0)}_m$, $\mu^{(0)}_m$ & expected-degree ratio and expected degree in the reduced graph & App.~\ref{app:tools}\\
$\Hc_1$ & everything except node $1$'s block & App.~\ref{app:tools}\\
$\Odeg^{(0)}$, $\Ods^{(0)}$, $\onedeg$ & degree event and degree-and-spectral event for the reduced graph; indicator of $\Odeg^{(0)}$ & Lem.~\ref{lem:A-conc}\\
$\tilde S$ & the sample as seen from the resampled node: $1\notin\tilde S$, $|\tilde S|\le n$ & \eqref{eq:ES}\\
$G$ & effect of node $1$ on $\sum_{i\in\tilde S}V_i$: full graph minus graph with node $1$ removed & Lem.~\ref{lem:A-loo}\\
$\lambda$, $\hat\ell$ & sensitivity of the sample sum to the input at each node; expected-weight version & Lem.~\ref{lem:A-loo}\\
$\tilde Y_m$, $\hat y_m$ & what node $1$ transmits to its neighbor $m$; expected-weight version & Lem.~\ref{lem:A-loo}\\
$\src$, $\hat\src$ & sources $(D^{(0)})^{-1}\one_{\tilde S}$ and $(\bar D^{(0)})^{-1}\one_{\tilde S}$ & Lem.~\ref{lem:A-dual}\\
$b$, $\hat b$ & dual vectors: degree-free sensitivities, $\lambda=D^{(0)}b$ & Lem.~\ref{lem:A-dual}\\
$\beta$, $\hat\beta$, $c_m$ & edge-level weights in $G$; bound on their discrepancy & Lem.~\ref{lem:A-dual}\\
$\tau$, $\Eg{\tau}$, $\Eg{\tau;\Omega}$ & weight averaged with the dual weights; its degree-biased mean $\E[\hat g^{(0)}_2\tau_2]$, without or with the indicator of a good event $\Omega$ & Lem.~\ref{lem:weight-transfer}\\
$r_m$ & degree ratio $g(\xi_m)/\hat g^{(0)}_m$ & Lem.~\ref{lem:A-two}\\
\bottomrule
\end{tabular}
\caption{Notation for the leave-one-node-out argument.}
\label{tab:not-loo}
\end{table}

By exchangeability we resample node $1$. Let $\mathcal V=\{2,\dots,N\}$ and $n_0=N-2$. Superscript ${(0)}$ marks objects of the graph with node $1$ deleted: reduced degrees $d^{(0)}$, $P^{(0)}$ with row and column $1$ zero, $R^{(0)}=(I-\Gamma P^{(0)})^{-1}$, $Y^{(0)}=R^{(0)}Z$, $\mu^{(0)}_m=\mu_m-p_{m1}=\rho n_0\hat g^{(0)}_m$ with $\hat g^{(0)}_m=n_0^{-1}\sum_{k\in\mathcal V\setminus m}w_{mk}$, $\hat P^{(0)}_{mk}=p_{mk}/\mu^{(0)}_m$, $M^{(0)}=(I-\Gamma\hat P^{(0)})^{-1}$, $y^{(0)}=M^{(0)}Z$, $V^{(0)}=Y^{(0)}-y^{(0)}$ and $\cE^{(0)}=P^{(0)}-\hat P^{(0)}$. Let $\Hc_1:=\sigma(\xi_k,k\ne1;\ U_{ik},i,k\ne1)$, where $U_{ik}$ are the edge uniforms of Definition~\ref{def:graphon}. Given $(\Hc_1,\xi_1)$, the edges $(A_{m1})_m$ are independent $\mathrm{Bern}(p_{m1})$. The events of Lemma~\ref{lem:A-conc} for the reduced graph are $\Odeg^{(0)}$, $\Osp^{(0)}$ and $\Ods^{(0)}=\Odeg^{(0)}\cap\Osp^{(0)}$; all three belong to $\Hc_1$. We abbreviate $\onedeg:=\one_{\Odeg^{(0)}}$.

\begin{lemma}[Reduction to $\Var\bar V_S$]\label{lem:A-suff}
$\bar Y_S-\mu_N=(\bar Y^W_S-\E\bar Y^W_S)+(\bar V_S-\E\bar V_S)$. Also $\Var(\avg_S\Delta_i)\le\Var\bar V_S$ and $\E\Var(\bar Y_S\mid\F)\le\Var\bar V_S$, where $\Delta_i:=\E[Y_i\mid\F]-Y^W_i$.
\end{lemma}

\begin{proof}
The decomposition holds because $Y=Y^W+V$ and $\E\bar Y_S=\mu_N$. For the inequalities, $Y^W$ is $\F$-measurable, so $\E[\bar V_S\mid\F]=\avg_S\Delta_i$ and $\Var(\bar Y_S\mid\F)=\Var(\bar V_S\mid\F)$. The law of total variance, $\Var\bar V_S=\Var(\E[\bar V_S\mid\F])+\E\Var(\bar V_S\mid\F)$, then gives both.
\end{proof}

\begin{lemma}[Leave-one-node-out identity]\label{lem:A-loo}
For $\tilde S\subseteq\mathcal V$, let $\lambda=(R^{(0)})^\top\one_{\tilde S}$ and $\hat\ell=(M^{(0)})^\top\one_{\tilde S}$. Both are $\Hc_1$-measurable. Everywhere,
\[
G:=\sum_{i\in\tilde S}\big(V_i-V^{(0)}_i\big)=\sum_{m\in\mathcal V}\gamma_m\Big[A_{m1}\beta_m\tilde Y_m-p_{m1}\hat\beta_m\hat y_m\Big],
\]
where $\beta_m=\lambda_m/(d^{(0)}_m+1)$, $\hat\beta_m=\hat\ell_m/\mu_m$, $\tilde Y_m=Y_1-(P^{(0)}Y)_m$ and $\hat y_m=Y^W_1-(\hat P^{(0)}Y^W)_m$.
\end{lemma}

\begin{proof}
\emph{The actual response.} From $Y=Z+\Gamma PY$ and $Y^{(0)}=Z+\Gamma P^{(0)}Y^{(0)}$, we get $(I-\Gamma P^{(0)})(Y-Y^{(0)})=\Gamma(P-P^{(0)})Y$, so $Y-Y^{(0)}=R^{(0)}\Gamma(P-P^{(0)})Y$. Column $1$ of $P^{(0)}$ vanishes, so $R^{(0)}e_1=e_1$, and the row-$1$ entry of $(P-P^{(0)})Y$ affects only coordinate $1\notin\tilde S$. For $m\ne1$, adding neighbor $1$ to $m$'s neighborhood changes the average by
\[
(PY)_m-(P^{(0)}Y)_m=\frac{A_{m1}}{d_m}\big(Y_1-(P^{(0)}Y)_m\big),
\]
since the new average is a weighted combination of the old one, with weight $d^{(0)}_m/d_m$, and $Y_1$, with weight $1/d_m$. When $A_{m1}=1$, $d_m=d^{(0)}_m+1$. This includes the case $d^{(0)}_m=0$, where $(P^{(0)}Y)_m=0$ by the zero-row convention and $(PY)_m=Y_1$. Summing over $i\in\tilde S$ and using $\sum_{i\in\tilde S}(R^{(0)})_{im}=\lambda_m$,
\[
\sum_{i\in\tilde S}(Y_i-Y^{(0)}_i)=\sum_m\lambda_m\gamma_m\frac{A_{m1}}{d^{(0)}_m+1}\tilde Y_m=\sum_m\gamma_mA_{m1}\beta_m\tilde Y_m .
\]

\emph{The expected-weight response.} The same computation with the weights $p_{m1}$ in place of $A_{m1}$, and $\mu_m=\mu^{(0)}_m+p_{m1}$ in place of $d_m$, gives $\sum_{i\in\tilde S}(Y^W_i-y^{(0)}_i)=\sum_m\gamma_mp_{m1}\hat\beta_m\hat y_m$. Subtracting the two displays gives the identity for $G$.
\end{proof}

\begin{proof}[Proof of \eqref{eq:ES} (vertex-block Efron--Stein)]
\emph{The blocks.} The blocks $\mathcal B_j=(\xi_j,(U_{jk})_{k<j})$, $j\in[N]$, are independent and determine $(\xi,A)$. Replacing $\mathcal B_j$ by an independent copy leaves $\xi_k$ for $k\ne j$, and all edges not incident to $j$, unchanged.

\emph{The change at one block.} Let $f=\sum_{i\in S}V_i$, and let $f^{(j)}$ be $f$ computed with $\mathcal B_j$ redrawn. Write $V^{(0,j)}$ for the network component of the graph with node $j$ deleted. It does not depend on $\mathcal B_j$: the reduced graph involves neither $\xi_j$ nor any edge at $j$, and since $(R^{(0,j)})_{ij}=\delta_{ij}$, $Z_j$ never enters $V^{(0,j)}_i$ for $i\ne j$. Splitting off node $j$,
\[
f=\one\{j\in S\}V_j+\sum_{i\in S\setminus j}V^{(0,j)}_i+G_j,\qquad G_j:=\sum_{i\in S\setminus j}\big(V_i-V^{(0,j)}_i\big),
\]
so $f-f^{(j)}=\one\{j\in S\}(V_j-V_j')+(G_j-G_j')$, where primes denote the values after redrawing; $V'_j$ and $G'_j$ have the same laws as $V_j$ and $G_j$.

\emph{Efron--Stein.} Using $(a+b)^2\le2a^2+2b^2$, $\E(V_j-V'_j)^2\le4\E V_j^2=4\E V_1^2$ and $\E(G_j-G'_j)^2\le4\E G_j^2$,
\[
\Var f\le\frac12\sum_j\E(f-f^{(j)})^2\le4n\,\E V_1^2+4\sum_{j=1}^N\E G_j^2 .
\]
Dividing by $n^2$ and multiplying by $n$, $n\Var\bar V_S\le4\E V_1^2+\frac4n\sum_j\E G_j^2$.

\emph{Exchangeability.} Relabeling the nodes so that $j$ becomes node $1$ does not change the joint law, and turns $S\setminus j$ into a set $\tilde S\not\ni1$ with $|\tilde S|\le n$, and $G_j$ into the quantity $G$ of Lemma~\ref{lem:A-loo}. Hence $\E G_j^2\le\sup_{\tilde S}\E G^2$ for every $j$, and \eqref{eq:ES} follows.
\end{proof}

The leave-one-node-out identity expresses the change in $\sum_{i\in\tilde S}V_i$
through the weights $\lambda=(R^{(0)})^\top\one_{\tilde S}$: $\lambda_m$ is the change in
the sample sum when the input at node $m$ increases by one, and collects the
discounted walks that start in $\tilde S$ and end at $m$. By detailed balance,
$d_i(P^t)_{im}=d_m(P^t)_{mi}$, and therefore $\lambda=D^{(0)}b$ with
$b=(I-P^{(0)}\Gamma)^{-1}\src$, $\src=(D^{(0)})^{-1}\one_{\tilde S}$. The entry $b_m$ is a walk
expectation: over walks started at $m$, the product of the coefficients $\gamma$ along
the walk times $\one\{\text{endpoint}\in\tilde S\}/d_{\text{endpoint}}$. We work with $b$
rather than $\lambda$ because $b$ is the degree-free part of the influence: it is
uniformly of order $1/\dmin$, and it contracts under the walk in every degree-weighted norm.
The expected-weight dual vector $\hat b$, built in the same way from $\hat P^{(0)}$,
differs from $b$ only through the degree fluctuations and the perturbation
$\cE^{(0)}$.
\begin{lemma}[Dual vectors]\label{lem:A-dual}
On $\Odeg^{(0)}$, set $\src=(D^{(0)})^{-1}\one_{\tilde S}$ and $\hat\src=(\bar D^{(0)})^{-1}\one_{\tilde S}$. Then $\lambda=D^{(0)}b$ with $b=(I-P^{(0)}\Gamma)^{-1}\src$, and $\hat\ell=\bar D^{(0)}\hat b$ with $\hat b=(I-\hat P^{(0)}\Gamma)^{-1}\hat\src$. Moreover:
\begin{itemize}
\item[(i)] $\sum_md^{(0)}_m|b_m|\le\frac{n}{1-\bar\gamma}$ and $\sum_m\mu^{(0)}_m|\hat b_m|\le\frac n{1-\bar\gamma}$;
\item[(ii)] $\sum_md^{(0)}_mb_m^2\le\frac{2n}{(1-\bar\gamma)^2\dmin}$ and $\sum_m\mu^{(0)}_m\hat b_m^2\le\frac{n}{(1-\bar\gamma)^2\dmin}$;
\item[(iii)] $\|b\|_\infty\le\frac2{(1-\bar\gamma)\dmin}$ and $\|\hat b\|_\infty\le\frac1{(1-\bar\gamma)\dmin}$;
\item[(iv)] on $\Ods^{(0)}$, $\sum_md^{(0)}_m(b_m-\hat b_m)^2\le C_1\varepsilon^2n/\dmin$;
\item[(v)] $|\beta_m|\le|b_m|$ and $|\beta_m-\hat\beta_m|\le c_m:=|b_m-\hat b_m|+\frac2{\dmin}\big(|b_m|+|\hat b_m|\big)$.
\end{itemize}
All of $b$, $\hat b$ and $c$ are $\Hc_1$-measurable, and $\hat b$ and $\hat\src$ are $\F$-measurable.
\end{lemma}

\begin{proof}
Write $d=d^{(0)}$, $\mu=\mu^{(0)}$, $P=P^{(0)}$, $\hat P=\hat P^{(0)}$ and $\cE=\cE^{(0)}$ within this proof.

\emph{The representation.} On $\Odeg^{(0)}$ every reduced degree is positive, and $d_mP_{mk}=A_{mk}$ is symmetric, so $P^\top=DPD^{-1}$. Since $\Gamma$ and $D$ are diagonal and commute, $(I-\Gamma P)^{-\top}=(I-P^\top\Gamma)^{-1}=D(I-P\Gamma)^{-1}D^{-1}$, and therefore
\[
\lambda=(I-\Gamma P)^{-\top}\one_{\tilde S}=D(I-P\Gamma)^{-1}D^{-1}\one_{\tilde S}=Db .
\]
The same computation with $\mu_m\hat P_{mk}=p_{mk}$ symmetric gives $\hat\ell=\bar D\hat b$.

\emph{(i) and (ii).} These follow from \ref{F:ops}(ii)--(iii) with $\nu=d$ for $b$ and $\nu=\mu$ for $\hat b$: the norm of $b$ in $\ell^s(d)$ is at most $(1-\bar\gamma)^{-1}$ times that of $\src$. For the sources, $d_k|\src_k|=\one\{k\in\tilde S\}$ gives $\|\src\|_{\ell^1(d)}=|\tilde S|\le n$, and $d_k\src_k^2=\one\{k\in\tilde S\}/d_k$ with $d_k\ge\dmin/2$ gives $\|\src\|^2_{\ell^2(d)}\le2n/\dmin$. For $\hat\src$, the same holds with $\mu_k\ge\dmin$, which gives the constants $n$ and $n/\dmin$.

\emph{(iii).} By \ref{F:ops}(i) and (iii), $\|b\|_\infty\le\|\src\|_\infty/(1-\bar\gamma)$ and
$\|\hat b\|_\infty\le\|\hat\src\|_\infty/(1-\bar\gamma)$, where
$\|\src\|_\infty=\max_{k\in\tilde S}1/d_k\le2/\dmin$ on $\Odeg^{(0)}$ and
$\|\hat\src\|_\infty=\max_{k\in\tilde S}1/\mu_k\le1/\dmin$.

\emph{(iv).} By definition, $(I-P\Gamma)b=\src$ and $(I-\hat P\Gamma)\hat b=\hat\src$; the second can be rewritten as $(I-P\Gamma)\hat b=\hat\src-\cE\Gamma\hat b$. Subtracting,
\[
(I-P\Gamma)(b-\hat b)=(\src-\hat\src)+\cE\Gamma\hat b ,
\]
so by \ref{F:ops}(iii), $\|b-\hat b\|_{\ell^2(d)}\le(1-\bar\gamma)^{-1}\big(\|\src-\hat\src\|_{\ell^2(d)}+\|\cE\Gamma\hat b\|_{\ell^2(d)}\big)$. For the first term, on $\Odeg^{(0)}$, $|\mu_k-d_k|\le\varepsilon\mu_k$, so
\[
\|\src-\hat\src\|^2_{\ell^2(d)}=\sum_{k\in\tilde S}d_k\Big(\frac1{d_k}-\frac1{\mu_k}\Big)^2=\sum_{k\in\tilde S}\frac{(\mu_k-d_k)^2}{d_k\mu_k^2}\le\varepsilon^2\sum_{k\in\tilde S}\frac1{d_k}\le\frac{2\varepsilon^2n}{\dmin}.
\]
For the second, on $\Ods^{(0)}$, Lemma~\ref{lem:A-ops} gives $\|\cE\Gamma\hat b\|_{\ell^2(\mu)}\le4\varepsilon\|\hat b\|_{\ell^2(\mu)}\le C\varepsilon(n/\dmin)^{1/2}$ by (ii), and $\ell^2(d)$ and $\ell^2(\mu)$ norms differ by at most a factor $\sqrt{3/2}$ on $\Odeg^{(0)}$.

\emph{(v).} Since $\lambda_m=d_mb_m$, $\beta_m=\lambda_m/(d_m+1)=b_m-b_m/(d_m+1)$, so $|\beta_m|\le|b_m|$. Similarly, with $\mu^{\rm full}_m=\mu_m+p_{m1}$, $\hat\beta_m=\hat\ell_m/\mu^{\rm full}_m=\hat b_m-\hat b_mp_{m1}/\mu^{\rm full}_m$. Hence
\[
|\beta_m-\hat\beta_m|\le|b_m-\hat b_m|+\frac{|b_m|}{d_m+1}+\frac{p_{m1}|\hat b_m|}{\mu^{\rm full}_m}\le|b_m-\hat b_m|+\frac2{\dmin}|b_m|+\frac1{2\dmin}|\hat b_m|,
\]
using $d_m+1\ge\dmin/2$, $\mu^{\rm full}_m\ge\dmin$ and $p_{m1}\le1/2$.

\emph{Measurability.} $b$, $\hat b$ and $c$ are functions of the reduced graph and the types other than $\xi_1$, hence $\Hc_1$-measurable; $\hat b$ and $\hat\src$ depend only on the types.
\end{proof}

The term $G_2$ involves the discrepancy $b-\hat b$ between the dual vectors of the actual
and the expected-weight graph. Lemma~\ref{lem:A-dual}(iv) bounds its size on the
spectral event, at the cost of a factor $\varepsilon^2\asymp\log N/(N\rho_N)$. The
logarithm can be avoided because the expected-weight objects $\hat b$ and $\hat\src$
are $\F$-measurable. The discrepancy is then driven by two quantities whose conditional
second moments given $\F$ can be computed directly: the degree fluctuations
$d_k-\mu_k$, with conditional variance at most $\mu_k$, and the perturbation
$\cE^{(0)}$ acting on the fixed vector $\Gamma\hat b$. This is the dual counterpart of
Lemma~\ref{lem:A-noise}. It gives a size of order $n/\dmin^2$ rather than
$\varepsilon^2n/\dmin$, and is used whenever the factor multiplying the size is itself
$\F$-measurable, as in Proposition~\ref{prop:A-G2}.

\begin{lemma}[Dual-vector discrepancy given the types]\label{lem:A-energy}
The spectral event is not needed if the discrepancy is bounded in conditional expectation given $\F$:
\begin{gather*}
\E\Big[\one_{\Odeg^{(0)}}\sum_md^{(0)}_m(b_m-\hat b_m)^2\Bigm|\F\Big]\le\frac{Cn}{\dmin^2},\\
\E\Big[\one_{\Odeg^{(0)}}\sum_md^{(0)}_mc_m^2\Bigm|\F\Big]\le\frac{Cn}{\dmin^2}.
\end{gather*}
\end{lemma}

\begin{proof}
Write $d=d^{(0)}$, $\mu=\mu^{(0)}$, and so on, as in the previous proof.

\emph{Reduction.} By the identity $(I-P\Gamma)(b-\hat b)=(\src-\hat\src)+\cE\Gamma\hat b$ from the proof of Lemma~\ref{lem:A-dual}(iv), the $\ell^2(d)$ contraction, and $(a+b)^2\le2a^2+2b^2$, on $\Odeg^{(0)}$
\[
\sum d(b-\hat b)^2\le2(1-\bar\gamma)^{-2}\big(\|\src-\hat\src\|^2_{\ell^2(d)}+\|\cE\Gamma\hat b\|^2_{\ell^2(d)}\big).
\]
Both terms are now bounded in conditional expectation given $\F$, without the spectral event.

\emph{The sources.} On $\Odeg^{(0)}$, $d_k\ge\mu_k/2$, so $\|\src-\hat\src\|^2_{\ell^2(d)}=\sum_{\tilde S}(\mu_k-d_k)^2/(d_k\mu_k^2)\le2\sum_{\tilde S}(\mu_k-d_k)^2/\mu_k^3$. Given $\F$, $\E[(d_k-\mu_k)^2\mid\F]=\Var(d_k\mid\F)\le\mu_k$, so
\[
\E\big[\onedeg\,\|\src-\hat\src\|^2_{\ell^2(d)}\bigm|\F\big]\le2\sum_{k\in\tilde S}\frac1{\mu_k^2}\le\frac{2n}{\dmin^2}.
\]

\emph{The perturbation.} Let $v=\Gamma\hat b$, which is $\F$-measurable. By \eqref{eq:A-Ealg} and $d_k\ge\mu_k/2$,
\begin{align*}
\E\big[\onedeg\,d_k(\cE v)_k^2\bigm|\F\big]&\le\frac2{\mu_k}\E\Big[\Big(\sum_jX_{kj}\big(v_j-(\hat Pv)_k\big)\Big)^2\Bigm|\F\Big]\\
&\le\frac2{\mu_k}\sum_jp_{kj}\big(v_j-(\hat Pv)_k\big)^2\le2(\hat Pv^2)_k ,
\end{align*}
since $p_{kj}=\mu_k\hat P_{kj}$ and the variance of $v_j$ under $\hat P_{k\cdot}$ is at most its second moment. Summing over $k$, and using $\mu_k\ge\dmin$ and the stationarity $\sum_k\mu_k\hat P_{kj}=\mu_j$,
\[
\sum_k(\hat Pv^2)_k\le\frac1{\dmin}\sum_k\mu_k(\hat Pv^2)_k=\frac1{\dmin}\sum_j\mu_jv_j^2\le\frac1{\dmin}\sum\mu\hat b^2 ,
\]
and by Lemma~\ref{lem:A-dual}(ii),
\[
\E\big[\onedeg\,\|\cE v\|^2_{\ell^2(d)}\bigm|\F\big]\le\frac2{\dmin}\sum\mu\hat b^2\le\frac{2n}{(1-\bar\gamma)^2\dmin^2}.
\]

\emph{The bound on $c$.} By the definition of $c_m$ and $(a+b)^2\le2a^2+2b^2$, $c_m^2\le2(b_m-\hat b_m)^2+16\dmin^{-2}(b_m^2+\hat b_m^2)$. The first part is bounded above, and Lemma~\ref{lem:A-dual}(ii) bounds the second deterministically by $Cn/\dmin^3$.
\end{proof}

The bounds on $G_1$, $G_2$ and $G_3$ involve sums such as $\sum_md_m\tau_m|b_m|$, in
which a weight $\tau_m$, such as a degree ratio or a local moment of the input, is
averaged with the dual weights. These are not averages over typical nodes. The dual
weight of $m$ reflects how often walks from the sample reach $m$, and walks reach
high-degree nodes more often: this is the size bias of network sampling, as in the
friendship paradox. Detailed balance makes the bias exact: $\sum_md_m\tau_m|b_m|$ is
bounded by the expected value of $\tau$ at the endpoints of walks started from the
sample. Sampled nodes are typical, and a walk from a typical node ends at a
degree-biased node, so the result is controlled by the degree-biased mean
$\E[g(\xi)\tau(\xi)]$, with no loss in $n$ and no dependence on walk length.
\begin{lemma}[Degree-biased averaging]\label{lem:weight-transfer}
Let $\tau_m\ge0$, $m\in\mathcal V$, be a permutation-equivariant function of $(\xi_{\mathcal V},A^{(0)})$, for example $\tau_m=\tau(\xi_m)$, and let $\Omega\subseteq\Odeg^{(0)}$ be an event determined by $(\xi_{\mathcal V},A^{(0)})$ and invariant under permutations of $\mathcal V$ (for instance $\Odeg^{(0)}$, $\Ods^{(0)}$, or the event $\OU$ of Appendix~\ref{app:flU}). Write
\[
\Eg{\tau;\Omega}:=\E\big[\one_\Omega\hat g^{(0)}_2\tau_2\big],\qquad \Eg{\tau}:=\E\big[\hat g^{(0)}_2\tau_2\big],
\]
so that $\Eg{\tau;\Omega}\le\Eg{\tau}$, and $\Eg{\tau}=\E[g(\xi)\tau(\xi)]$ when $\tau_m=\tau(\xi_m)$. On $\Odeg^{(0)}$, writing $d=d^{(0)}$ and $P=P^{(0)}$,
\begin{align}
\sum_md_m\tau_m|b_m|&\le\sum_{t\ge0}\bar\gamma^t\sum_{k\in\tilde S}(P^t\tau)_k,\label{eq:A-wt1}\\
\sum_md_m\tau_mb_m^2&\le\frac1{1-\bar\gamma}\sum_{t\ge0}\bar\gamma^t\sum_{k\in\tilde S}\frac{(P^t\tau)_k}{d_k},\label{eq:A-wt2}
\end{align}
and \eqref{eq:A-wt1} holds for $\hat b$ with $(\mu^{(0)},\hat P^{(0)})$ in place of $(d,P)$. Moreover
\[
\E\Big[\one_\Omega\sum_{k\in\tilde S}(P^t\tau)_k\Big]\le\frac{3n}{c_0}\Eg{\tau;\Omega},
\]
and consequently
\begin{gather*}
\E\Big[\one_\Omega\sum d\tau|b|\Big]\le Cn\,\Eg{\tau;\Omega},\\
\E\Big[\one_\Omega\sum d\tau b^2\Big]\le C\frac n{\dmin}\,\Eg{\tau;\Omega},\\
\E\Big[\one_\Omega\sum\mu^{(0)}\tau|\hat b|\Big]\le Cn\,\Eg{\tau;\Omega}.
\end{gather*}
We refer to the first two consequences as (DB1) and (DB2). The indicator is kept on the right-hand side so that bounds on $\tau$ that are themselves available only on $\Omega$, such as \eqref{eq:A-noiseU} for $\tau=\gn^2$ or Lemma~\ref{lem:A-twohop} for $\tau=X^{z}$ with $z$ built from $V^{(0)}$, can be inserted; when $\Eg{\tau}$ is bounded directly, we simply use $\Eg{\tau;\Omega}\le\Eg{\tau}$.
\end{lemma}

\begin{proof}
On $\Odeg^{(0)}$, every node has positive degree, so $P$ is row-stochastic, and
$d_mP_{mk}=A^{(0)}_{mk}$ is symmetric, so $P$ satisfies detailed balance with respect to $d$.
Iterating, $d_m(P^t)_{mk}=d_k(P^t)_{km}$ for every $t\ge0$. We use four properties in
turn.

\emph{Positivity.} Since $P_{mk}\ge0$ and $|\gamma_k|\le\bar\gamma$, for every vector $v$,
$|(P\Gamma v)_m|\le\bar\gamma(P|v|)_m$. Iterating and summing the Neumann series
$b=\sum_t(P\Gamma)^t\src$ gives
\[
|b_m|\le\sum_{t\ge0}\bar\gamma^t\,(P^t|\src|)_m .
\]

\emph{Detailed balance.} Multiplying by $d_m\tau_m$, summing over $m$ and exchanging
the roles of $m$ and $k$,
\[
\sum_md_m\tau_m(P^t|\src|)_m=\sum_{m,k}d_k(P^t)_{km}\tau_m|\src_k|=\sum_kd_k|\src_k|\,(P^t\tau)_k .
\]
Since $d_k|\src_k|=\one\{k\in\tilde S\}$, this gives \eqref{eq:A-wt1}: the weight $\tau$,
attached to the nodes $m$ that the walk reaches, has been moved to the sampled nodes,
where it is averaged by walks started from the sample. The same argument with $(\mu^{(0)},\hat P^{(0)})$ gives the version for $\hat b$.

\emph{Row-stochasticity.} For \eqref{eq:A-wt2}, Cauchy--Schwarz over $t$ with weights
$\bar\gamma^t$ gives $b_m^2\le(1-\bar\gamma)^{-1}\sum_t\bar\gamma^t(P^t|\src|)_m^2$, and
Jensen, valid because each row of $P^t$ is a probability vector, gives
$(P^t|\src|)_m^2\le(P^t\src^2)_m$. Detailed balance as above, with
$d_k\src_k^2=\one\{k\in\tilde S\}/d_k$, gives \eqref{eq:A-wt2}.

\emph{Stationarity.} For the expectation, $\tilde S$ is deterministic and the law of $(\xi_{\mathcal V},A^{(0)},\tau,\one_\Omega)$ is exchangeable over $\mathcal V$, because $\tau$ is permutation-equivariant and $\Omega$ is permutation-invariant. Hence
\[
\E\Big[\one_\Omega\sum_{k\in\tilde S}(P^t\tau)_k\Big]=|\tilde S|\,\E\Big[\one_\Omega\avg_{k\in\mathcal V}(P^t\tau)_k\Big].
\]
Let $\pi^d_k=d_k/\sum_ld_l$. For $u\ge0$, $\avg_ku_k\le(\bar d/d_{\min})\sum_k\pi^d_ku_k$, where $\bar d$ and $d_{\min}$ are the average and the minimum degree; and stationarity, $\sum_k\pi^d_kP_{kj}=\pi^d_j$, gives $\sum_k\pi^d_k(P^t\tau)_k=\sum_k\pi^d_k\tau_k$. On $\Omega\subseteq\Odeg^{(0)}$, $d_k\le\frac32\mu^{(0)}_k=\frac32\rho n_0\hat g^{(0)}_k$ and $d_{\min}\ge\dmin/2$, so
\[
\frac{\bar d}{d_{\min}}\sum_k\pi^d_k\tau_k=\frac1{d_{\min}}\avg_kd_k\tau_k\le\frac3{c_0}\avg_k\hat g^{(0)}_k\tau_k .
\]
Taking expectations, $\E[\one_\Omega\avg_k\hat g^{(0)}_k\tau_k]=\Eg{\tau;\Omega}$ by exchangeability.

\emph{Consequences.} Summing over $t$ with weights $\bar\gamma^t$ gives (DB1). For (DB2), additionally $1/d_k\le2/\dmin$ on $\Odeg^{(0)}$. For the version with $\hat b$, the same stationarity argument applies with $\pi^\mu_k=\mu^{(0)}_k/\sum_l\mu^{(0)}_l$, which is stationary for $\hat P^{(0)}$, and with $\mu^{(0)}_k=\rho n_0\hat g^{(0)}_k$ and $\min_k\mu^{(0)}_k\ge\dmin$ in place of the degree comparisons.
\end{proof}

\begin{lemma}[Averaging over the resampled node]\label{lem:A-two}
Let $r_m:=g(\xi_m)/\hat g^{(0)}_m$. On $\Odeg^{(0)}$, $g(\xi_m)\le\frac{2}{\rho n_0}r_md^{(0)}_m$. For nonnegative $\Hc_1$-measurable $x_m$,
\begin{gather*}
\E[w(\xi_m,\xi_1)\mid\Hc_1]=g(\xi_m),\\
\E\Big[\Big(\sum_mx_mw(\xi_m,\xi_1)\Big)^2\Bigm|\Hc_1\Big]\le\Big(\sum_mx_mg_m\Big)\Big(\sum_mx_mg_mK_2(\xi_m)\Big).
\end{gather*}
Moreover $\Eg{r}=1$, $\Eg{r^2}\le\E g^2/c_0$, and $\Eg{r^2K_2^2}\le\E[g^2K_2^2]/c_0\le\mathsf M_N/c_0^3$.
\end{lemma}

\begin{proof}
\emph{The degree comparison.} On $\Odeg^{(0)}$, $d^{(0)}_m\ge\mu^{(0)}_m/2=\rho n_0\hat g^{(0)}_m/2$, so $g(\xi_m)=r_m\hat g^{(0)}_m\le\frac2{\rho n_0}r_md^{(0)}_m$.

\emph{The averages over $\xi_1$.} Given $\Hc_1$, $\xi_1$ is uniform and independent of $\Hc_1$, so $\E[w(\xi_m,\xi_1)\mid\Hc_1]=\int w(\xi_m,v)\,dv=g(\xi_m)$. For the second moment, write $w(\xi_m,v)=g_m\,p(v\mid\xi_m)$ and apply Cauchy--Schwarz with weights $x_mg_m$:
\[
\Big(\sum_mx_mg_mp(v\mid\xi_m)\Big)^2\le\Big(\sum_mx_mg_m\Big)\Big(\sum_mx_mg_mp(v\mid\xi_m)^2\Big).
\]
Integrating over $v$, and using $\int p(v\mid\xi_m)^2dv=K_2(\xi_m)$, gives the claim.

\emph{The moments.} By exchangeability and $\E[\hat g^{(0)}_2\mid\xi_2]=g(\xi_2)$: $\Eg{r}=\E[\hat g^{(0)}_2\,g_2/\hat g^{(0)}_2]=\E g=1$; $\Eg{r^2}=\E[g_2^2/\hat g^{(0)}_2]\le\E g^2/c_0$; and $\Eg{r^2K_2^2}\le\E[g^2K_2^2]/c_0=\E[\kappa_2^2/g^2]/c_0\le\mathsf M_N/c_0^3$.
\end{proof}

\section{Proof of Theorem~\ref{thm:fluct}}\label{app:fluctuation}

\paragraph*{Conventions for this appendix} The leave-one-node-out objects of
Appendix~\ref{app:tools} are used throughout, and the reduced graph on
$\mathcal V=\{2,\dots,N\}$ is the default: in this appendix, $d$, $P$, $\hat g$,
$\mu$, $\hat P$, $\cE$, $R$ and $M$ without a superscript denote the reduced
objects $d^{(0)},P^{(0)},\hat g^{(0)},\mu^{(0)},\hat P^{(0)},\cE^{(0)},R^{(0)},M^{(0)}$
of Appendix~\ref{app:tools}, and $Y^{(0)},y^{(0)},V^{(0)}$ are the reduced
responses. Quantities indexed by node $1$ refer to the full graph, as node $1$ is
absent from the reduced one: $d_1$ is its degree,
$\hat g_1=(N-1)^{-1}\sum_{k\in\mathcal V}w(\xi_1,\xi_k)$,
$\mu_1=\rho(N-1)\hat g_1$ and $\hat P_{1m}=p_{m1}/\mu_1$. The responses $Y$, $Y^W$
and $V=Y-Y^W$ are those of the full graph, as everywhere else in the paper; where
the full-graph walk itself is needed we write $P^{\rm full}$ and
$\hat P^{\rm full}$. Events keep their superscripts ($\Odeg^{(0)}$,
$\Ods^{\rm full}$, and so on). A reader interested only in bounded inputs can
stop at the end of Appendix~\ref{app:flB}; Appendix~\ref{app:flU} is not used
elsewhere.

\paragraph*{The decomposition} By \eqref{eq:ES} it suffices to bound $\E V_1^2$ and $\frac Nn\sup_{\tilde S}\E G^2$. By Corollary~\ref{cor:A-V1}, $\E V_1^2\le C\E[gZ^2]/\dmin+CN^{-13}\E Z^2$ in both settings. By Lemma~\ref{lem:A-loo}, adding and subtracting $\gamma_mp_{m1}\beta_m\hat y_m$ and $\gamma_mA_{m1}\beta_m\hat y_m$, $G=G_1+G_2+G_3$ with
\begin{gather*}
G_1=\sum_m\gamma_m(A_{m1}-p_{m1})\beta_m\hat y_m,\\
G_2=\sum_m\gamma_mp_{m1}(\beta_m-\hat\beta_m)\hat y_m,\\
G_3=\sum_m\gamma_mA_{m1}\beta_m(\tilde Y_m-\hat y_m).
\end{gather*}
$G_1$ is the fluctuation of node $1$'s edges around their probabilities. $G_2$ is the discrepancy between the dual weights of the graph and of the expected-weight graph. $G_3$ is the feedback of the perturbation through the rest of the network. In setting (U), $G_3$ is split further, into the graph noise at node $1$'s neighbors, the graph noise at node $1$ itself (which has three sources), and the echo of node $1$'s own influence. Table~\ref{tab:A-map} lists every term, where it is bounded, and what it contributes to $n\Var(\bar V_S)$; the proofs below are organized by its rows.

\begin{table}[t]
\centering\footnotesize
\setlength{\tabcolsep}{4pt}
\begin{tabular}{>{\raggedright\arraybackslash}p{0.13\linewidth}>{\raggedright\arraybackslash}p{0.30\linewidth}>{\raggedright\arraybackslash}p{0.25\linewidth}>{\raggedright\arraybackslash}p{0.24\linewidth}}
\toprule
Term & What it is & Setting (B): bound on its share of $n\Var(\bar V_S)$ & Setting (U): bound on its share of $n\Var(\bar V_S)$\\
\midrule
$\E V_1^2$ & graph noise at a single node & $B^2/\dmin$ (Cor.~\ref{cor:A-V1}) & $K'/\dmin$ (Cor.~\ref{cor:A-V1})\\
\addlinespace
$G_1$ & node $1$'s edges fluctuate around their probabilities & $B^2/\dmin$ (Prop.~\ref{prop:A-flB}) & $K'/\dmin$ (Prop.~\ref{prop:A-G1})\\
$G_2$ & dual weights of the graph vs.\ of the expected-weight graph & $B^2\mathsf M_N(\varepsilon^2+\dmin^{-2})$ (Prop.~\ref{prop:A-flB}) & $K'/\dmin$ (Prop.~\ref{prop:A-G2})\\
$G_3$ & feedback of node $1$'s edges through the network & $B^2\mathsf M_N\varepsilon^2$ (Prop.~\ref{prop:A-flB}; uniform bound, the source of $\log N$) & split as below\\
\addlinespace
\quad $G_{3\gn}$ & graph noise in the neighbors' averages & & $K'/\dmin$ (Prop.~\ref{prop:A-G3a})\\
\quad $G^\flat$ & node $1$ averages its neighbors' graph noise & & $K'/\dmin$ (Prop.~\ref{prop:A-G3a})\\
\quad $G^\sharp$ & node $1$'s own edges act on its neighbors' expected responses & & $K'/\dmin$ (Prop.~\ref{prop:A-G3a})\\
\quad $\gamma_1\nu_1\Lambda^A$,\newline\quad $G_{3f}$ & echo of node $1$'s influence, and feedback $\vartheta$ & & $K'(\dmin^{-2}+W_N^4/N^2)$ (Prop.~\ref{prop:A-echo})\\
\addlinespace
bad events & concentration fails & $N^{-12}$ & $N^{-2}$\\
\midrule
total & Theorem~\ref{thm:fluct} & $B^2(1+\mathsf M_N)\dfrac{\log N}{N\rho_N}$ & $K'\Big(\dfrac1{N\rho_N}+\dfrac{W_N^4}{N^2}\Big)$\\
\bottomrule
\end{tabular}
\caption{Map of the proof of Theorem~\ref{thm:fluct}. Each entry is an upper bound up to a constant $C=C(\bar\gamma,c_0,\delta,K)$, after the Efron--Stein scaling by $N/n$ in \eqref{eq:ES}; $\varepsilon^2\asymp\log N/\dmin$ and $\dmin\asymp N\rho_N$. In setting (U), $V_1=\gamma_1(V^\flat+V^\sharp+\nu_1)$, $G^\flat=\gamma_1V^\flat\Lambda^A$, $G^\sharp=\gamma_1V^\sharp\Lambda^A$ and $G_3=V_1\Lambda^A-G_{3\gn}-G_{3f}$, as set out before Proposition~\ref{prop:A-G3a}.}
\label{tab:A-map}
\end{table}

\subsection{Setting (B)}\label{app:flB}
We first show that the network perturbation is negligible when the inputs are
bounded. By \eqref{eq:ES}, it suffices to show that a single node's graph noise,
$\E V_1^2$, and the effect of one node's edges on the sample, $\frac Nn\E G^2$, both
vanish. The first is immediate from Corollary~\ref{cor:A-V1}. For the second, we bound
the three parts of $G$ separately:
\begin{itemize}
\item $G_1$, the fluctuation of node $1$'s edges around their probabilities, is a sum
of independent centered terms given everything else, so its second moment is a sum
of variances;
\item $G_2$, the discrepancy between the dual weights of the actual and the
expected-weight graph, is controlled by the discrepancy between the dual vectors;
\item $G_3$, the feedback through the rest of the network, is small because, on a
high-probability event, the network noise is uniformly of order $\varepsilon B$ at
every node.
\end{itemize}
In each case, averaging over the type of node $1$ turns its edge weights into
expected degrees of its neighbors, and degree-biased averaging then reduces every bound
to the mixed norm $\mathsf M_N$. The uniform bound used for $G_3$ is the source of the
logarithm in the setting (B) rate.

\begin{proposition}\label{prop:A-flB}
Under (M1)--(M4) and $|Z|\le B$, for every $S$ with $|S|=n\le N$,
\[
n\Var(\bar V_S)\le CB^2(1+\mathsf M_N)\Big(\varepsilon^2+\frac1{\dmin}\Big)\le CB^2(1+\mathsf M_N)\frac{\log N}{N\rho_N}.
\]
In particular $n\Var(\bar V_S)\to0$ under (B3), since $\mathsf M_N\ge1$ forces $\log N=o(N\rho_N)$.
\end{proposition}

\begin{proof}
\emph{Preliminaries and the bad event.} By Lemma~\ref{lem:A-ops}, $\|Y^W\|_\infty,\|Y\|_\infty\le B_y$, so $|\hat y_m|\le2B_y$, $|V_i|\le2B_y$, and $|G|\le\sum_{\tilde S}(|V_i|+|V^{(0)}_i|)\le4nB_y$ always. Let $\OB$ be the intersection of the three events of Lemma~\ref{lem:A-conc} for the full graph and for the reduced graph, with $u=Y^W$ in the sup-norm event in both cases. Then $\Pr(\OB^c\mid\F)\le6N^{-14}$, and $\OB^c$ contributes at most $\frac Nn\cdot16n^2B_y^2\cdot6N^{-14}\le CB^2N^{-12}$ to $\frac Nn\E G^2$. Write $\oneds:=\one_{\Ods^{(0)}}$, which is $\Hc_1$-measurable; since $\OB\subseteq\Ods^{(0)}$, we may replace $\one_{\OB}$ by $\oneds$ in every upper bound below. We use Lemmas~\ref{lem:A-dual}, \ref{lem:weight-transfer} (with $\Omega=\Ods^{(0)}$, dropping the indicator on the right via $\Eg{\tau;\Ods^{(0)}}\le\Eg{\tau}$) and~\ref{lem:A-two}, with $K_g:=\E[g^2(1+K_2^2)]\le(1+c_0^{-2})\mathsf M_N$ by \ref{F:Mbig}.

\emph{$G_1$.} Given $(\Hc_1,\xi_1)$, $G_1$ is a sum of independent centered terms, with $\Var(A_{m1})\le p_{m1}=\rho w(\xi_m,\xi_1)$, $|\beta_m|\le|b_m|$ and $|\hat y_m|\le2B_y$. So on $\Ods^{(0)}$
\[
\E[G_1^2\mid\Hc_1,\xi_1]\le4\rho B_y^2\sum_mw(\xi_m,\xi_1)b_m^2 .
\]
Averaging over $\xi_1$ replaces $w(\xi_m,\xi_1)$ by $g(\xi_m)$ (Lemma~\ref{lem:A-two}), and $g(\xi_m)\le\frac2{\rho n_0}r_md_m$, so
\[
\E[G_1^2\oneds]\le\frac{8B_y^2}{n_0}\E\Big[\oneds\sum_md_mr_mb_m^2\Big]\le\frac{8B_y^2}{n_0}\cdot\frac{Cn}{\dmin}\Eg{r}=\frac{CB_y^2n}{N\dmin},
\]
by (DB2) and $\Eg{r}=1$. After the scaling by $N/n$, this is $CB^2/\dmin$.

\emph{$G_2$.} By Lemma~\ref{lem:A-dual}(v), $|\beta_m-\hat\beta_m|\le c_m$, so $|G_2|\le2B_y\rho\sum_mc_mw(\xi_m,\xi_1)$ with $c$ $\Hc_1$-measurable. By Lemma~\ref{lem:A-two} and $g_m\le\frac2{\rho n_0}r_md_m$,
\[
\E[G_2^2\oneds\mid\Hc_1]\le4B_y^2\rho^2\Big(\sum_mc_mg_m\Big)\Big(\sum_mc_mg_mK_2(\xi_m)\Big)\oneds\le\frac{16B_y^2}{n_0^2}\,\mathsf S_r\,\mathsf S_{rK_2}\,\oneds,
\]
where $\mathsf S_\tau:=\sum_md_mc_m\tau_m$. To bound $\mathsf S_\tau$, split $c_m$ according to its definition and apply Cauchy--Schwarz to the first part:
\[
\mathsf S_\tau\le\Big(\sum d(b-\hat b)^2\Big)^{1/2}\Big(\sum d\tau^2\Big)^{1/2}+\frac2{\dmin}\sum d\tau(|b|+|\hat b|).
\]
For the first term, Lemma~\ref{lem:A-dual}(iv) gives $\sum d(b-\hat b)^2\le C\varepsilon^2n/\dmin$ on $\Ods^{(0)}$, and $\E[\oneds\sum_md_m\tau_m^2]\le\frac32\rho n_0\,\E\sum_m\hat g_m\tau_m^2\le C\rho N^2\Eg{\tau^2}$. For the second term, Cauchy--Schwarz with weights $d|b|$ and Lemma~\ref{lem:A-dual}(i) give $(\sum d\tau|b|)^2\le\frac n{1-\bar\gamma}\sum d\tau^2|b|$, which has expectation at most $Cn^2\Eg{\tau^2}$ by (DB1); likewise for $\hat b$. So
\[
\E[\oneds\mathsf S_\tau^2]\le C\Eg{\tau^2}\Big(\varepsilon^2\frac{n\rho N^2}{\dmin}+\frac{n^2}{\dmin^2}\Big).
\]
By Cauchy--Schwarz in $\E$, $\E[\mathsf S_r\mathsf S_{rK_2}\oneds]\le(\E[\oneds\mathsf S_r^2]\,\E[\oneds\mathsf S_{rK_2}^2])^{1/2}$, and $\Eg{r^2}\Eg{r^2K_2^2}\le CK_g^2$ by Lemma~\ref{lem:A-two}. Hence, using $N/n_0^2\le C/N$ and $\rho N\le C\dmin$,
\[
\frac Nn\E[G_2^2\oneds]\le\frac{CB_y^2K_g}{nN}\Big(\varepsilon^2\frac{n\rho N^2}{\dmin}+\frac{n^2}{\dmin^2}\Big)\le CB_y^2K_g\Big(\varepsilon^2+\frac1{\dmin^2}\Big).
\]

\emph{$G_3$.} By Lemma~\ref{lem:A-loo} and Lemma~\ref{lem:A-ops}, $\tilde Y_m-\hat y_m=V_1-(PV)_m-(\cE Y^W)_m$. On $\OB$, $\|V\|_\infty\le C\varepsilon B_y$, because $V=R_P\Gamma\cE^{\rm full}Y^W$ with $\|R_P\Gamma\|_{\infty\to\infty}\le\bar\gamma/(1-\bar\gamma)$ everywhere, and $|\cE^{\rm full}Y^W|\le2\varepsilon\|Y^W\|_\infty$ entrywise on the full-graph degree and sup-norm events; likewise $\|\cE Y^W\|_\infty\le2\varepsilon B_y$ on the reduced-graph events. So $|\tilde Y_m-\hat y_m|\le2\|V\|_\infty+\|\cE Y^W\|_\infty\le C\varepsilon B_y$, and
\begin{gather*}
|G_3|\le C\varepsilon B_y\Labs,\\
\Labs:=\sum_mA_{m1}|b_m| .
\end{gather*}
Given $(\Hc_1,\xi_1)$, $\E[\Labs^2\mid\Hc_1,\xi_1]=(\sum_mp_{m1}|b_m|)^2+\sum_m\Var(A_{m1})b_m^2\le\rho^2(\sum_mw_{m1}|b_m|)^2+\rho\sum_mw_{m1}b_m^2$. Averaging over $\xi_1$ with Lemma~\ref{lem:A-two}, and converting $g_m$ to $d_mr_m$ as for $G_2$,
\[
\E[\Labs^2\oneds\mid\Hc_1]\le\frac{4}{n_0^2}\Big(\sum dr|b|\Big)\Big(\sum drK_2|b|\Big)+\frac2{n_0}\sum dr\,b^2 .
\]
By Cauchy--Schwarz, (DB1) and (DB2), the expectation is at most $C(K_gn^2/N^2+n/(N\dmin))$. Hence $\frac Nn\E[G_3^2\one_{\OB}]\le C\varepsilon^2B_y^2(K_gn/N+1/\dmin)\le C\varepsilon^2B_y^2K_g$.

\emph{Collecting.} Combining the three bounds with \eqref{eq:ES}, $\E V_1^2\le CB^2/\dmin$ and the bad-event bound gives the claim, since $\varepsilon^2\asymp\log N/\dmin$ and $\dmin\asymp N\rho_N$.
\end{proof}

\subsection{Setting (U)}\label{app:flU}
\begin{table}[htbp]
\centering\small
\begin{tabular}{p{0.22\linewidth}p{0.55\linewidth}p{0.15\linewidth}}
\toprule
Symbol & Meaning & Defined in\\
\midrule
$K'$ & finite combination of the moments in (U2) and (U4) & App.~\ref{app:fluctuation}\\
$x_m$, $\kappa_m$ & $|\beta_m|$; $\kappa_2(\xi_m)$ & App.~\ref{app:fluctuation}\\
$\Ohg$, $\OU$, $\oneU$, $\Odegone$ & all $\hat g_m\ge g_m/2$; the good event $\Ods^{(0)}\cap\Ohg$ and its indicator; node $1$'s degree within a factor $2$ of its mean & App.~\ref{app:fluctuation}\\
$\mathrm{Bad}$ & $\OU^c\cup\Odegone^c$ & App.~\ref{app:fluctuation}\\
$\mathsf d$, $\vartheta$ & change in the expected-weight responses when node $1$ is added; its largest entry & Lem.~\ref{lem:A-fb}\\
$\bar A$ & average $|y^{(0)}|$ over $\mathcal V$ & Lem.~\ref{lem:A-fb}\\
$u$, $\gn$ & $\hat Py^{(0)}$; graph noise in a node's neighborhood average & \eqref{eq:A-noiseU}\\
$\hat W$ & degree-weighted squared edge profile $\sum_mw(\xi_m,\cdot)^2/\hat g_m$ & Lem.~\ref{lem:A-dw}\\
$\nb{z}(v)$ & average of $z$ over the potential neighbors of a node of type $v$; $\nb{(y^{(0)})^2}$ is the case $z=(y^{(0)})^2$ & Lem.~\ref{lem:A-fb}\\
$X^z_m$ & two-hop average of $z$ from $m$ through node $1$ & Lem.~\ref{lem:A-twohop}\\
$\tau^F_m$ & degree-weighted node-$1$ average of $F$ & Lem.~\ref{lem:A-deg1}\\
$\echo$, $f_m$ & change in the actual responses when node $1$ is added (the echo, the actual-graph counterpart of $\mathsf d$); echo and feedback seen at $m$ & App.~\ref{app:fluctuation}\\
$V^\flat$, $V^\sharp$, $\nu_1$ & noise at node $1$: from its neighbors' noise, from its own edges, and from the echo & App.~\ref{app:fluctuation}\\
$\Lambda^A$, $\Labs$ & $\sum_m\gamma_mA_{m1}\beta_m$; its bound $\sum_mA_{m1}x_m$ & App.~\ref{app:fluctuation}\\
$s_{ij}$, $t_{ij}$ & $\sum_mp_{m1}x_m^iz_m^j$; $\sum_mp_{m1}x_m^i|a_m|^j$ & Lem.~\ref{lem:A-bern}\\
$\mathcal S$, $\Sigma$ & size of node $1$'s inputs; size of the echo & Prop.~\ref{prop:A-echo}\\
\bottomrule
\end{tabular}
\caption{Notation for setting (U) of the network fluctuation.}
\label{tab:not-U}
\end{table}

Throughout this subsection, assume (M1)--(M4) and (U1)--(U4). We write $K'$ for a finite combination of $\E[gZ^4]$, $\E g^4$, $\mathsf M_N$, $\mathfrak D$, $\mathcal H$ and $\mathfrak G$, which is bounded by (U2) and (U4), and we take $N$ large enough that $W_N^2\le N$.

\paragraph*{Notation and events} Let $x_m:=|\beta_m|\le|b_m|$ and $\kappa_m:=\kappa_2(\xi_m)$. Define
\begin{gather*}
\Ohg:=\{\hat g_m\ge g_m/2\ \forall m\in\mathcal V\},\\
\OU:=\Ods^{(0)}\cap\Ohg,\\
\oneU:=\one_{\OU},\\
\Odegone:=\{\mu_1/2\le d_1\le2\mu_1\}.
\end{gather*}
By Lemmas~\ref{lem:A-lt}, \ref{lem:A-conc} and \ref{F:M3}, $\Pr(\OU^c)\le N^{-10}$ (once $C_*\ge44$) and $\Pr(\Odegone^c\mid\Hc_1,\xi_1)\le2e^{-\dmin/12}$. The event $\OU$ is determined by $(\xi_{\mathcal V},A^{(0)})$ and is invariant under permutations of $\mathcal V$, so Lemma~\ref{lem:weight-transfer} applies with $\Omega=\OU$; throughout this subsection (DB1) and (DB2) refer to that choice. On $\OU$, Lemma~\ref{lem:A-dual} gives $\sum d|b|\le Cn$, $\sum db^2\le Cn/\dmin$ and $\|b\|_\infty\le C/\dmin$. It also gives $\dmin/2\le d_m\le\frac32\mu_m$ and $g_m\le2\hat g_m\le4d_m/(\rho n_0)$.

\paragraph*{Outline} The structure is the same as in setting (B): we bound
$\E V_1^2$ and the three parts of $G$. What changes is that, with heavy-tailed inputs,
the network noise is no longer uniformly small, so every factor that was bounded by a
constant in setting (B) must instead be controlled through its moments. The argument
proceeds in four steps.
\begin{enumerate}
\item \emph{Moments of the graph noise.} Lemma~\ref{lem:A-noise} and
\eqref{eq:A-noiseU} bound the second and fourth moments of the noise at a typical
node, with no logarithm. Lemma~\ref{lem:A-fb} bounds the feedback $\vartheta$: how
much adding node $1$ changes the expected-weight responses of the others. This is the
first place where (U1) enters.
\item \emph{Averaging over node $1$.} Lemmas~\ref{lem:A-dw}--\ref{lem:A-deg1} compute
what remains after averaging over the type of node $1$. Each lemma handles one way in
which node $1$'s neighbors can share structure, such as a common neighbor, a repeated
edge, or a high-degree neighbor, and yields one of the functionals in (U4).
\item \emph{The terms $G_1$ and $G_2$.} These follow the setting (B) proof, with
moments in place of uniform bounds. For $G_2$, the discrepancy between the dual vectors is bounded
conditionally on the latent types (Lemma~\ref{lem:A-energy}), which avoids the
logarithm.
\item \emph{The term $G_3$.} The feedback is split into the noise at node $1$'s
neighbors, the noise at node $1$ itself, and the echo of node $1$'s influence
returning through its neighbors. The resulting products of sums over node $1$'s edges
are expanded exactly with Lemma~\ref{lem:A-bern}, and each coincidence pattern is
bounded as recorded in Table~\ref{tab:A-patterns}.
\end{enumerate}
Each part contributes at most $CK'/\dmin$, up to terms carrying the factor
$\|w_N\|_\infty^4/N^2$, which gives Proposition~\ref{prop:A-flU}.

\begin{lemma}[Expected-weight identities]\label{lem:A-id}
\begin{itemize}
\item[(a)] $\sum_m\hat g_mh_m(\hat Pv)_m=\sum_kv_k\hat g_k(\hat Ph)_k$.
\item[(b)] $\sum_k\hat g_k|(\hat Pv)_k|^s\le\sum\hat g|v|^s$ for $s\ge1$, and $\sum\hat g|y^{(0)}|^s\le(1-\bar\gamma)^{-s}\sum\hat g|Z|^s$.
\item[(c)] On $\Odeg^{(0)}$, $\sum_k\hat g_k(Pv)_k\le3\sum_k\hat g_kv_k$ for $v\ge0$.
\end{itemize}
\end{lemma}

\begin{proof}
(a) Both sides equal $\sum_{m,k}\hat g_m\hat P_{mk}h_mv_k$, because $\hat g_m\hat P_{mk}=w_{mk}/n_0$ is symmetric. (b) The first inequality is Jensen row by row followed by stationarity of $\hat g$; the second is \ref{F:ops}(ii)--(iii) with $\nu=\hat g$. (c) On $\Odeg^{(0)}$, $\hat g_k\le2d_k/(\rho n_0)$ and $d_k\le\frac32\rho n_0\hat g_k$. So, using stationarity $\sum_kd_k(Pv)_k=\sum_kd_kv_k$,
\[
\sum_k\hat g_k(Pv)_k\le\frac2{\rho n_0}\sum_kd_k(Pv)_k=\frac2{\rho n_0}\sum_kd_kv_k\le3\sum_k\hat g_kv_k .
\]
\end{proof}

By Lemma~\ref{lem:A-noise} with $\mathcal V=\{2,\dots,N\}$ and exchangeability, for $s\in\{2,4\}$ and $\gn:=PV^{(0)}+\cE y^{(0)}$,
\begin{equation}\label{eq:A-noiseU}
\E\big[\oneU\hat g_2\big(|V^{(0)}_2|^s+|\gn_2|^s\big)\big]\le C\dmin^{-s/2}\E[g|Z|^s],\qquad\E\big[\oneU\hat g_2(Y^{(0)}_2)^4\big]+\E\big[\hat g_2(y^{(0)}_2)^4\big]\le C\E[gZ^4].
\end{equation}
The indicator is essential in the first display, since $V^{(0)}$ and $\gn$ are controlled only on the degree event; the bound on $y^{(0)}$ holds without it.

The leave-one-node-out argument compares quantities with and without node $1$. On
the expected-weight graph, adding node $1$ changes every other node's response
slightly, because node $1$ becomes a weighted neighbor of everyone. For a node $m$,
node $1$ enters the neighborhood average with weight $p_{m1}/\mu_m$, which is of order
$\|w_N\|_\infty/N$. It contributes the difference between its own response and $m$'s
current average, and that change then propagates through the network. The
\emph{feedback} $\vartheta$ is the largest resulting change in any node's
expected-weight response. It is small because each node's share in another's average
is small. The bound depends on $\|w_N\|_\infty$, since a hub of large weight has a
larger share, and on the size of node $1$'s input and of the typical response, since
those set the size of the difference it contributes. In setting (B) all of these are
bounded and the feedback is of order $1/(N\rho_N)$; with heavy-tailed inputs, the
next lemma bounds it explicitly.

\begin{lemma}[Feedback]\label{lem:A-fb}
Let $\mathsf d:=Y^W_{\mathcal V}-y^{(0)}$ and $\vartheta:=\|\mathsf d\|_\infty$. If $CW_N\le N$, then
\begin{gather*}
\vartheta\le\frac{CW_N}{N}\big(|Z_1|+W_N\bar A\big),\\
\bar A:=\avg_{k\in\mathcal V}|y^{(0)}_k|,\\
\E\bar A^4\le C\E[gZ^4].
\end{gather*}
Moreover $(Y^W_1)^2\le3\big(Z_1^2+\nb{(y^{(0)})^2}(\xi_1)+\vartheta^2\big)$, where $\nb{(y^{(0)})^2}(v):=((N-1)c_0)^{-1}\sum_{k\in\mathcal V}w(v,\xi_k)(y^{(0)}_k)^2$ is the average of $(y^{(0)})^2$ over the potential neighbors of a node of type $v$.
\end{lemma}

\begin{proof}
\emph{The equation for $\mathsf d$.} For $m\in\mathcal V$, the expected-weight version of the computation in Lemma~\ref{lem:A-loo} gives $(\hat P^{\rm full}Y^W)_m=(\hat PY^W)_m+\varpi_m\hat y_m$, with $\varpi_m=p_{m1}/(\mu_m+p_{m1})\le W_N/(c_0n_0)$. Subtracting the equations $Y^W_m=Z_m+\gamma_m(\hat P^{\rm full}Y^W)_m$ and $y^{(0)}_m=Z_m+\gamma_m(\hat Py^{(0)})_m$ gives $(I-\Gamma\hat P)\mathsf d=\Gamma(\varpi\hat y)$, so $\mathsf d=M\Gamma(\varpi\hat y)$ and, by \ref{F:ops}, $\vartheta\le\frac{\bar\gamma}{1-\bar\gamma}\max_m\varpi_m|\hat y_m|$.

\emph{Bounding $\hat y_m$.} By definition, $\hat y_m=Y^W_1-(\hat PY^W)_m$, and $Y^W_{\mathcal V}=y^{(0)}+\mathsf d$, so $|\hat y_m|\le|Y^W_1|+|(\hat Py^{(0)})_m|+\vartheta$. Similarly $Y^W_1=Z_1+\gamma_1(\hat P_1Y^W_{\mathcal V})$ gives $|Y^W_1|\le|Z_1|+|\hat P_1y^{(0)}|+\vartheta$. Since every row of $\hat P$ and $\hat P_1$ is bounded entrywise by $W_N/(c_0n_0)$, both averages are at most $(2W_N/c_0)\bar A$. Combining, $\vartheta\le\frac{CW_N}N\big(|Z_1|+W_N\bar A+\vartheta\big)$, and since $CW_N/N\le1/2$, the $\vartheta$ on the right can be absorbed.

\emph{Moments and $Y^W_1$.} By Jensen, $\bar A^4\le\avg_k(y^{(0)}_k)^4\le c_0^{-1}\avg_k\hat g_k(y^{(0)}_k)^4$, and Lemma~\ref{lem:A-id}(b) gives $\E\bar A^4\le C\E[gZ^4]$. Finally, $(Y^W_1)^2\le3(Z_1^2+(\hat P_1y^{(0)})^2+\vartheta^2)$ and $(\hat P_1y^{(0)})^2\le\hat P_1(y^{(0)})^2\le\nb{(y^{(0)})^2}(\xi_1)$, since $\hat P_{1k}=w(\xi_1,\xi_k)/((N-1)\hat g_1)$ and $\hat g_1\ge c_0$.
\end{proof}

\begin{lemma}[Degree-weighted averaging over $\xi_1$]\label{lem:A-dw}
Let $\hat W(v):=\sum_{m\in\mathcal V}w(\xi_m,v)^2/\hat g_m$. For $\Hc_1$-measurable $x\ge0$ and $F\ge0$, and any $\Hc_1$-measurable $u$, on $\Odeg^{(0)}$:
\begin{gather*}
\E\Big[F(\xi_1)\Big(\sum_mx_mw(\xi_m,\xi_1)\Big)^2\Bigm|\Hc_1\Big]\le\frac2{\rho n_0}\Big(\sum dx^2\Big)\int F\hat W,\\
\E\Big[\Big(\sum_mx_m|u_m|w(\xi_m,\xi_1)\Big)^2\Bigm|\Hc_1\Big]\le\frac2{\rho n_0}\Big(\sum dx^2\Big)\sum_m\frac{u_m^2\kappa_m}{\hat g_m}.
\end{gather*}
\end{lemma}

\begin{proof}
For any numbers $\omega_m$, Cauchy--Schwarz with weights $d_m$ gives
\[
\Big(\sum_mx_m\omega_m\Big)^2\le\Big(\sum_md_mx_m^2\Big)\Big(\sum_m\frac{\omega_m^2}{d_m}\Big),
\]
and on $\Odeg^{(0)}$, $1/d_m\le2/(\rho n_0\hat g_m)$. For the first display, take $\omega_m=w(\xi_m,v)$, multiply by $F(v)$ and integrate over $v$; this gives $\int F\hat W$. For the second, take $\omega_m=|u_m|w(\xi_m,v)$ and integrate over $v$, using $\int w(\xi_m,v)^2dv=\kappa_m$.
\end{proof}

\begin{lemma}[$\hat W$-moments]\label{lem:A-Wmom}
For $\Hc_1$-measurable $z\ge0$, let $\nb{z}(v):=((N-1)c_0)^{-1}\sum_{k\in\mathcal V}w(v,\xi_k)z_k$, and let $u:=\hat Py^{(0)}$. Let $\Omega$ be any event (we use $\OU$ and the sure event). Then:
\begin{itemize}
\item[(a)] $\E\int Z^2\hat W\le CN(\mathsf M_N\E[gZ^4])^{1/2}$;
\item[(b)] $\E[\one_\Omega\int\nb{z}\hat W]\le CN\,\E[\one_\Omega\avg_k\hat g_kz_k^2]^{1/2}\big(\mathcal H+W_N^2\mathsf M_N/N\big)^{1/2}$;
\item[(c)] $\E[\sum_mu_m^2\kappa_m/\hat g_m]\le CN(\E[gZ^4]\mathsf M_N)^{1/2}$; and the two quantities
\[
\E\Big[\oneU\sum_m\gn_m^2\kappa_m/\hat g_m\Big],\qquad\E\Big[\oneU\sum_m(V^{(0)}_m)^2\kappa_m/\hat g_m\Big]
\]
obey the same bound with an extra factor $\dmin^{-1}$.
\end{itemize}
\end{lemma}

\begin{proof}
(a) Using $\hat g\ge c_0$, $\E\hat W(v)\le c_0^{-1}N\kappa_2(v)$, so by Cauchy--Schwarz, $\E\int Z^2\hat W\le(N/c_0)\int Z^2\kappa_2\le(N/c_0)\|Z\|^2_{L^4}\mathsf M_N^{1/2}$, and $\|Z\|_{L^4}^4\le\E[gZ^4]/c_0$.

(b) \emph{Rewriting.} Exchanging sum and integral, $\int\nb{z}\hat W=((N-1)c_0)^{-1}\sum_kz_kH_k$ with $H_k:=\int w(v,\xi_k)\hat W(v)\,dv$. By Cauchy--Schwarz with weights $\hat g_k$,
\[
\int\nb{z}\hat W\le\frac1{(N-1)c_0}\Big(\sum_k\hat g_kz_k^2\Big)^{1/2}\Big(\sum_kH_k^2/\hat g_k\Big)^{1/2}.
\]
\emph{The second moment of $H$.} Since $\hat g\ge c_0$, $H_2\le c_0^{-1}\sum_{m\in\mathcal V}\varkappa(\xi_m,\xi_2)$, where $\varkappa(x,z):=\int w(x,v)^2w(v,z)\,dv\le W_N\kappa_2(x)$ and $\int\varkappa(x,z)\,dx=\mathfrak h(z)$. Squaring the sum over $m$: the off-diagonal pairs $m\ne m'$, both different from $2$, average over $\xi_m,\xi_{m'}$ to at most $N^2\mathfrak h(\xi_2)^2$, which integrates to $N^2\mathcal H$; the diagonal pairs and the term $m=2$ contribute at most $NW_N^2\mathsf M_N$, using $\varkappa\le W_N\kappa_2$. Hence $\E H_2^2\le Cc_0^{-2}\big(N^2\mathcal H+NW_N^2\mathsf M_N\big)$. Multiplying the display by $\one_\Omega$, taking expectations, and using Cauchy--Schwarz in $\E$ (the indicator stays with the first factor) together with exchangeability gives (b).

(c) By Jensen, $u^2\le\hat P(y^{(0)})^2$. Lemma~\ref{lem:A-id}(a) with $h=\kappa/\hat g^2$ moves $\hat P$ onto $\kappa/\hat g^2$:
\[
\sum_m\frac{u_m^2\kappa_m}{\hat g_m}\le\sum_m\hat g_m\frac{\kappa_m}{\hat g_m^2}(\hat P(y^{(0)})^2)_m=\sum_k(y^{(0)}_k)^2\hat g_k\Big(\hat P\frac{\kappa}{\hat g^2}\Big)_k ,
\]
and Cauchy--Schwarz with weights $\hat g$, followed by Lemma~\ref{lem:A-id}(b) with $s=2$, bounds this by $(\sum\hat g(y^{(0)})^4)^{1/2}(\sum\kappa^2/\hat g^3)^{1/2}$. Taking expectations, the second display of \eqref{eq:A-noiseU} (which needs no indicator) and $\E\kappa_2^2=\mathsf M_N$ give the claim. For $\gn$ and $V^{(0)}$, apply Cauchy--Schwarz directly to $\sum_m\gn_m^2\kappa_m/\hat g_m$, and use the fourth-moment bound in \eqref{eq:A-noiseU}, which carries the factor $\dmin^{-2}$ and the indicator $\oneU$; its square root gives $\dmin^{-1}$.
\end{proof}

Many of the terms in $G_3$ multiply the edge from node $1$ to a neighbor $m$, which
carries the dual weight at $m$, by node $1$'s average over its neighbors $k$ of a
quantity $z_k$, such as the graph noise or the squared responses. Averaging over the
type of node $1$ links $m$ and $k$ through their common neighbor: the result is the
two-step weight $W_2(\xi_m,\xi_k)=\int w(\xi_m,v)w(v,\xi_k)\,dv$, and the averaged
quantity is $X^z_m$ below. This pattern, a path $m\to1\to k$ through the resampled
node, is the most frequent one in Table~\ref{tab:A-patterns}. The next lemma bounds it
by the size of $z$ and by $\mathfrak D$, a degree-weighted second moment of the
two-step kernel, which measures how much nodes interact through common neighbors.

\begin{lemma}[Two-hop weights]\label{lem:A-twohop}
For $\Hc_1$-measurable $z\ge0$, let $X^z_m:=((N-1)c_0)^{-1}\sum_kz_kW_2(\xi_m,\xi_k)=\E[w(\xi_m,\xi_1)\nb{z}(\xi_1)\mid\Hc_1]$. Then, for every event $\Omega$ determined by $(\xi_{\mathcal V},A^{(0)})$ and invariant under permutations of $\mathcal V$ (including the sure event and $\OU$),
\[
\E[\one_\Omega\hat g_2X^z_2]\le C\,\E\big[\one_\Omega\avg_k\hat g_kz_k^2\big]^{1/2}\big(\mathfrak D+W_N^2/N\big)^{1/2}.
\]
\end{lemma}

\begin{proof}
\emph{Rewriting.} By exchangeability (here the permutation invariance of $\Omega$ is used), $\E[\one_\Omega\hat g_2X^z_2]=\E[\one_\Omega\avg_m\hat g_mX^z_m]=\E\big[\one_\Omega((N-1)c_0)^{-1}\sum_kz_k\hat H_k\big]$ with $\hat H_k:=\avg_m\hat g_mW_2(\xi_m,\xi_k)$. By Cauchy--Schwarz, first in $k$ with weights $\hat g_k$ and then in $\E$ (keeping the indicator with the factor $\sum_k\hat g_kz_k^2$), it suffices to show $\E[\hat H_2^2]\le C(\mathfrak D+W_N^2/N)$.

\emph{The second moment.} By Jensen, $\hat H_2^2\le\avg_m\hat g_m^2W_2(\xi_m,\xi_2)^2$. For $m\ne2$, write $\hat g_m\le\hat g'_m+W_N/n_0$, so that $\hat g_m^2\le2\hat g'^2_m+2W_N^2/n_0^2$, where $\hat g'_m$ omits node $2$ and is therefore independent of $\xi_2$ given $\xi_m$. Then $\E[\hat g'^2_m\mid\xi_m]\le g_m^2+\kappa_m/(n_0-1)$, and the leading term averages to $2\iint g(x)^2W_2(x,z)^2\,dx\,dz=2\mathfrak D$. Using $W_2(x,z)\le W_Ng(x)$, the two corrections are at most $C\frac{W_N^2}{N}(\E\kappa_2^2)^{1/2}(\E g^4)^{1/2}$ and $C\frac{W_N^4}{N^2}\E g^2\le C\frac{W_N^2}{N}\E g^2$ (by $W_N^2\le N$). The term $m=2$ has weight $1/(N-1)$ and is at most $CW_N^2\E g^4/N$.
\end{proof}

\begin{lemma}[Degree-weighted node-$1$ average]\label{lem:A-deg1}
For $F\ge0$, let $\tau^F_m:=\int\hat g_{(v)}F(v)w(\xi_m,v)\,dv$ with $\hat g_{(v)}=(N-1)^{-1}\sum_{k\in\mathcal V}w(v,\xi_k)$, so that $\E[\hat g_1F(\xi_1)w(\xi_m,\xi_1)\mid\Hc_1]=\tau^F_m$. Then $\E[\hat g_2\tau^F_2]\le\big(1+\frac{2W_N}{(N-1)c_0}\big)\int g^2(Tg)F$. For $F=Z^2$, this is at most $C\mathfrak G^{1/2}\E[gZ^4]^{1/2}$.
\end{lemma}

\begin{proof}
\emph{Expansion.} Write $\hat g_2=n_0^{-1}\sum_{j\in\mathcal V\setminus2}w(\xi_2,\xi_j)$
and $\hat g_{(v)}=(N-1)^{-1}\sum_{k\in\mathcal V}w(v,\xi_k)$, so that
\[
\hat g_2\tau^F_2=\frac1{n_0(N-1)}\sum_{j\ne2}\sum_{k}\int w(\xi_2,\xi_j)\,w(v,\xi_k)\,F(v)\,w(\xi_2,v)\,dv .
\]
We split the double sum into three cases.

\emph{Distinct indices, $k\ne j$ and $k\ne2$.} The types $\xi_2,\xi_j,\xi_k$ are
independent, and averaging over $\xi_j$ and $\xi_k$ replaces $w(\xi_2,\xi_j)$ by
$g(\xi_2)$ and $w(v,\xi_k)$ by $g(v)$. Averaging over $\xi_2$ then gives
\[
\int\Big(\int g(x)w(x,v)\,dx\Big)g(v)F(v)\,dv=\int g(v)^2(Tg)(v)F(v)\,dv,
\]
since $\int g(x)w(x,v)\,dx=g(v)(Tg)(v)$.

\emph{A shared neighbor, $k=j$.} Node $2$ and the point $v$ are joined through the
same node $j$, which gives the factor $w(\xi_2,\xi_j)w(v,\xi_j)$. Bounding
$w(v,\xi_j)\le W_N$ and averaging over $\xi_j$ gives at most $W_Ng(\xi_2)$, and
averaging over $\xi_2$ leaves $W_N\int g(Tg)F$.

\emph{Node $2$ itself, $k=2$.} Node $2$ appears in $\hat g_{(v)}$, which gives the
factor $w(v,\xi_2)w(\xi_2,v)=w(\xi_2,v)^2$. Bounding one factor by $W_N$ and averaging
$\hat g_2$ over the other nodes leaves $W_N\int g(Tg)F$.

Each of the last two cases has relative weight at most $1/(N-1)$, and
$\int g(Tg)F\le c_0^{-1}\int g^2(Tg)F$, which gives the factor $1+2W_N/((N-1)c_0)$.

\emph{The case $F=Z^2$.} By Cauchy--Schwarz,
$\int g^2(Tg)Z^2=\int g^{3/2}(Tg)\cdot g^{1/2}Z^2\le\mathfrak G^{1/2}\E[gZ^4]^{1/2}$.
\end{proof}
\paragraph*{Bad events and the first Efron--Stein term} The bad event
$\mathrm{Bad}:=\OU^c\cup\Odegone^c$ collects the ways in which the comparisons used
in the proof can fail:
\begin{itemize}
\item $\OU^c$ means that, in the graph without node $1$, some degree is far from
its expected value, or the adjacency noise fails the spectral bound, or some node's
empirical expected degree $\hat g_m$ falls below half its population value $g_m$; and
\item $\Odegone^c$ means that node $1$'s own degree is more than a factor of two away
from its expected value.
\end{itemize}
Each of these is a large deviation of a sum of many independent terms, so
$\Pr(\mathrm{Bad})\le N^{-10}$ under (M4) once $C_*\ge132$ (the binding constraint is $2e^{-\dmin/12}\le2N^{-C_*/12}$). On the bad event we use no structure at
all, only the crude bound $|G|\le\sum_{\tilde S}(|Y|+|Y^{(0)}|+|Y^W|+|y^{(0)}|)\le Cn\max_i|Z_i|$. By Cauchy--Schwarz, and $\E\max_iZ_i^4\le N\E Z^4$,
\[
\E[G^2\one_{\rm Bad}]\le Cn^2(\E\max_iZ_i^4)^{1/2}\Pr(\mathrm{Bad})^{1/2}\le Cn^2(N\E Z^4)^{1/2}N^{-5}\le Cn^2N^{-4};
\]
the polynomially small probability makes this contribution negligible. On the good event, $G^2\le3(G_1^2+G_2^2+G_3^2)$.

\begin{proposition}[$G_1$]\label{prop:A-G1}
$\frac Nn\E[G_1^2\oneU]\le CK'/\dmin$.
\end{proposition}

\begin{proof}
\emph{Conditional variance.} Given $(\Hc_1,\xi_1)$, $G_1$ is a sum of independent centered terms, so $\E[G_1^2\mid\Hc_1,\xi_1]\le\rho\sum_mw_{m1}x_m^2\hat y_m^2$.

\emph{Bounding $\hat y_m$.} With $u_m=(\hat Py^{(0)})_m$, the proof of Lemma~\ref{lem:A-fb} gives $|\hat y_m|\le|Y^W_1|+|u_m|+\vartheta$, and Lemma~\ref{lem:A-fb} bounds $Y^W_1$ and $\vartheta$:
\[
\hat y_m^2\le9\big(Z_1^2+\nb{(y^{(0)})^2}(\xi_1)+u_m^2+2\vartheta^2\big),\qquad\vartheta^2\le\frac{CW_N^2}{N^2}\big(Z_1^2+W_N^2\bar A^2\big).
\]

\emph{Averaging over $\xi_1$.} Each term is an $\Hc_1$-measurable quantity times a function of $\xi_1$, multiplied by $w_{m1}$. Averaging over $\xi_1$: $\E[w_{m1}Z_1^2\mid\Hc_1]=g_m(TZ^2)(\xi_m)$; $\E[w_{m1}\nb{(y^{(0)})^2}(\xi_1)\mid\Hc_1]=X^{(y^{(0)})^2}_m$ by the definition of $X^z_m$; and $\E[w_{m1}\mid\Hc_1]=g_m$. This gives $\E[G_1^2\mid\Hc_1]\le C\rho\sum_mx_m^2\tau_m$, with
\[
\tau_m=g_m(TZ^2)(\xi_m)\Big(1+\frac{CW_N^2}{N^2}\Big)+X^{(y^{(0)})^2}_m+g_mu_m^2+\frac{CW_N^4}{N^2}g_m\bar A^2 .
\]

\emph{Degree-biased averaging.} On $\OU$, $x_m\le|b_m|$ and $1\le2d_m/\dmin$, so $\rho\sum x^2\tau\le\frac{2\rho}{\dmin}\sum db^2\tau$, and (DB2) with $\Omega=\OU$ gives $\E[\oneU\rho\sum x^2\tau]\le C\frac{\rho n}{\dmin^2}\Eg{\tau;\OU}$. Each of the four parts of $\tau$ is bounded directly, without the indicator, so we use $\Eg{\tau;\OU}\le\Eg{\tau}=\E[\hat g_2\tau_2]$ and bound the four parts of $\E[\hat g_2\tau_2]$ as follows:
\begin{itemize}
\item $\E[\hat g_2g_2(TZ^2)(\xi_2)]\le C\int g\,TZ^2\,d\pi\le C(\E g^3)^{1/2}\E[gZ^4]^{1/2}$, by Cauchy--Schwarz and stationarity;
\item $C\E[gZ^4]^{1/2}(\mathfrak D+W_N^2/N)^{1/2}$, by Lemma~\ref{lem:A-twohop} with $\Omega$ the sure event and the second display of \eqref{eq:A-noiseU};
\item $\E[\hat g_2g_2u_2^2]\le\E[\avg_k(y^{(0)}_k)^2\hat g_k(\hat Pg)_k]\le C\E[gZ^4]^{1/2}(\E g^3)^{1/2}$, by Jensen, Lemma~\ref{lem:A-id}(a) and Cauchy--Schwarz;
\item $(\E[\hat g_2^2g_2^2])^{1/2}(\E\bar A^4)^{1/2}\le CK'$.
\end{itemize}
Multiplying by $N/n$ gives $C\rho NK'/\dmin^2\le CK'/\dmin$, since $\rho N\le C\dmin$.
\end{proof}

\begin{proposition}[$G_2$]\label{prop:A-G2}
$\frac Nn\E[G_2^2\oneU]\le CK'/\dmin$.
\end{proposition}

\begin{proof}
\emph{Pointwise bound.} By Lemma~\ref{lem:A-dual}(v) and the bound on $\hat y_m$ in the previous proof, $|G_2|\le\rho\sum_mw_{m1}c_m(|Y^W_1|+|u_m|+\vartheta)$. Separating the factors that depend on $\xi_1$ from those that depend on $m$, and using Lemma~\ref{lem:A-fb},
\begin{gather*}
G_2^2\le C\rho^2\Big[F(\xi_1)\Big(\sum_mw_{m1}c_m\Big)^2+\Big(\sum_mw_{m1}c_m|u_m|\Big)^2\Big],\\
F(v):=\Big(1+\frac{CW_N^2}{N^2}\Big)Z(v)^2+\nb{(y^{(0)})^2}(v)+\frac{CW_N^4}{N^2}\bar A^2 .
\end{gather*}

\emph{Averaging over $\xi_1$.} By Lemma~\ref{lem:A-dw}, applied to both terms with $x=c$, $\E[G_2^2\mid\Hc_1]\oneU\le C\frac{\rho}{n_0}\big(\sum dc^2\big)\mathcal Q\,\oneU$, with $\mathcal Q:=\int F\hat W+\sum_mu_m^2\kappa_m/\hat g_m$.

\emph{Conditioning on the types.} The key point is that $\mathcal Q$ is $\F$-measurable, as is $\Ohg$, while $\sum dc^2$ is not. Since $\OU\subseteq\Odeg^{(0)}\cap\Ohg$, Lemma~\ref{lem:A-energy} applies conditionally on $\F$:
\[
\E[G_2^2\oneU]\le C\frac{\rho}{n_0}\E\Big[\mathcal Q\one_{\Ohg}\,\E\Big[\one_{\Odeg^{(0)}}\sum dc^2\Bigm|\F\Big]\Big]\le C\frac{\rho n}{n_0\dmin^2}\E\mathcal Q .
\]
This step avoids the spectral event, and with it the logarithm.

\emph{The functional $\mathcal Q$.} By Lemma~\ref{lem:A-Wmom}(a)--(c), with $\Omega$ the sure event and $z=(y^{(0)})^2$ in (b), $\E\mathcal Q\le CNK'$. The term $\bar A^2\int\hat W$ is handled by Cauchy--Schwarz with $\E\bar A^4$ and $\E(\sum\kappa/\hat g)^2\le CN^2\mathsf M_N$. Multiplying by $N/n$ gives $C\rho NK'/\dmin^2\le CK'/\dmin$.
\end{proof}

\paragraph*{Decomposition of $G_3$}
The term $G_3$ is the feedback of node $1$'s edges through the network. Each neighbor
$m$ of node $1$ receives from node $1$ the difference $\tilde Y_m$ between node $1$'s
response and $m$'s neighborhood average, while the expected-weight graph delivers the
corresponding quantity $\hat y_m$. The discrepancy $\tilde Y_m-\hat y_m$ splits into
three parts, each with a simple meaning:
\begin{itemize}
\item $V_1$, the graph noise at node $1$ itself: the difference between node $1$'s
actual response and its expected-weight response;
\item $\gn_m$, the graph noise in $m$'s neighborhood average, in the graph without node
$1$;
\item $f_m$, the change caused by node $1$'s own presence: the echo of node $1$'s
influence, which returns to its neighbors through the network, together with the
feedback $\vartheta$ of Lemma~\ref{lem:A-fb}.
\end{itemize}
The noise at node $1$ is decomposed once more, by where it comes from. $V^\flat$ is
node $1$ averaging the graph noise of its neighbors. $V^\sharp$ is the randomness of
node $1$'s own edges acting on its neighbors' expected-weight responses. $\nu_1$ is
the echo and feedback at node $1$. Each part is a product of sums over node $1$'s
edges, and is bounded by expanding its moments (Lemma~\ref{lem:A-bern}).

Let $\echo:=Y_{\mathcal V}-Y^{(0)}$ be the \emph{echo}, the actual-graph counterpart of $\mathsf d=Y^W_{\mathcal V}-y^{(0)}$, and recall $\gn=PV^{(0)}+\cE y^{(0)}$. Writing $Y_{\mathcal V}=Y^{(0)}+\echo$ and $Y^W_{\mathcal V}=y^{(0)}+\mathsf d$ in the definitions of $\tilde Y_m$ and $\hat y_m$, and using $Y^{(0)}=y^{(0)}+V^{(0)}$,
\[
\tilde Y_m-\hat y_m=V_1-\gn_m-f_m,\qquad f_m=(P(\echo-\mathsf d))_m+(\cE\mathsf d)_m,\qquad|f_m|\le(P|\echo|)_m+3\vartheta .
\]
Hence $G_3=V_1\Lambda^A-G_{3\gn}-G_{3f}$, with $\Lambda^A=\sum_m\gamma_mA_{m1}\beta_m$ satisfying $|\Lambda^A|\le\Labs:=\sum_mA_{m1}x_m$, $G_{3\gn}=\sum\gamma_mA_{m1}\beta_m\gn_m$ and $G_{3f}=\sum\gamma_mA_{m1}\beta_mf_m$. On $\Odegone$, since $V_1=\gamma_1[(PY)_1-(\hat P^{\rm full}Y^W)_1]$,
\begin{gather*}
V_1=\gamma_1\big[V^\flat+V^\sharp+\nu_1\big],\qquad V^\flat=\frac1{d_1}\sum_mA_{m1}V^{(0)}_m,\\
V^\sharp=\frac1{d_1}\sum_mX_{m1}a_m,\qquad\nu_1=\frac1{d_1}\sum_mA_{m1}\echo_m-\hat P_1\mathsf d,
\end{gather*}
where $a_m=y^{(0)}_m-\zeta$, $\zeta=\hat P_1y^{(0)}$, $\sum_mp_{m1}a_m=0$ and $|\hat P_1\mathsf d|\le\vartheta$. For $\Hc_1$-measurable $z\ge0$, $(\hat P_1z)^2\le\hat P_1z^2\le\nb{z^2}(\xi_1)$, and averaging $w(\xi_m,\xi_1)\nb{z}(\xi_1)$ over $\xi_1$ gives $X^z_m$.

We apply Lemma~\ref{lem:A-bern} conditionally on $(\Hc_1,\xi_1)$, with $p_m=p_{m1}$, $s_{ij}=\sum_mp_{m1}x_m^iz_m^j$, $s_{01}/\mu_1=\hat P_1z$, $\mu_1\ge \dmin$ and $x\le C/\dmin$. Each resulting pattern is averaged over $\xi_1$ and bounded by a functional in (U2) or (U4), as recorded in Table~\ref{tab:A-patterns}.
\begin{table}[H]
\centering\small
\begin{tabular}{p{0.36\linewidth}p{0.33\linewidth}p{0.25\linewidth}}
\toprule
Pattern (after dividing by $\mu_1^2$) & Averaging and transfer & Functional\\
\midrule
$s_{10}^2\nb{z^2}$, $Z_1^2s_{10}^2$ & Lemma~\ref{lem:A-dw} & $\mathcal H$, $\mathsf M_N$ (Lemma~\ref{lem:A-Wmom}(a),(b))\\
$s_{20}\nb{z^2}$, $s_{10}\nb{z^2}$, $s_{11}^2$, $s_{12}s_{10}$, $t_{11}^2$, $t_{12}s_{10}$ & $X^{z^2}_m$; Lemma~\ref{lem:weight-transfer} & $\mathfrak D$ (Lemma~\ref{lem:A-twohop})\\
$s_{21}s_{01}$ & $\le\frac12(s_{22}\mu_1+s_{20}s_{02})$ & reduces to rows above and below\\
$s_{22}$, $t_{22}$ (after splitting $a^2$), $\sum p x^2 z$ & $\rho\sum gx^2z\le\frac4{n_0}\sum db^2z$; (DB2) & $\E[g|Z|^s]$ via \eqref{eq:A-noiseU}\\
$Z_1^2s_{20}$ & $\tau=g\,TZ^2$; (DB2) & $\E g^3$, $\E[gZ^4]$\\
$s_{10}\mu_1Z_1^2$ & $\tau^{Z^2}_m$; (DB1) & $\mathfrak G$ (Lemma~\ref{lem:A-deg1})\\
$s_{10}s_{01}$, $s_{11}$, $\sum pxZ_1^2$ & $X^z_m$ (off-diagonal), $\kappa_m$ (diagonal), $g(z+TZ^2)$; (DB1) & $\mathfrak D$, $\mathsf M_N$, $\E g^3$, $\E[gZ^4]$\\
$\vartheta$-terms & Lemma~\ref{lem:A-fb} & factor $W_N^4/N^2$; $\E\bar A^4$, $\mathsf M_N$\\
\bottomrule
\end{tabular}
\caption{Pattern bookkeeping for $G_3$ in setting (U). Every pattern retains a factor $s_{10}$ or $s_{20}$ (or $t_{\cdot}$), which carries the dual weight. Every functional is one of those in (U2) or (U4). The number of patterns is finite and independent of $N$. Throughout, (DB1), (DB2), Lemma~\ref{lem:A-twohop} and \eqref{eq:A-noiseU} are used with the indicator of $\OU$ in place, as the $z$-moments of $V^{(0)}$ and $\gn$ require.}
\label{tab:A-patterns}
\end{table}

\paragraph*{A worked example} The table compresses a computation that is the same
for every pattern; we carry out one instance in full. Take $G^\flat=\gamma_1V^\flat\Lambda^A$,
whose bound on $\Odegone$ is $(G^\flat)^2\le4\mu_1^{-2}(\sum_mA_{m1}x_m)^2(\sum_mA_{m1}z_m)^2$
with $z=|V^{(0)}|$, and consider the leading pattern $s_{10}^2s_{01}^2$ of the first
display of Lemma~\ref{lem:A-bern}, the one in which all four indices are distinct.
\begin{enumerate}
\item \emph{Rewrite in terms of node $1$'s averages.} Since $s_{01}=\sum_mp_{m1}z_m=\mu_1\hat P_1z$
and $\hat P_1$ is a probability vector, Jensen and $\hat g_1\ge c_0$ give
\[
4\mu_1^{-2}s_{10}^2s_{01}^2=4s_{10}^2(\hat P_1z)^2\le4s_{10}^2\,\hat P_1(z^2)\le4s_{10}^2\,\nb{z^2}(\xi_1).
\]
This is the entry ``$s_{10}^2\nb{z^2}$'' in the first row of Table~\ref{tab:A-patterns}.
\item \emph{Average over the type of node $1$.} Here $s_{10}=\rho\sum_mx_mw(\xi_m,\xi_1)$ and
$\nb{z^2}$ is an $\Hc_1$-measurable function of the type, so the first display of
Lemma~\ref{lem:A-dw}, with $F=\nb{z^2}$, gives
\[
\E\big[4s_{10}^2\,\nb{z^2}(\xi_1)\bigm|\Hc_1\big]\le\frac{8\rho}{n_0}\Big(\sum_md_mx_m^2\Big)\int\nb{z^2}\hat W .
\]
The averaging has turned node $1$'s edge weights into the degree-weighted profile $\hat W$.
\item \emph{The dual weight.} On $\OU\subseteq\Odeg^{(0)}$, Lemma~\ref{lem:A-dual}(ii) bounds
$\sum_md_mx_m^2\le\sum_md_mb_m^2\le Cn/\dmin$ deterministically. This is the factor that
carries the size $n$ of the sample and one power of $1/\dmin$; it is why every pattern
must retain a factor $s_{10}$ or $s_{20}$.
\item \emph{The functional.} By Lemma~\ref{lem:A-Wmom}(b) with $\Omega=\OU$ and $z^2$ in place of $z$,
\[
\E\Big[\oneU\int\nb{z^2}\hat W\Big]\le CN\,\E\big[\oneU\avg_k\hat g_kz_k^4\big]^{1/2}\big(\mathcal H+W_N^2\mathsf M_N/N\big)^{1/2}\le\frac{CNK'}{\dmin},
\]
where the last step is \eqref{eq:A-noiseU} with $s=4$, $\E[\oneU\avg_k\hat g_k|V^{(0)}_k|^4]\le C\dmin^{-2}\E[gZ^4]$, together with $W_N^2\le N$. This is where the smallness comes from: the fourth moment of the graph noise is of order $\dmin^{-2}$, and it is available only on $\OU$, which is why the indicator is carried along.
\item \emph{Collect and rescale.} Multiplying the three factors,
\[
\E\big[\oneU\one_{\Odegone}\,4\mu_1^{-2}s_{10}^2s_{01}^2\big]\le\frac{8\rho}{n_0}\cdot\frac{Cn}{\dmin}\cdot\frac{CNK'}{\dmin}\le\frac{C\rho nK'}{\dmin^2},
\]
and the Efron--Stein factor $N/n$ in \eqref{eq:ES} turns this into $C\rho NK'/\dmin^2\le CK'/\dmin$, since $\rho N\le2\dmin/c_0$.
\end{enumerate}
The other patterns differ only in step 1 (which coincidence occurs, hence which of
$s_{20}$, $s_{11}$, $s_{22}$ appears and whether $\nb{z^2}$ or $X^{z^2}_m$ results after
averaging) and in step 4 (which functional of (U4) bounds the result). Each coincidence costs a factor $\|x\|_\infty\le C/\dmin$ or $\mu_1^{-1}\le\dmin^{-1}$
relative to the leading pattern, provided the factor carrying the dual weight is
kept; where a crude bound would lose it, as for $t_{22}$ in $G^\sharp$, the pattern is
split before averaging (proof of Proposition~\ref{prop:A-G3a}). With this, none is
larger than the leading pattern.

The next proposition bounds the parts of $G_3$ that come from graph noise, as opposed
to the echo of node $1$'s own influence, which is treated in
Proposition~\ref{prop:A-echo}. There are three such parts: the noise at node $1$'s
neighbors ($G_{3\gn}$), node $1$ averaging that noise ($G^\flat$), and the randomness of
node $1$'s own edges acting on its neighbors' expected responses ($G^\sharp$). Each is
a product of two sums over node $1$'s edges, one carrying the dual weights and one
carrying the noise. The smallness comes from the noise. For $G_{3\gn}$ and $G^\flat$, it
enters through its second and fourth moments, which are of order $\dmin^{-1}$ and $\dmin^{-2}$
by \eqref{eq:A-noiseU}. For $G^\sharp$, the edges are centered, so their products
contribute only through coincidences, and each of these carries an extra factor $1/\dmin$.
In every case the result is at most $CK'/\dmin$ after the Efron--Stein scaling by $N/n$.
\begin{proposition}[$G_{3\gn}$, $G^\flat$, $G^\sharp$]\label{prop:A-G3a}
With $G^\flat=\gamma_1V^\flat\Lambda^A$ and $G^\sharp=\gamma_1V^\sharp\Lambda^A$,
\[
\frac Nn\E\big[(G_{3\gn}^2+(G^\flat)^2+(G^\sharp)^2)\oneU\one_{\Odegone}\big]\le\frac{CK'}{\dmin}.
\]
\end{proposition}

\begin{proof}
\emph{$G_{3\gn}$.} Given $(\Hc_1,\xi_1)$, $G_{3\gn}$ is linear in node $1$'s edges with $\Hc_1$-measurable coefficients $\gamma_m\beta_m\gn_m$, so by the last display of Lemma~\ref{lem:A-bern},
\[
\E[G_{3\gn}^2\mid\Hc_1,\xi_1]\le2\Big(\sum_mp_{m1}x_m|\gn_m|\Big)^2+2\sum_mp_{m1}x_m^2\gn_m^2 .
\]
For the first term, average over $\xi_1$ with the second display of Lemma~\ref{lem:A-dw} ($u=\gn$); this gives $\frac{C\rho}{n_0}(\sum dx^2)\sum_m\gn_m^2\kappa_m/\hat g_m$, and Lemma~\ref{lem:A-dual}(ii) and Lemma~\ref{lem:A-Wmom}(c) bound its expectation (with $\oneU$) by $C\rho nK'/\dmin^2$. For the second term, averaging over $\xi_1$ gives $\rho\sum gx^2\gn^2\le\frac4{n_0}\sum db^2\gn^2$ on $\OU$; then (DB2) with $\Omega=\OU$ and $\tau=\gn^2$, followed by \eqref{eq:A-noiseU} for $\Eg{\gn^2;\OU}=\E[\oneU\hat g_2\gn_2^2]$, bound its expectation by $C\frac{n}{N\dmin}\cdot\dmin^{-1}K'$.

\emph{$G^\flat$.} On $\Odegone$, $d_1\ge\mu_1/2$, so $(G^\flat)^2\le4\mu_1^{-2}(\sum Ax)^2(\sum Az)^2$ with $z=|V^{(0)}|$. Apply the first display of Lemma~\ref{lem:A-bern} and Table~\ref{tab:A-patterns}. The $z$-moments enter through $\E[\oneU\avg\hat gz^4]\le C\dmin^{-2}\E[gZ^4]$, from \eqref{eq:A-noiseU}; the leading pattern is carried out in full in the worked example after Table~\ref{tab:A-patterns}, and contributes at most $C\rho nK'/\dmin^2$ before multiplication by $N/n$.

\emph{$G^\sharp$.} On $\Odegone$, $(G^\sharp)^2\le4\mu_1^{-2}(\sum Ax)^2(\sum Xa)^2$. Apply
the second display of Lemma~\ref{lem:A-bern}. Since $a^2\le2(y^{(0)})^2+2\zeta^2$,
$\zeta^2\le\nb{(y^{(0)})^2}(\xi_1)$ and $\sum_mp_{m1}(y^{(0)}_m)^2\le\mu_1\nb{(y^{(0)})^2}(\xi_1)$, we have
$t_{02}\le4\mu_1\nb{(y^{(0)})^2}(\xi_1)$, hence $t_{02}/\mu_1^2\le\frac4{\dmin}\nb{(y^{(0)})^2}(\xi_1)$. With
$t_{11}^2\le s_{10}t_{12}$ and $t_{12}\le\|x\|_\infty t_{02}$, the patterns
$s_{10}^2t_{02}$, $s_{20}t_{02}$, $t_{11}^2$ and $t_{12}s_{10}$ are those of $G^\flat$
with $z^2$ replaced by $(y^{(0)})^2$ and an extra factor $1/\dmin$. The missing factor
$\dmin^{-2}$ in $\E[\avg\hat g(y^{(0)})^4]\le C\E[gZ^4]$ is thereby compensated. The pattern
$t_{22}$ must be split before averaging, so that the dual weight is not lost:
$t_{22}\le2\sum_mp_{m1}x_m^2(y^{(0)}_m)^2+2\zeta^2s_{20}$. The first part averages to
$\rho\sum gx^2(y^{(0)})^2\le\frac4{n_0}\sum db^2(y^{(0)})^2$ and is handled by (DB2).
The second is at most $\nb{(y^{(0)})^2}(\xi_1)s_{20}$, which averages to
$\rho\sum x^2X^{(y^{(0)})^2}$. (Bounding $t_{22}\le\|x\|^2_\infty t_{02}$ instead would
lose the dual weight and give only $CN/(n\dmin^3)$.)

In each case multiplication by $N/n$ gives at most $C\rho NK'/\dmin^2\le CK'/\dmin$.
\end{proof}

The remaining parts of $G_3$ come from the \emph{echo}: the influence of node $1$
returning to itself and to its neighbors through the network. When node $1$ is added
to the graph, each of its neighbors $m$ receives a new term in its neighborhood
average, so its response changes, by $\echo_m$. These changes propagate: they alter the
responses of other nodes, including other neighbors of node $1$, and node $1$'s own
response through its neighborhood average. The simplest instance is a backtracking
path $1\to m\to1$, in which node $1$'s input reaches $m$ and comes straight back; longer
cycles through node $1$ contribute in the same way.

The echo is self-referential, because the size of the change at the neighbors depends
on the differences $\tilde Y_m$ that node $1$ transmits, and these depend in turn on the
change. The proposition therefore first proves an a priori inequality, the size of
the echo is at most a constant times the size of the inputs plus a small multiple
of itself, and then rearranges it. The small multiple comes from a counting fact: a
path that leaves node $1$ and returns must pass through common neighbors of node $1$
and some other node $k$, and there are at most $d_k$ of these. This co-degree bound is
deterministic, holds for every graph, and makes the multiple of order $\dmin^{-2}$. No
condition on hubs is needed for this step.

\begin{proposition}[Echo terms]\label{prop:A-echo}
On $\OU\cap\Odegone$, let $\mathcal S:=d_1Z_1^2+\sum_mA_{m1}z_m$ with $z:=(Y^{(0)})^2+P(Y^{(0)})^2$, and $\Sigma:=\sum_mA_{m1}\tilde Y_m^2$. Then:
\begin{itemize}
\item[(a)] $\Sigma\le C\mathcal S$, for $N$ large enough that $\dmin^2\ge2C'$ with $C'=C'(\bar\gamma)$ below;
\item[(b)] $\nu_1^2\le\frac{C}{\dmin^2d_1}\mathcal S+2\vartheta^2$ and $G_{3f}^2\le\frac C{\dmin^3}\Labs\mathcal S+18\vartheta^2\Labs^2$;
\item[(c)] $\frac Nn\E\big[(\nu_1^2\Labs^2+G_{3f}^2)\oneU\one_{\Odegone}\big]\le CK'\big(\dmin^{-2}+W_N^4/N^2\big)$.
\end{itemize}
\end{proposition}

\begin{proof}
(a) \emph{The echo as a resolvent.} For $m\in\mathcal V$, $(P^{\rm full}Y)_m=(PY)_m+\varpi^A_m\tilde Y_m$ with $\varpi^A_m=A_{m1}/(d_m+A_{m1})$, by the computation in the proof of Lemma~\ref{lem:A-loo}. Subtracting the equations for $Y_{\mathcal V}$ and $Y^{(0)}$, as there, gives $\echo=R\Gamma(\varpi^A\tilde Y)$. By the $\ell^2(d)$ contraction (\ref{F:ops}),
\[
\sum_kd_k\echo_k^2\le C\sum_md_m(\varpi^A_m\tilde Y_m)^2\le C\sum_m\frac{A_{m1}\tilde Y_m^2}{d_m}\le\frac C{\dmin}\Sigma .
\]
\emph{The transmitted differences.} Since $Y=Y^{(0)}+\echo$ on $\mathcal V$, $\tilde Y_m=Y_1-(PY^{(0)})_m-(P\echo)_m$, and by $(a+b+c)^2\le3(a^2+b^2+c^2)$ and Jensen, $\tilde Y_m^2\le3[Y_1^2+(P(Y^{(0)})^2)_m+(P\echo^2)_m]$.

\emph{The co-degree bound.} Since $\sum_mA_{m1}A_{mk}\le d_k$ (node $k$ has at most $d_k$ neighbors in common with node $1$) and $d_m\ge\dmin/2$,
\[
\sum_mA_{m1}(P\echo^2)_m=\sum_k\echo_k^2\sum_m\frac{A_{m1}A_{mk}}{d_m}\le\frac2{\dmin}\sum_kd_k\echo_k^2\le\frac C{\dmin^2}\Sigma .
\]
\emph{Node $1$'s response.} Also $Y_1=Z_1+\gamma_1d_1^{-1}\sum_mA_{m1}Y_m$ with $Y_m=Y^{(0)}_m+\echo_m$, so by Jensen,
\[
d_1Y_1^2\le2d_1Z_1^2+4\sum_mA_{m1}(Y^{(0)}_m)^2+4\sum_m\echo_m^2\le2d_1Z_1^2+4\sum_mA_{m1}(Y^{(0)}_m)^2+\frac C{\dmin^2}\Sigma ,
\]
using $\sum_m\echo_m^2\le(2/\dmin)\sum_md_m\echo_m^2$.

\emph{The fixed point.} Summing the bound on $\tilde Y_m^2$ over the neighbors of node $1$ and inserting the last two displays gives $\Sigma\le C\mathcal S+\frac{C'}{\dmin^2}\Sigma$, with $C'$ depending only on $\bar\gamma$. Since $\Sigma<\infty$, (a) follows once $\dmin^2\ge2C'$. This is the only fixed-point step. It uses the deterministic co-degree bound and (M4), and does not use (U1).

(b) By Jensen, $(d_1^{-1}\sum_mA_{m1}\echo_m)^2\le d_1^{-1}\sum_m\echo_m^2\le\frac C{\dmin^2d_1}\Sigma$, and the term $\hat P_1\mathsf d$ is at most $\vartheta$; this gives the bound on $\nu_1^2$ by (a). For $G_{3f}$, the part with $(P|\echo|)_m$ satisfies, by Cauchy--Schwarz with weights $A_{m1}x_m$,
\[
\Big(\sum_mA_{m1}x_m(P|\echo|)_m\Big)^2\le\Labs\|x\|_\infty\sum_mA_{m1}(P\echo^2)_m\le\frac C{\dmin^3}\Labs\Sigma ,
\]
and the $\vartheta$-parts contribute at most $3\vartheta\Labs$. Then (a) applies.

(c) \emph{Reduction to patterns.} On $\Odegone$, $d_1^{-1}\mathcal S\le Z_1^2+\frac2{\mu_1}\sum A_{m1}z_m$. So the expectations to bound are $\E[\Labs^2\sum Az]$ and $\E[\Labs\mathcal S]$, which the third and fourth displays of Lemma~\ref{lem:A-bern} express through the patterns of Table~\ref{tab:A-patterns}. Here $\E[\oneU\avg\hat gz^2]\le C\E[gZ^4]$, by Lemma~\ref{lem:A-id}(c) and \eqref{eq:A-noiseU}.

\emph{The pattern $s_{10}s_{01}$.} In averaging $s_{10}s_{01}=\rho^2\sum_{m,k}x_mz_kw_{m1}w_{k1}$ over $\xi_1$, the terms
$k\ne m$ give $\rho^2(N-1)c_0\sum_mx_mX^z_m$, handled by (DB1) with $\Omega=\OU$ and Lemma~\ref{lem:A-twohop}. The diagonal $k=m$ gives
$\rho^2\sum_mx_mz_m\kappa_m\le C\rho^2\dmin^{-1}\sum_mz_m\kappa_m$, whose expectation is at
most $C\rho^2NK'/\dmin$ by Cauchy--Schwarz with $\E[\oneU\sum\hat gz^2]$ and $\E\sum\kappa^2/\hat g$.

\emph{Collecting.} The $\nu_1$-part carries the prefactor $\dmin^{-2}$, and the $G_{3f}$-part carries $\dmin^{-3}$ together with $\E[\Labs\mathcal S]\le C\rho^2Nn K'/\dmin$, which comes from the $s_{10}\mu_1Z_1^2$ pattern (Lemma~\ref{lem:A-deg1}). Both are at most $CK'/\dmin^2$ after multiplication by $N/n$.

\emph{The feedback terms.} Finally, $\vartheta^2\Labs^2\le\frac{CW_N^2}{N^2}(Z_1^2+W_N^2\bar A^2)\Labs^2$ by Lemma~\ref{lem:A-fb}. The $Z_1^2\Labs^2$ part is as above, with the extra factor $W_N^2/N^2$. For the $\bar A^2$ part, use $\E[\Labs^2\mid\Hc_1]\le\frac{2\rho}{n_0}(\sum dx^2)\int\hat W+\rho\sum gx^2$ and Cauchy--Schwarz with $\E\bar A^4$; this carries the factor $W_N^4/N^2$.
\end{proof}

\begin{proposition}[Setting (U)]\label{prop:A-flU}
Under (M1)--(M4) and (U1)--(U4), with $C_*$ large enough, for every $S$ with $|S|=n\le N$,
\[
n\Var(\bar V_S)\le CK'\Big(\frac1{\dmin}+\frac{W_N^4}{N^2}\Big)\to0 .
\]
\end{proposition}

\begin{proof}
By \eqref{eq:ES}, it suffices to bound $\E V_1^2$ and $\frac Nn\E G^2$. Corollary~\ref{cor:A-V1} gives $\E V_1^2\le CK'/\dmin$. For $G$, the bad event contributes at most $\frac Nn\cdot Cn^2N^{-4}\le CN^{-2}$. On the good event, $G=G_1+G_2+G_3$ and $G_3=V_1\Lambda^A-G_{3\gn}-G_{3f}$, with $|V_1\Lambda^A|\le(|V^\flat|+|V^\sharp|+|\nu_1|)\Labs$. Propositions~\ref{prop:A-G1}, \ref{prop:A-G2}, \ref{prop:A-G3a} and~\ref{prop:A-echo} bound each part.
\end{proof}

\begin{proof}[Proof of Theorem~\ref{thm:fluct}] Follows immediately from
Propositions~\ref{prop:A-flB} and~\ref{prop:A-flU}.
\end{proof}

\begin{remark}[Where (U1) enters the network half]\label{rem:A-U1net}
In setting (U), $W_N$ enters only through the feedback $\vartheta$ (Lemma~\ref{lem:A-fb}) and the coincidence corrections in Lemmas~\ref{lem:A-Wmom}(b), \ref{lem:A-twohop} and~\ref{lem:A-deg1}. Each of these carries a factor $W_N^2/N$ or $W_N^4/N^2$. The echo fixed point, Proposition~\ref{prop:A-echo}(a), does not use (U1).
\end{remark}

\section{Variance estimation}\label{app:svar}

\begin{lemma}[Variance perturbation]\label{lem:A-pert}
For vectors $a,b$ on $S$, write $\hVar_S(a)=\avg_S(a_i-\bar a_S)^2$. Then $|\hVar_S(a+b)-\hVar_S(a)|\le2\hVar_S(a)^{1/2}(\avg_Sb_i^2)^{1/2}+\avg_Sb_i^2$.
\end{lemma}

\begin{proof}
Expand $\hVar_S(a+b)=\hVar_S(a)+2\widehat{\Cov}_S(a,b)+\hVar_S(b)$. By Cauchy--Schwarz, $|\widehat{\Cov}_S(a,b)|\le\hVar_S(a)^{1/2}\hVar_S(b)^{1/2}$, and $\hVar_S(b)\le\avg_Sb_i^2$ because centering reduces the mean square.
\end{proof}

\begin{proof}[Proof of Proposition~\ref{prop:svar}]
\emph{Decomposition.} Write $Y_i=y^\circ(\xi_i)+e_i$ with $e_i=V_i+\Phi_i+R_i$, by \eqref{eq:split} and \eqref{eq:per-node}. The vector $e$ is permutation-equivariant, so by exchangeability and $(a+b+c)^2\le3(a^2+b^2+c^2)$,
\[
\E\avg_Se_i^2=\E e_1^2\le3(\E V_1^2+\E\Phi_1^2+\E R_1^2).
\]
We bound the three terms in turn.
\begin{itemize}
\item By Corollary~\ref{cor:A-V1}, $\E V_1^2\le C\E[gZ^2]/(N\rho_N)+CN^{-13}\E Z^2$. This argument uses neither a sup bound on $Z$ nor the spectral event.
\item Given $\xi_1$, $\Phi_1$ is an average of $N-1$ i.i.d.\ terms $\phi(\xi_1,\xi_l)$, which are centered by Lemma~\ref{lem:A-pop}(d). So $\E\Phi_1^2\le\|\phi\|^2_{L^2}/(N-1)\le C\Theta_2/N$.
\item By Lemma~\ref{lem:A-transfer}, $\E R_1^2\le C_Rr_N^2/N$.
\end{itemize}
Hence
\[
\E\avg_{i\in S}\big(Y_i-y^\circ(\xi_i)\big)^2\le C\Big(\frac{\E[gZ^2]}{N\rho_N}+\frac{\|\phi\|^2_{L^2}}{N}+\frac{r_N^2}{N}\Big)+CN^{-13}\E Z^2\to0 .
\]

\emph{The oracle term.} The variables $h_i:=y^\circ(\xi_i)$ are i.i.d.\ with $\sup_N\E|h_1|^{2+\delta}<\infty$. In (B) this holds because $|h|\le B_y$. In (U) it holds because $\|y^\circ\|_{L^4(dx)}\le c_0^{-1/4}\|y^\circ\|_{L^4(\pi)}\le C\E[gZ^4]^{1/4}$; this is a use of (M2) through the transfer $L^4(\pi)\to L^4(dx)$. A truncation argument gives the weak law for triangular arrays with uniformly bounded $(1+\eta)$-moments. Applied to $h$ and to $h^2$, it gives $\hVar_S(h)-\Var h_1\to0$ in probability.

\emph{Conclusion.} Lemma~\ref{lem:A-pert}, with $a=h$ and $b=e$, shows that $\hVar_S(Y)-\hVar_S(h)\to0$ in probability, and $s_S^2=\frac n{n-1}\hVar_S(Y)$ with $n/(n-1)\to1$.
\end{proof}

\begin{proof}[Proof of Corollary~\ref{cor:t}]
When $n/N\to0$, Lemma~\ref{lem:variance-bounds} gives $\sigma^2_{n,N}-\Var H_{\rm root}\to0$ and, under (M5), $\liminf\Var H_{\rm root}>0$. Proposition~\ref{prop:svar} then gives $s_S^2/\sigma^2_{n,N}\to1$ in probability. Combine with Theorem~\ref{thm:clt} and Slutsky's lemma.
\end{proof}

\section{Weighted means}\label{app:weighted}

Fix one input $Z$ and one weight $X=X(\xi)$; the vector case follows by Cram\'er--Wold, since every step below is linear in the pair $(X,Z)$ for a fixed component. Write $\bar U_S=\avg_SX_iY_i$.

\begin{lemma}[Signal part with unbounded weights]\label{lem:A-wsignal}
Assume the conditions of Theorem~\ref{thm:clt}, and that $\int X^2\Phi_2^2\le K$ and $\E|X(\xi)|^{2+\delta}\le K$, with $\E|X(\xi)|^{4(2+\delta)/(2-\delta)}\le K$ in setting (U). Then
\begin{gather*}
\avg_SX_iY^W_i-\E\avg_SX_iY^W_i=L^X_S+o_P(n^{-1/2}),\\
L^X_S:=\avg_{i\in S}\big(X_iy^\circ_i-\E X y^\circ\big)+\frac1N\sum_lH^X_{\rm rest}(\xi_l),
\end{gather*}
with $H^X_{\rm rest}(v)=\int X(x)\phi(x,v)\,dx$. Moreover $n\Var L^X_S=(1-n/N)\Var H^X_{\rm root}+(n/N)\Var H^X_{\rm full}$, and $L^X_S$ satisfies the Lyapunov condition whenever $\liminf n\Var L^X_S>0$. All of these conditions hold when $X$ is bounded.
\end{lemma}

\begin{proof}
\emph{Decomposition.} By \eqref{eq:per-node}, $X_iY^W_i=X_iy^\circ_i+X_i\Phi_i+X_iR_i$, and we treat the three terms in turn.

\emph{Remainder.} By Cauchy--Schwarz and Lemma~\ref{lem:A-transfer},
\[
\sqrt n\,\E\Big|\avg_{i\in S}X_iR_i\Big|\le\sqrt n\Big(\E\avg_{i\in S}X_i^2\Big)^{1/2}\Big(\E\avg_{i\in S}R_i^2\Big)^{1/2}\le C\|X\|_{L^2}r_N\to0 .
\]

\emph{Linear term.} The average $\avg_SX_i\Phi_i$ is the generalized $U$-statistic with kernel $\phi^X(x,v)=X(x)\phi(x,v)$. This kernel satisfies $\int\phi^X(x,v)\,dv=0$ and $\|\phi^X\|^2_{L^2}=\int X^2\Phi_2^2\le K$, so the proof of Proposition~\ref{prop:linear} applies verbatim, with $H^X_{\rm rest}$ in place of $H_{\rm rest}$.

\emph{Lyapunov.} The summands are bounded in $L^{2+\delta}$ as in Lemma~\ref{lem:A-lyap}. For the root term, $\|Xy^\circ\|_{L^{2+\delta}}\le C$: in (B) because $y^\circ$ is bounded, and in (U) by H\"older with $\|y^\circ\|_{L^4}$ and the moment of $X$ of order $4(2+\delta)/(2-\delta)$. For the rest term, Minkowski and H\"older give
\[
\|H^X_{\rm rest}\|_{L^{2+\delta}}\le\int|X|\Phi_{2+\delta}\,dx\le\|X\|_{L^{(2+\delta)/(1+\delta)}}\|\Phi_{2+\delta}\|_{L^{2+\delta}}\le C\|X\|_{L^{2+\delta}}\Theta_{2+\delta}^{1/(2+\delta)} .
\]
\end{proof}

\begin{lemma}[Network part, bounded weights]\label{lem:A-wnet}
If $|X|\le B_X$, then $n\Var(\avg_SX_iV_i)$ obeys the bounds of Propositions~\ref{prop:A-flB} and~\ref{prop:A-flU}, multiplied by $B_X^2$.
\end{lemma}

\begin{proof}
In Lemma~\ref{lem:A-loo}, replace $\one_{\tilde S}$ by $X\one_{\tilde S}$; for $i\ne1$, $X_i$ is $\Hc_1$-measurable, so the identity holds with the same proof. The dual vectors become $b=(I-P^{(0)}\Gamma)^{-1}\src_X$ with $\src_X=(D^{(0)})^{-1}X\one_{\tilde S}$, and $|\src_X|\le B_X|\src|$ entrywise. Every bound in Lemmas~\ref{lem:A-dual}, \ref{lem:A-energy} and~\ref{lem:weight-transfer} passes through $|b|\le\sum_t\bar\gamma^tP^t|\src_X|$, through $|\src_X-\hat\src_X|\le B_X|\src-\hat\src|$, or through $\|\Gamma\hat b\|$, each of which scales by $B_X$. In \eqref{eq:ES}, the change at node $j$ now includes $\delta_j(X_jV_j)$, with $\E(\delta_jX_jV_j)^2\le4B_X^2\E V_1^2$.
\end{proof}

\begin{proof}[Proof of Corollary~\ref{cor:delta}]
\emph{CLT.} Combine Lemmas~\ref{lem:A-wsignal} and~\ref{lem:A-wnet} with Cram\'er--Wold, and apply the delta method.

\emph{Variance estimation.} Write $X_iY_i=X_iy^\circ_i+X_ie_i$ and apply Lemma~\ref{lem:A-pert} coordinatewise, and its polarized version to the off-diagonal entries. Since $\E[X^2e^2]\le B_X^2\E e_1^2\to0$, the sample covariance converges to $\Cov\mathbf H_{\rm root}$. This equals $\Sigma_{n,N}+o(1)$ when $n/N\to0$, by Lemma~\ref{lem:variance-bounds}(b) applied to each linear combination.
\end{proof}

\section{The finite-population bias}\label{app:bias}

Until Proposition~\ref{prop:A-biasU}, we work in setting (B), with $\Delta_i:=\E[Y_i\mid\F]-Y^W_i$ and $\Ebar:=\E[\cE\mid\F]$, the conditional mean of the perturbation given the types. Let $\bar\pi_i=\mu_i/\sum\mu=\hat\pi_i$, and recall that $\Ods^{\rm full}$ is the intersection of the degree and spectral events of Lemma~\ref{lem:A-conc} for the full graph, so that $\Pr((\Ods^{\rm full})^c\mid\F)\le2N^{-14}$. Since $\hat g_l\le W_N\le1/(2\rho)\le N/2$ (using $\rho\ge1/N$), we have $\sum_l\hat g_l\le N^2/2$, which is used to bound the off-event contributions below.

The bias $\mu_N-\mu^\circ_N$ comes mainly from $\E\Delta_1$, where
$\Delta=\E[Y\mid\F]-Y^W$ is the mean effect of replacing edge probabilities by edges.
Taking the conditional expectation of $Y=R_PZ$ directly is not possible, because the
resolvent $R_P$ is a nonlinear function of the random graph. The expansion below
writes $R_P$ around the expected-weight resolvent $M$, which is fixed given the latent
positions, as a series in the perturbation $\cE=P-\hat P$. Each term can then be averaged
over the edges separately. We stop at second order for the following reason. The
first-order term involves only $\E[\cE\mid\F]$, which is small: replacing degrees by
their expectations in the normalization is unbiased up to a relative error of order
$1/N$ (Lemma~\ref{lem:A-B}). The leading bias therefore comes from the second-order
term, where two perturbations are multiplied. Their product has nonzero mean whenever
both involve the same edge, through backtracking or through the edge's contribution to
two degrees, and this is of order $1/(N\rho_N)$, which our bound captures up to a
logarithm. The third-order remainder is controlled crudely by operator norms on the
good event.

\begin{lemma}[Second-order resolvent expansion]\label{lem:A-resexp}
Everywhere, $Y-Y^W=M\Gamma\cE Y^W+M\Gamma\cE M\Gamma\cE Y^W+R_P\Gamma(\cE M\Gamma)^2\cE Y^W$. Hence $\Delta=T_1+T_2+T_3$ with
\[
T_1=M\Gamma\Ebar Y^W,\qquad T_2=M\Gamma\,\E[\cE M\Gamma\cE\mid\F]\,Y^W,\qquad T_3=\E\big[R_P\Gamma(\cE M\Gamma)^2\cE Y^W\bigm|\F\big].
\]
\end{lemma}

\begin{proof}
By Lemma~\ref{lem:A-ops}, $Y-Y^W=R_P\Gamma\cE Y^W$. The resolvent identity $R_P=M+R_P\Gamma\cE M$ follows from $R_P^{-1}-M^{-1}=-\Gamma(P-\hat P)=-\Gamma\cE$, multiplied by $R_P$ on the left and $M$ on the right. Applying it twice,
\[
R_P=M+M\Gamma\cE M+R_P\Gamma\cE M\Gamma\cE M .
\]
Multiplying by $\Gamma\cE Y^W$ on the right gives the expansion of $Y-Y^W$. Then take $\E[\cdot\mid\F]$, using that $M$, $\Gamma$ and $Y^W$ are $\F$-measurable.
\end{proof}

\begin{lemma}[Reciprocal moment of a Bernoulli sum]\label{lem:A-recip}
Let $J=\sum_mI_m$ with independent $I_m\sim\mathrm{Bern}(q_m)$, $q_m\le1/2$, $\bar q=\sum_mq_m\ge1$ and $q_2=\sum_mq_m^2$. Then $\E[1/(1+J)]=1/\bar q-\theta$ with $0\le\theta\le e^{-\bar q}/\bar q+2q_2/\bar q^3$.
\end{lemma}

\begin{proof}
\emph{Integral representation.} Since $1/(1+j)=\int_0^1s^j\,ds$, and $\E s^J=\prod_m(1-q_m+q_ms)$, substituting $s\mapsto1-s$ gives
\[
\E\Big[\frac1{1+J}\Big]=\int_0^1\prod_m(1-q_m(1-s))\,ds=\int_0^1\prod_m(1-q_ms)\,ds .
\]
\emph{Upper bound.} By $1-u\le e^{-u}$, the integrand is at most $e^{-\bar qs}$, and $\int_0^1e^{-\bar qs}ds\le1/\bar q$; so $\theta\ge0$.

\emph{Lower bound.} By $\log(1-u)\ge-u-u^2$ on $[0,1/2]$, the integrand is at least $e^{-\bar qs-q_2s^2}\ge e^{-\bar qs}(1-q_2s^2)$. Integrating, $\int_0^1e^{-\bar qs}ds=(1-e^{-\bar q})/\bar q$ and $\int_0^1s^2e^{-\bar qs}ds\le2/\bar q^3$, which gives the bound on $\theta$.
\end{proof}

\begin{lemma}[The first-order bias is of relative order $1/N$]\label{lem:A-B}
Under (M3)--(M4), for $i\ne j$,
\[
|\Ebar_{ij}|\le\hat P_{ij}\big(2\hat P_{ij}+16s_i+2\tilde e\big),\qquad s_i:=\sum_m\hat P_{im}^2,\qquad\tilde e:=e^{-\dmin/2},
\]
and $\sum_j|\Ebar_{ij}|\le18s_i+2\tilde e$.
\end{lemma}

\begin{proof}
\emph{Representation.} $P_{ij}=A_{ij}/d_i$, and on $\{A_{ij}=1\}$, $d_i=1+J^{(j)}_i$ with $J^{(j)}_i=\sum_{m\ne i,j}A_{im}$ independent of $A_{ij}$ given $\F$. So $\E[P_{ij}\mid\F]=p_{ij}\E[1/(1+J^{(j)}_i)\mid\F]$. The mean of $J^{(j)}_i$ is $\bar q_{ij}=\mu_i-p_{ij}\ge\mu_i/2\ge1$, by (M3) and (M4).

\emph{Applying Lemma~\ref{lem:A-recip}.} Since $\hat P_{ij}=p_{ij}/\mu_i$,
\[
\Ebar_{ij}=p_{ij}\Big(\frac1{\bar q_{ij}}-\theta_{ij}\Big)-\frac{p_{ij}}{\mu_i}=\hat P_{ij}\Big(\frac{\mu_i}{\bar q_{ij}}-1-\mu_i\theta_{ij}\Big)=\hat P_{ij}\Big(\frac{p_{ij}}{\bar q_{ij}}-\mu_i\theta_{ij}\Big).
\]
Here $p_{ij}/\bar q_{ij}\le2p_{ij}/\mu_i=2\hat P_{ij}$, and Lemma~\ref{lem:A-recip} gives $\mu_i\theta_{ij}\le2e^{-\mu_i/2}+16\sum_mp_{im}^2/\mu_i^2=2e^{-\mu_i/2}+16s_i$, using $\bar q_{ij}\ge\mu_i/2$. Summing over $j$ and using $\sum_j\hat P_{ij}=1$ gives the second claim.
\end{proof}

\begin{proof}[Proof of Proposition~\ref{prop:bias}, setting (B)]
\emph{Decomposition.} Write $\mu_N-\mu^\circ_N=\E\Delta_1+(\E Y^W_1-\mu^\circ_N)$. By Lemma~\ref{lem:A-transfer} with $S=[N]$, $|\E Y^W_1-\mu^\circ_N|\le Cr_N/\sqrt N$. For $\E\Delta_1$ we bound $\E\Delta_1^2=\E[N^{-1}\sum_i\Delta_i^2]$, by exchangeability. For any vector $v$, since $\hat g_i\ge c_0$,
\begin{equation}\label{eq:A-pernode}
\frac1N\sum_iv_i^2\le\frac1{c_0}\Big(\frac1N\sum_l\hat g_l\Big)\|v\|^2_{\ell^2(\bar\pi)}.
\end{equation}

\emph{$T_2$ and $T_3$.} On $\Ods^{\rm full}$, $\|\cE\|_{\ell^2(\bar\pi)}\le4\varepsilon$ by Lemma~\ref{lem:A-ops}, $\|M\Gamma\|_{\ell^2(\bar\pi)}\le\bar\gamma/(1-\bar\gamma)$, and $\|R_P\|_{\ell^2(\bar\pi)}\le3/(1-\bar\gamma)$. The last bound holds because $P$ is $d$-reversible and $d_i/\mu_i\in[\frac12,\frac32]$, so the norms of $\ell^2(\pi^d)$ and $\ell^2(\bar\pi)$ are equivalent with constant $\sqrt3$. Each of $T_2$ and $T_3$ contains at least two factors of $\cE$, so on $\Ods^{\rm full}$ its integrand has $\ell^2(\bar\pi)$ norm at most $C\varepsilon^2\|Y^W\|_{\ell^2(\bar\pi)}$. Off $\Ods^{\rm full}$, all operators are bounded in sup norm, and $\|Y^W\|_\infty\le B_y$. Since $\|\E[V\mid\F]\|_{\ell^2(\bar\pi)}\le\E[\|V\|_{\ell^2(\bar\pi)}\mid\F]$,
\[
\|T_2\|_{\ell^2(\bar\pi)}+\|T_3\|_{\ell^2(\bar\pi)}\le C\varepsilon^2\|Y^W\|_{\ell^2(\bar\pi)}+CN^{-14}B .
\]
By \eqref{eq:A-pernode}, $\sum\hat g(Y^W)^2\le(1-\bar\gamma)^{-2}\sum\hat gZ^2$ and $\sum_l\hat g_l\le N^2/2$, $\E[N^{-1}\sum_i(T_2+T_3)_i^2]\le CB^2(\varepsilon^4+N^{-27})$.

\emph{$T_1$.} By \eqref{eq:A-pernode} and the contraction of $M\Gamma$, $\E[N^{-1}\sum_i(T_1)_i^2]\le C\E[\hat g_1(\Ebar Y^W)_1^2]$. Split $Y^W=y^\circ(\xi_\cdot)+e$ with $e=\Phi+R$; the two parts behave differently, because $\Ebar$ is small on functions of the own type but only of order $\varepsilon$ in general.
\begin{itemize}
\item \emph{Own-type part.} By Lemma~\ref{lem:A-B} and $|y^\circ|\le B_y$, $|(\Ebar y^\circ)_1|\le B_y\sum_j|\Ebar_{1j}|\le B_y(18s_1+2\tilde e)$. Moreover $s_1=\sum_m\hat P_{1m}^2$ and $\hat P_{1m}=w_{1m}/((N-1)\hat g_1)$, so $\hat g_1s_1^2\le c_0^{-3}(N-1)^{-4}(\sum_mw_{1m}^2)^2$. Given $\xi_1=x$, separating the diagonal of the square,
\[
\E\Big[\Big(\sum_mw(x,\xi_m)^2\Big)^2\Bigm|x\Big]\le N^2\kappa_2(x)^2+N\int w(x,v)^4dv .
\]
Since $w^4\le w^2/(4\rho^2)$ by (M3), integrating over $x$ gives $\E[\hat g_1(\Ebar y^\circ)_1^2]\le CB^2\big(\mathsf M_N/N^2+\mathsf M_N^{1/2}/(N^3\rho^2)+\tilde e^2\big)$.
\item \emph{Global part.} Write $\Ebar=\E[\cE\one_{\Ods^{\rm full}}\mid\F]+\E[\cE\one_{(\Ods^{\rm full})^c}\mid\F]$. The first has $\ell^2(\bar\pi)$ norm at most $4\varepsilon$, and the second has sup norm at most $4N^{-14}$. With $\|e\|_\infty\le2B_y$ and $\sum_l\hat g_l\le N^2/2$,
\[
\E[\hat g_1(\Ebar e)_1^2]\le C\varepsilon^2\E[\hat g_1e_1^2]+CB^2N^{-27}.
\]
It remains to bound $\E[\hat g_1e_1^2]\le2\E[\hat g_1R_1^2]+2\E[\hat g_1\Phi_1^2]$. First, $\E[\hat g_1R_1^2]\le Cr_N^2/N$ by Lemma~\ref{lem:A-transfer}. Second, given $\xi_1=x$, $\hat g_1\Phi_1^2$ is a product of two averages over the other nodes, and $\Phi_1$ is centered given $x$, so only the index patterns with $l=l'$ survive in $\E[\hat g_1\Phi_1^2\mid x]$:
\[
\E[\hat g_1\Phi_1^2\mid x]\le\frac{g(x)\Phi_2(x)^2}{N-1}+\frac{\int w(x,v)\phi(x,v)^2dv}{(N-1)^2}.
\]
By Lemma~\ref{lem:A-rowcol}, $\phi^2\le CB^2\psi^2$ and $\|\psi\|_\infty\le C/\rho$, so $\iint w\phi^2\le CB^2\rho^{-1}\int gK_2\le CB^2\mathsf M_N^{1/2}/\rho$. Together with $\int\Phi_2^2d\pi\le C\Theta_2\le CB^2\mathsf M_N^{1/2}$, this gives $\E[\hat g_1\Phi_1^2]\le CB^2\mathsf M_N^{1/2}\big(N^{-1}+(N^2\rho)^{-1}\big)$.
\end{itemize}

\emph{Collecting.} Three elementary comparisons are used: $\varepsilon^2\ge1/N$, $(N\rho)^{-1}\le1$, and, for the remainder contribution, $\varepsilon^2r_N^2/N\le CB^2\mathsf M_N/(\rho N^2)$, which follows from Proposition~\ref{prop:A-remB} with $\varepsilon\le1/4$ and $N\rho\ge1$. (This last term is not dominated by $\mathsf M_N/N^2$, and is kept separately.) Also $\mathsf M_N^{1/2}/(N^3\rho^2)\le\mathsf M_N/(\rho N^2)$ by $\mathsf M_N\ge1$ and $N\rho\ge1$. Hence
\begin{gather*}
\E\Delta_1^2\le CB^2\Big(\varepsilon^4+\frac{\varepsilon^2\mathsf M_N^{1/2}}{N}+\frac{\mathsf M_N}{\rho N^2}\Big),\\
|\E\Delta_1|\le(\E\Delta_1^2)^{1/2}\le CB\Big(\varepsilon^2+\frac{\varepsilon\mathsf M_N^{1/4}}{\sqrt N}+\frac1N\Big(\frac{\mathsf M_N}{\rho}\Big)^{1/2}\Big)\le CB\Big(\varepsilon^2+\frac1N\Big(\frac{\mathsf M_N}{\rho}\Big)^{1/2}\Big),
\end{gather*}
using $\varepsilon\mathsf M_N^{1/4}/\sqrt N\le\varepsilon^2+\mathsf M_N^{1/2}/N\le\varepsilon^2+(\mathsf M_N/\rho)^{1/2}/N$ in the last step. By Proposition~\ref{prop:A-remB}, $r_N/\sqrt N\le CB(\mathsf M_N/\rho_N)^{1/2}/N$; each of its three terms is bounded this way, using $\mathsf M_N\ge1$ and $N\rho_N\ge1$. Since $\varepsilon^2\asymp\log N/(N\rho_N)$, we conclude that
\[
|\mu_N-\mu^\circ_N|\le CB\Big(\frac{\log N}{N\rho_N}+\frac1N\Big(\frac{\mathsf M_N}{\rho_N}\Big)^{1/2}\Big).
\]
The bound in setting (U) is Proposition~\ref{prop:A-biasU} below.
\end{proof}

\begin{proposition}[Finite-population bias, setting (U)]\label{prop:A-biasU}
Under (M1)--(M4) and (U1)--(U4),
\[
|\mu_N-\mu^\circ_N|\le C\Big(\frac{\log N}{N\rho_N}+\frac{W_N}N\Big).
\]
\end{proposition}

\begin{proof}
We follow the proof of Proposition~\ref{prop:bias} in setting (B), indicating the changes.
Lemmas~\ref{lem:A-resexp}--\ref{lem:A-B} use only (M3)--(M4).

\emph{$T_2$ and $T_3$.} The $\ell^2(\bar\pi)$ bounds are unchanged, with
$\E[N^{-1}\sum\hat gZ^2]=\E[gZ^2]$. On $(\Ods^{\rm full})^c$, sup norms are at most
$C\max_j|Z_j|$, and $\max_jZ_j^2\le\sum_jZ_j^2$. This gives
$\E[N^{-1}\sum(T_2+T_3)_i^2]\le C\varepsilon^4\E[gZ^2]+CN^{-26}\E Z^2$.

\emph{Own-type part of $T_1$.} By (U1), $\hat P_{1j}\le W_N/(c_0(N-1))$, so
$\sum_j\hat P_{1j}^2|y^\circ_j|\le\frac{CW_N}N(\hat P|y^\circ|)_1$ and $s_1\le CW_N/N$.
Lemma~\ref{lem:A-B} then gives
$|(\Ebar y^\circ)_1|\le\big(\frac{CW_N}N+2\tilde e\big)(\hat P|y^\circ|)_1$. Hence,
by Jensen and the $\hat g$-stationarity of $\hat P$,
\[
\E[\hat g_1(\Ebar y^\circ)_1^2]\le\Big(\frac{CW_N^2}{N^2}+C\tilde e^2\Big)
\E\big[\avg_i\hat g_i(y^\circ_i)^2\big]\le\frac{CW_N^2}{N^2}.
\]

\emph{Global part of $T_1$.} As before,
$\E[\hat g_1(\Ebar e)_1^2]\le C\varepsilon^2\E[\hat g_1e_1^2]+CN^{-26}\E Z^2$,
where now $\|e\|_\infty\le C\max_j|Z_j|$. By \eqref{eq:A-phiU} and $w\le W_N$,
$\int w(x,v)\phi(x,v)^2dv\le W_N\Phi_2(x)^2$, so
$\E[\hat g_1\Phi_1^2]\le C\Theta_2(N^{-1}+W_NN^{-2})$. Also
$\E[\hat g_1R_1^2]\le Cr_N^2/N$, with $r_N^2\le C(W_N^2/N+W_N^3/N^2)$ by
Proposition~\ref{prop:A-remU}.

\emph{Collecting.} We obtain $\E\Delta_1^2\le C(\varepsilon^4+\varepsilon^2/N+W_N^2/N^2)$,
so $|\E\Delta_1|\le C(\varepsilon^2+W_N/N)$, using $\varepsilon^2\ge1/N$. Finally
$r_N/\sqrt N\le C(W_N/N+W_N^{3/2}/N^{3/2})\le CW_N/N$.
\end{proof}

\begin{proof}[Proof of Corollary~\ref{cor:fixed-centering}]
In setting (B), Proposition~\ref{prop:bias} gives $\sqrt n|\mu_N-\mu^\circ_N|\le CB\big(\sqrt n\log N/(N\rho_N)+(\mathsf M_N/(N\rho_N))^{1/2}\big)\to0$, using $n\le N$ and (B3). In setting (U), Proposition~\ref{prop:A-biasU} gives $\sqrt n|\mu_N-\mu^\circ_N|\le C\big(\sqrt n\log N/(N\rho_N)+W_N/\sqrt N\big)\to0$, using (U1). Combine with Theorem~\ref{thm:clt}, or with Corollary~\ref{cor:t} when $n/N\to0$.
\end{proof}

\end{document}